%Authors:Dale Alspach

%Title:Tensor products and independent sums of $\Cal L_p$-spaces,
% $1<p<\infty$

%Filename:alspachprdscriptLp.atx
%TeX: AMSTeX
%Length:202005
%Received Date:4/8/96
%SubjectClass: 46B20
%Abstract:Two methods of constructing infinitely many isomorphically
%distinct $\Cal L_p$-spaces have been published. In this article we show
%that these constructions yield very different spaces and in the process
%develop
%methods for dealing with these spaces from the isomorphic viewpoint. We
%use these methods to give a complete isomorphic classification of the
%spaces $R_p^\alpha$ constructed by Bourgain, Rosenthal, and Schechtman and
%to show that $X_p\otimes X_p$ is not isomorphic to a complemented subspace
%of any $R_p^\alpha.$

%

%Citation:
%Use AMSTeX
%Last edited on 3/29/96
\input amstex
\documentstyle{amsppt}
\magnification=1200
%\leftheadtext{ \jobname \number\month--\number\day--\number\year}
%\rightheadtext{ \jobname \number\month--\number\day-- \number\year}

\define\tiny{\hskip.1em\relax}
\define\NxN{\Bbb N \times \Bbb N}
\define\N{\Bbb N}
\define\R{\Bbb R}
\define\Lp#1{L_p^{#1}}

\define\om{\omega}
\define\Om{\Omega}
\define\rp#1{R_p^{#1}}
\define\Rq#1{R_q^{#1}}

\topmatter
\title Tensor products and independent sums of $\Cal L_p$-spaces,
$1<p<\infty$
\endtitle
\rightheadtext{Tensor products and independent sums}
\leftheadtext{Dale Alspach}
\author Dale E. Alspach \endauthor
\affil Oklahoma State University \endaffil
\address Department of Mathematics,  Stillwater, OK 74078
\endaddress
\email alspach\@math.okstate.edu \endemail
\thanks Research supported in part by NSF grant DMS-9301506\endthanks
%  The following items provide publication information for the AMS-J logo
%Un-comment if using ams-j.sty
%\cvol{00}
%\cvolno{00}
%\cvolyear{0000}
%\cyear{0000}
%\cmonth{XXX}

\subjclass Primary 46B20 Secondary 46E30 \endsubjclass
\keywords projection, complemented, ordinal index \endkeywords
\abstract Two methods of constructing infinitely many isomorphically
distinct $\Cal L_p$-spaces have been published. In this article we show
that these constructions yield very different spaces and in the process develop
methods for dealing with these spaces from the isomorphic viewpoint. We
use these methods to give a complete isomorphic classification of the
spaces $R_p^\alpha$ constructed by Bourgain, Rosenthal, and Schechtman and
to show that $X_p\otimes X_p$ is not isomorphic to a complemented subspace
of any $R_p^\alpha.$
\endabstract

\toc
\head 0. Introduction \endhead
\head 1. The constructions of $\Cal L_p$-spaces \endhead
\head 2. Isomorphic properties of $(p,2)-$\tiny sums and the spaces $R_p^\alpha$
 \endhead
\head 3. The isomorphic classification of $R_p^\alpha$, $\alpha<\omega_1$
\endhead
\head 4. Isomorphisms from $X_p\otimes X_p$ into $(p,2)-$\tiny sums
\endhead
\head 5. Selection of bases in $X_p\otimes X_p$ \endhead
\head 6. $X_p\otimes X_p$-preserving operators on  $X_p\otimes X_p$ \endhead
\head 7. Isomorphisms of $X_p\otimes X_p$ onto complemented subspaces of 
$(p,2)-$\tiny sums \endhead
\head 8. $X_p\otimes X_p$ is not in the scale $R_p^\alpha,\alpha<\omega_1$
\endhead
\head 9. Final remarks and open problems \endhead
\specialhead  References \endspecialhead
%\hskip 7pt %
\endtoc
\endtopmatter
\document
\newpage
\head 0. Introduction \endhead

The purpose of this paper is to investigate the relationship
among three apparently
different constructions of $\Cal L_p$-spaces. We will show that two of the
methods 
produce the same isomorphic classes but that third method produces a
fundamentally different class of spaces. In particular the construction
due to Bourgain, Rosenthal and Schechtman \cite{BRS} will be shown to
produce different spaces than those Schechtman produced to show that
there are infinitely many isomorphically distinct $\Cal
L_p$-spaces. In order to explore the gap between the two
constructions we resurrect a 1974 construction of $\Cal
L_p$-spaces due to the author \cite{A1} that was presented in some seminars
at Ohio State but was not published at that time. (See \cite{F} for
a complete exposition and related results.) All of the methods of
construction make use of Rosenthal's fundamental space $X_p$ and
thus have a probabilistic or distributional character that makes
the passage to the isomorphic level difficult. One
consequence of this work is to show that modifications of the ideas
of Rosenthal can be used to work with these more complex spaces
within the isomorphic framework.

In Section 1.~we will review the constructions and the basic
properties. First we will describe some of the results in Rosenthal's paper.

\proclaim{Theorem 0.1} (Rosenthal's Inequality, \cite{R,Theorem 3} or
\cite{JSZ}) Let $2<p<\infty.$ Then there exists a 
constant $K_p$ depending only on $p$ such that if $f_1,\dots,f_n$ are 
independent, mean zero random variables in $L_p$, then
$$\align
\frac12 \max&\biggl \{\biggl (\sum_{i=1}^n \int \biggl |f_i\biggr |^p 
dx\biggr )^{1/p} , \biggl (\sum_{i=1}^n \int \biggl |f_i\biggr |^2
dx\biggr )^{1/2}\biggr\} \\
&\leq \biggl ( \int \biggl
|\sum_{i=1}^n f_i\biggr |^p dx
\biggr )^{1/p} \\
&\leq K_p\max\biggl \{\biggl (\sum_{i=1}^n \int \biggl |f_i\biggr |^p
dx\biggr )^{1/p} , \biggl (\sum_{i=1}^n \int \biggl |f_i\biggr |^2
dx\biggr)^{1/2}\biggr\}.
\endalign$$
\endproclaim

Using this inequality Rosenthal showed that there was a complemented subspace 
of $L_p$ which he called $X_p$ that was different (isomorphically) from the 
other complemented subspaces known at the time. In its sequential form 
$X_{p,(w_n)}$
is the completion of the space of sequences of real numbers $(a_n)$ with only 
finitely many $a_n$ non-zero under the norm
$$ \|(a_n)\|=\max \biggl \{\biggl (\sum_{i=1}^\infty |a_n|^p\biggr )^{1/p},
\biggl (\sum_{i=1}^\infty |a_n|^2 w_n^2\biggr )^{1/2}\biggr \}
$$
where $(w_n)$ is a bounded sequence of positive numbers such that for every
$\epsilon >0$, $$ \sum_{w_n<\epsilon} w_n^{2p/(p-2)}=\infty.\tag*$$

Using this space and a ``bad'' $l_2$-space formed by taking the definition of 
$X_{p,(w_n)}$ as above except that $w_n$ is a constant independent of $n$, 
Rosenthal defined a small list of additional $\Cal L_p$-spaces.
The spaces he defined were $B_p=(\sum_k X_{p,(w_{n}^k)})_p$ , where
$w_{n}^k=w_k$ for all $n$ and $\lim w_k =0,$ $B_p\oplus X_p,$
$X_p\oplus (\sum l_2)_p$, $X_p \oplus B_p$ and $(\sum X_p)_p.$ These
spaces, $L_p$, $l_p,$ $l_2 \oplus l_p$, and $X_p$ were the only known $\Cal
L_p$-spaces known at the time of Rosenthal's paper. Some of the isomorphic
relations between these spaces were determined by in Rosenthal. Others
were known to the author but never published. The current state of
knowledge of these smaller $\Cal L_p$-spaces can be found in the
dissertation  of G. Force \cite{F}.

One of the most
interesting parts of Rosenthal's paper is his proof that the space 
$X_{p,(w_n)}$  does not depend on the sequence $(w_n)$ as long as (*) is 
satisfied. We will make use of his ideas in this paper, so we record the basic 
formulas here.

\proclaim{Proposition 0.2} \cite{R,Lemma 7} Let $E_1, E_2, \dots$ be a sequence of 
disjoint finite subsets of $\Bbb N$. For each $j\in \Bbb N$ let
$$f_j=\biggl (\sum_{n\in E_j} w_n^{2p/(p-2)}\biggr)^{-1/p} \sum_{n\in E_j}
w_n^{2/(p-2)}e_n$$
where $(e_n)$ is the standard coordinate basis of $X_{p,(w_n)}$ and let
$$
f_j^*(\sum a_n e_n)= 
\biggl (\sum_{n\in E_j} w_n^{2p/(p-2)}\biggr)^{-(p-1)/p}
\sum_{n\in E_j} a_n w_n^{2(p-1)/(p-2)}.
$$
Then $(f_j)$ is 1-equivalent to the standard basis of $X_{p,(w_j')}$ where
$$w_j'=\biggl (\sum_{n\in E_j} w_n^{2p/(p-2)}\biggr)^{(p-2)/(2p)}$$ for each $j$ 
and $Px=\sum_j f_j^*(x) f_j$ is a contractive projection onto the closed linear
span of $(f_j).$
\endproclaim

Using Proposition 0.2 Rosenthal showed that every space $X_{p,(w_n)}$ 
 is isomorphic to a com\-ple\-men\-ted subspace of
$X_{p,(w_{n,k})}$, where $w_{n,k}\downarrow 0$ as $n\rightarrow \infty$ 
and $w_{n,k}$ does not depend on $k$, and that $X_{p,(w_{n,k})}$ is
complemented in every $X_{p,(w_n)}$ such 
that $(w_n)$ satisfies (*). Rosenthal then used a special sum of spaces and
a version of the Pelczynski decomposition method to show that all of the
spaces $X_{p,(w_n)}$ such 
that $(w_n)$ satisfies (*) are isomorphic.

In Section 2 we generalize Rosenthal's methods so that we may consider
spaces formed by replacing the scalars in the definition of $X_p$ by a
sequence of subspaces of $L_p.$ We prove a result analogous to Proposition
0.2 and develop the tools to use a decomposition method argument similar to
Rosenthal's. A crucial notion here is the use of a restricted class of
operators which are bounded in the $p$ and $2$ norms. This idea has been
used previously in working with $X_p.$ (See \cite{A2}, \cite{AC}
and \cite{JO}.) These
tools allow us to show that  for $\omega \leq \alpha <\omega_1$
the spaces $R_p^{\alpha+k}$, $k\in \N,$ are
isomorphic.

In order to distinguish the spaces $R_p^{\alpha}$ for $\alpha$ a limit
ordinal, we use some ideas of Schechtman \cite{S} but with a different basic
space. Whereas Schechtman used spaces of the form
$l_{r_1}\otimes l_{r_2} \otimes \dots
\otimes l_{r_k}$, we use the space $D_p$ originally defined in \cite{A1}.
Section 3 is devoted to completing the classification of the spaces
$R_p^\alpha.$

The remaining sections of the paper are devoted to the relationship between
$X_p\otimes X_p$ and $(p,2)-$\tiny sums of subspaces of $L_p.$ The goal is to
show that the two structures are esssentially incompatible and thus deduce
in Section 8 
that $X_p \otimes X_p$ is not isomorphic to a complemented subspace of any
$R_p^\alpha.$ In Section 4 we examine isomorphic embeddings of $X_p\otimes
X_p$ into $(p,2)-$\tiny sums of subspaces of $L_p$ and show 
that there are some restrictions on their behavior. Section 5 is devoted
to developing methods of passing to subsequences of the natural basis of
$X_p \otimes X_p$ which are bases for isomorphs of $X_p \otimes X_p$ so
that gliding hump style arguments can be used. The approach used is quite
general and may be of independent interest. We prove a general sufficient
condition for an
operator on $X_p\otimes X_p$ to be an isomorphism on a copy of $X_p\otimes
X_p$ in Section 6. We also look at the isomorphic types of certain natural
subspaces of $X_p\otimes X_p$ such as the span of the diagonal and lower
triangle basis vectors. With these methods in hand
we start looking at complemented embeddings of $X_p\otimes
X_p$ into $(p,2)-$\tiny sums of subspaces of $L_p$ in Section 7. The main
technical results are proved there.

In Section 9 we make some remarks about directions for further work and
list some open problems.

For the most part we use standard notation from Banach space theory as may
be found in the books of Lindenstrauss and Tzafriri, \cite{LT}, \cite{LTI},
\cite{LTII}.
We use the expression $a \sim b$ to denote equivalence of numerical
quantities up to multiplicative constants, i. e., there exist positive
numbers $K_1, K_2$ such that $K_1^{-1}a\leq b \leq K_2 a.$ We will assume
that the scalar field is the real numbers throughout. Unless otherwise
noted $p>2$ and $q<2$ is the conjugate index to $p$. If $F$ is a set, $|F|$
is its cardinality.

We would like to thank the Mathematical Sciences Research Institute at
Berkeley for its
support during which a portion of this work was completed.

\newpage

\head 1. The Constructions of $\Cal L_p$-spaces \endhead

In \cite{S} Schechtman used a simple tensor product to construct
infinitely many isomorphically distinct $\Cal L_p$-spaces. 
\allowlinebreak This tensor is
defined only for subspaces of $L_p.$

\definition{Definition 1.1}
Let $X$ and $Y$ be subspaces of $L_p[0,1]$ and define
$X \otimes Y$ to be the closed linear span of $\{x(t)y(s):x \in X
\text{ and }y\in Y\}$ in $L_p([0,1]\times[0,1])$ with the usual
product measure.
\enddefinition

This tensor product depends on the representation of $X$ and $Y$, a
priori. Note that if $T$ and $S$ are bounded operators on
$L_p[0,1]$ then we may define $T\otimes S$ on $L_p([0,1]\times[0,1])$
by $[T\otimes S](x\otimes y) = (Tx)\otimes (Sy)$ for all $x,y \in
L_p[0,1]$. A straight-forward calculation using the Fubini Theorem
shows that  $T\otimes S$ is well-defined and that 
$\|T\otimes S\|\leq \|T\|\cdot\|S\|$. It also follows from standard
techniques (integration against the Rademacher functions) that if
$(x_i)$ and $(y_i)$ are unconditional basic sequences in $ L_p[0,1]$
then $(x_i\otimes y_j)_{i,j}$ is an unconditional basis for
$[x_i]\otimes[y_j]$. (See \cite{S} or the proof of Lemma 1.2 below.)

Schechtman defines spaces $\otimes^n X_p =
X_p\otimes X_p \otimes \dots \otimes
X_p$ ($n$-times). The remarks above show that $\otimes^n X_p$ is
complemented in $L_p([0,1]^n)$ where the projection is $P\otimes
P \otimes \dots \otimes 
P$ ($n$-times) and $P$ is a projection from $L_p[0,1]$ onto $X_p$. Thus it is
immediate that $\otimes^n X_p$ is a $\Cal L_p$-space.
Unfortunately in this representation the norm of the projection
goes to $\infty$ with $n$. Indeed, it was communicated to me by
Schechtman from Pisier, that the norm of the projection must be at
least as large as the product of the smallest
norms of projections onto the
factors. To see this we have by \cite{TJ,Lemma 32.3}
that for $X \subset Y$ the relative projection constant
$$\multline
\lambda(\otimes^n
X,\otimes^n Y))=
\sup\{|\text{tr}(u:\otimes^n X
\rightarrow \otimes^n Y)|:u\in B(\otimes^n Y,\otimes^n Y), \\
\nu(u)\leq 1 \text{ and }u(\otimes^n X)\subset \otimes^n X\},
\endmultline$$
where tr denotes the trace and $\nu$ is the nuclear norm.
Because $\text{tr}(\otimes_1^n P)=\prod_1^n \text{tr}(P)$ and
$\nu(\otimes_1^n P)\leq \prod_1^n \nu(P)$, it follows that the
relative projection constant of $\otimes^n X_p$ in
$L_p([0,1]^n)$ is no better than the $n$th power of the relative
projection constant for $X_p$ in $L_p[0,1].$

Because of Rosenthal's Inequality it is possible to represent
$\otimes^n X_p$ as a sequence space relative to the
unconditional basis $(x_{k_1}\otimes x_{k_2} \otimes \dots \otimes
x_{k_n})_{k_j\in \Bbb N, 1\leq j \leq n}$ and explicitly compute a
formula for an equivalent sequence space norm. (In \cite{S} the case $n=2$ is
given.) The next lemma will
allow us to do the computation inductively. 

\proclaim{Lemma 1.2} Let $(x_n)$ be a normalized sequence of mean zero
independent random variables in $L_p[0,1]$, $2<p<\infty,$
and let $(y_k)$ be an
unconditional basic sequence in $L_p[0,1]$ and $Y=[y_k].$ Then for all
$(a_{n,k})$ in $\R,$
$$\multline
\|\sum_n \sum_k a_{n,k} x_n\otimes y_k \| \\
\sim 
\max \biggl\{\biggl( \sum_n 
\| \sum_k a_{n,k} y_k \|_p^p    
 \biggr )^{\frac1p},     
\biggl( \int \biggl\|\sum_k  [  \sum_n a_{n,k} \|x_n\|_2 r_n(\omega)]y_k
\biggr\|_Y^{p}d\omega\biggr )^{\frac1p}\biggr\},
\endmultline$$
where $r_n$ is the $n$th Rademacher function.
\endproclaim
\demo{Proof}
Let $f_n = \sum_k a_{n,k} y_k$ for all $n$. Then for each $t$,
$(x_n\cdot f_n(t))$ is a sequence of mean zero independent random
variables in $L_p[0,1]$ and thus by Rosenthal's Inequality,
$$\align 
\frac12 \max&\biggl \{\biggl (\sum_n \int \biggl
|x_n(s) f_n(t)\biggr |^p
ds\biggr )^{\frac1p} , \biggl (\sum_n \int \biggl |x_n(s) f_n(t)\biggr |^2
ds\biggr )^{\frac12}\biggr\} \\
&\leq \biggl ( \int \biggl
|\sum_n x_n(s) f_n(t)\biggr |^p ds
\biggr )^{\frac1p} \\
&\leq K_p\max\biggl \{\biggl (\sum_n \int \biggl
|x_n(s) f_n(t)\biggr |^p
ds\biggr )^{\frac1p} , \biggl (\sum_n \int \biggl |x_n(s) f_n(t)\biggr |^2
ds\biggr)^{\frac12}\biggr\}.
\endalign$$
Therefore
$$\align
&\biggl (  \iint \biggl
|\sum_n x_n(s) f_n(t)\biggr |^p ds\, dt
\biggr )^{\frac1p} \\
&\sim
\biggl (\max \biggl \{ \sum_n \iint \biggl 
|x_n(s) f_n(t)\biggr |^p
ds\, dt,
\int \biggl (\sum_n \int \biggl |x_n(s) f_n(t)\biggr |^2
ds\biggr)^{\frac p2}dt\biggr\}
\biggr )^{\frac1p}\\
&\sim 
\max \biggl \{\biggl ( \sum_n 
\|x_n\|_p^p\| f_n \|_p^p   
 \biggr )^{\frac1p},    
\biggl ( \int \biggl (\sum_n  \|x_n\|_2^2 | f_n(t) |^2 
\biggr)^{\frac p2}dt\biggr )^{\frac1p}\biggr\}\\
&\sim
\max \biggl \{\biggl ( \sum_n 
\| \sum_k a_{n,k} y_k \|_p^p    
 \biggr )^{\frac1p},     
\biggl ( \int \biggl (\sum_n  \|x_n\|_2^2 | \sum_k a_{n,k} y_k(t) |^2 
\biggr)^{\frac p2}dt\biggr )^{\frac1p}\biggr\}\\
&\sim 
\max \biggl \{\biggl ( \sum_n 
\| \sum_k a_{n,k} y_k \|_p^p    
 \biggr )^{\frac1p},     
\biggl ( \iint \biggl |\sum_n  \|x_n\|_2  \sum_k a_{n,k} y_k(t) r_n(\omega)
\biggr|^{p}d\omega\,dt\biggr )^{\frac1p}\biggr\}\\
\intertext{(by Khintchin's Inequality)}
&\sim
\max \biggl \{\biggl ( \sum_n 
\| \sum_k a_{n,k} y_k \|_p^p    
 \biggr )^{\frac1p},     
\biggl ( \int \biggl \|\sum_k \biggl [  \sum_n a_{n,k} \|x_n\|_2 r_n(\omega)
\biggr ]y_k
\biggr\|_Y^{p}d\omega\biggr )^{\frac1p}\biggr\}
\endalign$$
\qed\enddemo

Using this lemma we can now compute the sequence space norm of the
tensor product of finitely many copies of $X_p.$ Below we use the
convention that $\prod_{s\in \emptyset} f(s)=1.$

\proclaim{Proposition 1.3} For each $i,$ $1 \leq i \leq n$ let  $(x^i_k)_{k\in
\Bbb N}$ be a sequence of normalized
mean zero independent random variables in
$L_p[0,1]$ and let $w^i_k=\|x^i_k\|_2$ for all $i$ and $k.$
Then for all $(a(k_1,\dots,k_n))_{(k_1,\dots,k_n)\in \N^n}$ in $\R,$
$$ \multline
\biggl \| \sum_{k_1}\dots \sum_{k_n} a(k_1,\dots,k_n)
x^1_{k_1}\otimes \dots \otimes x^n_{k_n}\biggr \| 
\sim \max_{F} \biggl\{ \biggl ( \sum_{(k_s)_{s\in F}} \\
\biggl (
\sum_{(k_s)_{s\notin F}} a(k_{1},\dots,k_{n})^2 
W^2(F,k_{1},
\dots,k_{n})\biggr)^{\frac p2}\biggr)^{\frac1p}\biggr \},
\endmultline$$
where the maximum is taken over all 
subsets $F$ of the first $n$ natural numbers,
and $W^2(F,k_{1}, \dots,k_{n})= \prod_{s\notin F} (w^s_{k_{s}})^2.$
\endproclaim

\demo{Proof} The proof is by induction on the number of factors, $n$.
For $n=1$ the assertion is immediate from Rosenthal's Inequality.
Now assume it for $n$ factors and consider $n+1$ factors. By Lemma 1.2 with
$(y_k)$ replaced by $(x^1_{k_1}\otimes \dots \otimes x^n_{k_n})_{(k_1,\dots,
k_n)\in \Bbb N^n}$, we have that 
$$\align
&\biggl\|\sum_{k_{n+1}}\sum_{k_1}\dots \sum_{k_n} a(k_1,\dots,k_n,k_{n+1})
x^1_{k_1}\otimes \dots \otimes x^n_{k_n}\otimes x^{n+1}_{k_{n+1}}\biggr \|\\
&\sim \max \biggl \{\biggl ( \sum_{k_{n+1}}
\max_{F\subset \{1,2,\dots,n\}} \biggl\{ \biggl ( \sum_{(k_s)_{s\in F}}    \\
&\biggl (
\sum_{(k_{s})_{s\notin F} }
a(k_{1},\dots,k_{n},k_{n+1})^2
W^2(F,k_{1},
\dots,k_{n})\biggr)^{\frac p2}\biggr)^{\frac1p}\biggr \}^p
 \biggr )^{\frac1p},\\
&
\biggl ( \int \biggl \|\sum_{k_1,\dots,k_n}
 [  \sum_{k_{n+1}} a(k_1,\dots,k_n,k_{n+1})
 \|x^{n+1}_{k_{n+1}}\|_2
r_{k_{n+1}}
(\omega)]x^1_{k_1}\otimes \dots \otimes x^n_{k_n}
\biggr\|_p^{p}d\omega\biggr )^{\frac1p}\biggr\}
\endalign$$
Interchanging the summation over $k_{n+1}$ and the
{\sl max} in
the first expression produces the required expressions for which a
subset of $\{1,\dots,n+1\}$ would contain $n+1.$
Next we will use the following consequence of
Kahane's Inequality \cite{W,III.A.18} and the inductive hypothesis
to rewrite the second expression and
obtain the others. 

$$\align
\int |\sum_n [\sum_k a_{n,k} r_k(\omega)]^2|^{p/2}\,d\omega 
&= \int \|\sum_k [\sum_n a_{n,k}r_k(\omega)e_n]\|_{\ell_2}^p\,d\omega\\
&\sim (\int \|\sum_k [\sum_n
a_{n,k}r_k(\omega)e_n]\|_{\ell_2}^2\,d\omega)^{\frac{p}2}\\
&=(\sum_n \sum_k a_{n,k}^2)^{\frac{p}2},
\endalign$$ 
where $(e_n)$ is the usual unit vector basis of $\ell_2.$

$$\align
&\biggl ( \int \biggl \|\sum_{k_1,\dots,k_n}
 [  \sum_{k_{n+1}} a(k_1,\dots,k_n,k_{n+1})
 \|x^{n+1}_{k_{n+1}}\|_2  
r_{k_{n+1}} 
(\omega)]x^1_{k_1}\otimes \dots \otimes x^n_{k_n} 
\biggr\|_p^{p}d\omega\biggr )^{\frac1p}\\
&\sim
\biggl(\int\biggl (
\max_{F} \biggl\{ \biggl ( \sum_{(k_s)_{s\in F}}
\biggl (
\sum_{(k_{s})_{s\notin F} }
[\sum_{k_{n+1}}
a(k_{1},\dots,k_{n},k_{n+1})
\|x^{n+1}_{k_{n+1}}\|_2 r_{k_{n+1}} 
(\omega)]^2
\\
&
W^2(F,k_{1},
\dots,k_{n})\biggr)^{\frac p2}\biggr)^{\frac1p}\biggr \}\biggr
)^p d\omega\biggr )^{\frac1p}\\
&\aligned
\sim 
\max_{F}\biggl \{ \biggl (\sum_{(k_s)_{s\in F} } 
\int \biggl (\sum_{(k_s)_{s\notin F}}
&[\sum_{k_{n+1}} 
a(k_{1},\dots,k_{n},k_{n+1})
\|x^{n+1}_{k_{n+1}}\|_2 r_{k_{n+1}}  
(\omega)]^2 
\\ 
&
W^2(F,k_{1}, 
\dots,k_{n})\biggr)^{\frac p2} d\omega \biggr )^{\frac1p}\biggr \}
\endaligned\\
\intertext{(by interchanging the {\sl max} and the integral)}
&\sim \max_{F}\biggl \{ \biggl (\sum_{(k_s)_{s\in F}}  
\biggl (\sum_{(k_s)_{s\notin F}}\\ 
&\sum_{k_{n+1}}  
a(k_{1},\dots,k_{n},k_{n+1})^2
\|x^{n+1}_{k_{n+1}}\|^2_2
W^2(F,k_1,  
\dots,k_{n})\biggr)^{\frac p2}\biggr
)^{\frac1p}\biggr \}
\endalign$$
\qed\enddemo

Next we will briefly describe the construction of $\Cal L_p$-spaces
given in \cite{BRS}. Our exposition is slightly different, but the basic
ideas are the same.

\definition{Definition 1.4} Suppose that for each $n\in \Bbb N$,
$X_n$ is a subspace of $L_p(\Omega_n,\mu_n)$ for some
probability measure $\mu_n$.
Let $\Omega = \prod_{n=1}^\infty \Omega_n$ with the product measure
$ \mu=\prod \mu_n,$ and for each $n$ let $p_n$ be canonical map from
$\Omega$ onto $\Omega_n$ and $j_n(f)=f\circ p_n$ for all $f\in\Cal
L_p(\Omega_n,\mu_n)$. 
Let $X_0$ be the space of
constant functions on $\Omega$ and $j_0$ be the inclusion
of $X_0$ into $L_p(\Omega,\mu).$
Let $(\sum X_n)_I$ denote the closed linear span of $\cup_{n=0}^\infty j_n(X_n)$
and $(\sum' X_n)_I$ denote the closed linear
span of $\cup_{n=1}^\infty j_n(X_n)$ in $L_p(\Omega,\mu).$
 We will call $(\sum' X_n)_I$ the {\it
independent sum} of $\{X_n\}_{n=1}^\infty$ and $(\sum X_n)_I$ the {\it
complete independent
sum} of $\{X_n\}_{n=1}^\infty.$
\enddefinition

\remark{Remark 1.5} We will always choose the 
spaces $(X_n)_{n\geq 1}$ to be contained in the
mean zero functions. This will guarantee (See Lemma H.)
that the complete independent sum
has a natural unconditional decomposition into the spaces $(X_n)_{n\geq 0}.$
In \cite{BRS} the independent sum was used with spa\-ces, $(X_n),$
 containing the constants
and thus the independent summands were not necessarily direct summands.
\endremark \medskip

In the construction of $\Cal L_p$-spaces in \cite{BRS} a finite $\ell_p$ sum is
used. In this exposition we replace that approach by
using the tensor product.
For that purpose it is convenient to introduce the notation
$L_p^n$ for the space $L_p([0,1],\Cal D_n,\lambda)$ where
$\Cal D_n$ is the $\sigma$-algebra generated
by the intervals $I_k^n=[k2^{-n},(k+1)2^{-n})$ for
$k=0,1,\dots,2^n-1$ and $\lambda$ is Lebesgue measure. (Thus $L_p^n$ is
isometric to $\ell_p^{2^n}.$)

For each $\alpha < \omega_1$ we will define a subspace
$R_p^\alpha$ of $L_p(\mu)$ for some probability measure $\mu$.
The procedure is inductive. Let
$R_p^0=L_p^0$, the space of constant functions. Now suppose that
$R_p^\alpha$ has been defined. Define $R_p^{\alpha+1}=L_p^1\otimes
R_p^\alpha$. For a limit ordinal $\beta$ let $R_p^\beta =
(\sum_{n=1}^\infty
R_{p,0}^{\alpha_n})_I$ where $R_{p,0}^{\alpha_n}$ is the set of mean zero
functions in $R_p^{\alpha_n}$ and $(\alpha_n)_{n\in \N}$
 is an enumeration
of the ordinals $\gamma$, $0 <\gamma < \beta.$

\remark{Remark 1.6} Note that because $L_p^n\otimes L_p^m$ is isometric to
$L_p^{m+n}$, $R_p^{\alpha+k}$ is isometric to $L_p^k\otimes
R_p^\alpha,$ for any $\alpha<\omega_1.$ In the definition of
$R_p^\beta$ for $\beta$  a limit ordinal a more usual definition
would be to let $\alpha_n$ increase to $\beta$ and let
$R_p^\beta = (\sum
R_p^{\alpha_n})_I$. However no isomorphic theory of complete
independent sums was developed in \cite{BRS}, so at this stage we do
not know whether or not these spaces are all isomorphic.
\endremark \medskip

\newpage

\head 2. Isomorphic Properties of $(p,2)-$\tiny Sums and the Spaces $R_p^\alpha$
 \endhead

In \cite{BRS} it was shown by use of an ordinal index that uncountably
many of the spaces $R_p^\alpha$ are isomorphically distinct.
Unfortunately the nature of the proof is such that it does not provide any
additional information on which $\alpha's$ yield isomorphically different
spaces. One consequence of the results in this paper is that we
will see which ones are in fact different in a direct
fashion.

We will now use Rosenthal's inequality to get a little more
information about the spaces $R_p^\alpha$ for $\alpha$ a limit
ordinal. The approach used here is similar to that in \cite{A1} and \cite{F}.

\proclaim{Lemma 2.1} For each $n\in \N$ let $X_n$ be a subspace of
$L_p(\Om_n,\mu_n)$ which contains the constants and let $X_{n,0}$ denote the
the subspace of $X_n$ of all mean 0 functions in $X_n$. Then $(\sum X_n)_I$ is
isomorphic to $(\sum' X_{n,0})_I \oplus \Lp 0.$ Consequently,
$(\sum X_n)_I$ is isomorphic to the space 
$$\multline Z=\{(f_n)_{n=0}^\infty:f_n
\in X_{n,0} \text{ for } n \in \N, f_0\in \Lp 0,\\
\text{ and }
\|(f_n)\|=\max\{ (\sum \|f_n\|_p^p)^{1/p},(\sum
\|f_n\|_2^2)^{1/2}\}<\infty\}.
\endmultline$$
\endproclaim

\demo{Proof} For each $n\in \N$, $X_n=X_{n,0}\oplus[1_{\Om_n}].$
By definition $$(\sum X_n)_I=[\cup_{n=0}^\infty j_n X_n]
=[\cup_{n=1}^\infty j_n X_{n,0}\cup \{j_n(1_{\Om_n})
:n = 0,1,2,\dots\}].$$ Because $j_n(1_{\Om_n})=1_\Om$ for every $n$, it
follows that $(\sum X_n)_I=[\cup j_n(X_{n,0})]\oplus[1_\Om]=(\sum'
X_{n,0})_I\oplus \Lp 0.$ The final assertion is an immediate consequence of
Rosenthal's Inequality.
\qed\enddemo

We will have frequent use for spaces such as $Z$ above so we will
use the following notation for norm that occurs there.

\definition{Definition 2.2} Let $(X_n)$ be a sequence of subspaces of
$L_p(\Omega,\mu)$ for some probability measure $\mu,$ and let
$(w_n)$ be a sequence of real numbers, $0\leq w_n \leq 1.$ For any
sequence $(x_n)$ such that $x_n\in X_n$ for all $n$, let 
$$ \|(x_n)\|_{p,2,(w_n)}=\max\{ (\sum \|x_n\|_p^p)^{1/p},(\sum
\|x_n\|_2^2 w_n^2)^{1/2} \}
$$
and let 
$$X=(\sum X_n)_{p,2,(w_n)}=
\{(x_n):x_n \in X_n \text{ for all }n \text{ and
}\|(x_n)\|_{p,2,(w_n)}<\infty\}.$$
We will
say that $X$ is the $(p,2,(w_n))$-sum of $\{X_n\}$. In the special case
that $w_n=1$ for all $n$, we will sometimes write $(\sum X_n)_{p,2}$
instead of $(\sum X_n)_{p,2,(1)}.$
\enddefinition

In order to deal with these spaces with a norm defined in terms of an
$L_2$ norm and an $L_p$ norm
and with similar subspaces of $L_p$, we will use the terms $(p,2)-$\tiny bounded,
$(p,2)-$\tiny isomorphism, etc., to indicate that the map is bounded, is an
isomorphism, etc., in both norms. For example, if $T$ is a map from $(\sum
X_n)_{p,2,(w_n)}$ into $L_p$, we would say that $T$ is $(p,2)-$\tiny bounded, if
there exists a constant $K$ such that $\|T(x_n)\|_2\leq K (\sum
\|x_n\|_2^2 w_n^2)^{1/2}$ and $\|T(x_n)\|_p\leq K \|(x_n)\|_{p,2,(w_n)}$
for all $(x_n) \in (\sum X_n)_{p,2,(w_n)}.$ Note that this is slightly
different than the usage of this terminology in other papers, \cite{AC},
\cite{JO}.
This concept of a $(p,2)-$\tiny bounded operator is actually implicit in
\cite{R}
where it is used in estimating the norms of operators and in proving that
the spaces $X_{p,(w_n)}$, where $(w_n)$ satisfies (*) are all isomorphic.
He actually shows that the spaces are $(p,2)-$\tiny isomorphic. 

\remark{Remark 2.3} By Rosenthal's Inequality one can produce a
subspace of $L_p$ isomorphic to the $(p,2,(w_n))$-sum of $(X_n)$
in the following way. First symmetrize each space $X_n$ by mapping
$f\in X_n$ to $Sf \in L_p(\Omega',\nu)$ where $\Omega'=\Omega \times
\{0,1\}$ and
$\nu=\mu\times (\delta_{0}+\delta_{1})/2,$  by
$Sf(\omega,j)=f(\omega)(-1)^j.$ Next apply an isometry $J_n$ from
$L_p(\Omega',\nu)$ into $L_p(\Omega'\times[0,1],\nu\times
\lambda),$ where $\lambda$ is normalized Lebesgue measure on
$[0,1]$
to adjust the ratio between the $L_p$ and $L_2$ norms, where $J_n
f(\omega',\cdot)=f(\omega')1_{[0,\gamma_n]}\gamma_n^{-1/p}$ where
$\gamma_n=w_n^{2p/(p-2)}.$ For each $n\in\N$
let $\Omega_n=\Omega'\times [0,1]$ and
$\mu_n=\nu\times\lambda$ and identify the target space of $J_n$
with $L_p(\Omega_n,\nu_n).$ Finally let
$\Omega_\infty=\prod_{n=1}^\infty \Omega_n$,
$\mu_\infty=\prod_{n=1}^\infty \nu_n$ and $\pi_n$ denote the natural
projection from $\Omega_\infty$ onto $\Omega_n.$ Then for each
$n\in\N$, $T_n=(J_n S f)\circ \pi_n$
is an isometry from $X_n$ into $L_p(\Omega_\infty,\mu_\infty)$ such
that for any sequence $(f_n)$ with $f_n\in X_n$ for all $n$, $(T_n
f_n)$ is a sequence of independent mean zero random variables.
By Rosenthal's inequality $\|\sum T_n f_n\|_p$ is equivalent to 
$\|(T_n f_n)\|_{p,2,(1)}$ and by the definition of $T_n$ we have
that $\|(T_n f_n)\|_{p,2,(1)}=\|(f_n)\|_{p,2,(w_n)}$ because
$\|T_nf_n\|_2=\|f_n\|_2 w_n$ for all $n \in \N.$ Thus $[T_n X_n]$
is $(p,2)-$\tiny isomorphic to the $(p,2,(w_n))$-sum of $\{X_n\}.$
\endremark \medskip

 Let us now examine the spaces $\rp{\om k},$ for $k=1,2,\dots.$
First let $k=1$. For each $n\in \N$, $\Lp n=[1]\oplus[1_{I_k}-2^{-n}1_{I_0}:k
=1,2,\dots,2^n-1],$ where $I_k^n=[k2^{-n},(k+1)2^{-n}),$ for $k =
0,1,\dots,2^n-1,$ and $n\in \N.$
Because of Rosenthal's Inequality and Lemma H it is natural to compare
$\rp{\omega }$ with $(\sum \Lp n)_{p,2,(1)}$. 
Let $1_n=(x_k)\in (\sum \Lp n)_{p,2,(1)}$,
where $x_n=1_{[0,1]}$ and $x_k=0$ for $k\neq n.$
Observe that the operator $P(x_n)=((\int x_n)1_n)$  on $(\sum \Lp
n)_{p,2,(1)}$ is a (p,2)-norm
1 projection onto $[1_n:n\in \N]$ with kernel $(\sum \Lp
n{}_0)_{p,2,(1)}.$ Also we have that $[1_n:n\in \N]$ is $(p,2)-$\tiny isometric
to $X_{p,2,(1)}$ which is isomorphic (but not $(p,2)-$\tiny isomorphic) to $\ell_2.$

Let us explicitly compute the norm of an element $x\in (\sum \Lp
n)_{p,2,(1)}$. Let 
$$x=(\sum_{k=0}^{2^n-1} a_k^n 1_{I_k^n})_{n\in \N}.$$
Then
$$\|x\|_{p,2,(1)}=\max \{(\sum_{n=1}^\infty \sum_{k=0}^{2^n-1}
|a_k^n|^p2^{-n})^{1/p},(\sum_{n=1}^\infty \sum_{k=0}^{2^n-1}|a_k^n|^2
2^{-n})^{1/2}\}.$$
If we define $b_k^n=a_k^n2^{-n/p}$
and $w_k^n= 2^{-n(p-2)/(2p)}$ for
each $n,k,$ then
$$\|x\|=\max\{(\sum_{n=1}^\infty
\sum_{k=0}^{2^n-1}|b_k^n|^p)^{1/p},(\sum_{n=1}^\infty
\sum_{k=0}^{2^n-1}|b_k^n|^2 (w_k^n)^2)^{1/2}\}.$$
Thus $(\sum \Lp
n)_{p,2,(1)}$ is $(p,2)-$\tiny isometric to $X_{p,(w_k^n)},$ an isomorph of
$X_p.$ By Rosenthal's argument $X_p\oplus_{p,2}
X_{p,2,(1)}$ is $(p,2)-$\tiny isomorphic to
$X_p$, it follows that $(\sum L_p^n)_I$ is $(p,2)-$\tiny isomorphic to $(\sum \Lp
n)_{p,2,(1)}$.

Now we will prove that the spaces $\rp{\omega k}$ are 
$(p,2)-$\tiny iso\-mor\-phic to iterated $(p,2)-$\tiny sums of $X_p.$

\proclaim{Theorem 2.4} For any countable ordinal $\alpha $ and $j \in
\N$,
$\rp{\alpha +\omega j}$ is $(p,2)-$\tiny isomorphic to 
$$\underbrace{(\sum (\sum \dots(\sum \rp\alpha
)_{p,2,(w_{n,k})})\dots)_{p,2,(w_{n,k})})_{p,2,(w_{n,k})}}_{\text{$j$
times}},$$
where $w_{n,k}=2^{-n(p-2)/(2p)}$ for $0\leq k<2^n,n=1,2,\dots.$
\endproclaim

Before we begin the proof we need to make note of some properties
of $(p,2)-$\tiny iso\-mor\-phisms. Below, given a sequence of operators $(T_n)$,
$T_n:X_n\rightarrow Y_n,$ we will use the notation $\sum \oplus T_n$ denote
the operator defined from $(\sum X_n)$ to $(\sum Y_n)$ which is defined by
$T((x_n)_{n\in \N})=(T(x_n))_{n \in \N},$ for all finitely non-zero
sequences with $x_n\in X_n$ for all $n.$

\proclaim{Lemma 2.5} If $\{X_n\}$ and $\{Y_n\}$ are two sequences of
subspaces of $L_p(\Omega,\mu)$ for some probability measure $\mu,$
and there is a constant $K$ and $(p,2)-$\tiny continuous operators
$\{T_n\}$ such that $\|T_n\|_{p,2}\leq K$ for all $n\in \N$ then
the operator $T=\sum \oplus T_n$ is a $K$ $(p,2)-$\tiny bounded operator
from
$(\sum X_n)_{p,2,(w_n)}$ into $(\sum
Y_n)_{p,2,(w_n)}$. Consequently, 
\roster
\item if each $T_n$ is a
$(p,2)-$\tiny isomorphism and $\|T_n^{-1}\|_{p,2}\leq K$ as well then $T$
is a $(p,2)-$\tiny isomorphism.
\item if each $T_n$ is a projection, then $\sum \oplus T_n$ is a
$(p,2)-$\tiny bounded projection 
onto $$(\sum T_n X_n)_{p,2,(w_n)}.$$
\endroster
\endproclaim

The proof of Lemma 2.5 is straightforward and we leave it to the reader.
The next lemma is similar to Proposition 0.2 except that the scalars are
replaced by a subspace of $L_p$.

\proclaim{Lemma 2.6} Suppose that $\{X_n\}$ is a sequence of
subspaces of $L_p(\Omega,\mu)$ for some probability measure $\mu,$
$Y$ is a subspace of $L_p(\Omega,\mu)$, and $K$ is a constant
such that for every $n\in\N$ there is a projection $P_n$ from $X_n$
onto a subspace $Y_n$ such that $\|P_n\|_{p,2}\leq K$ and there is a
$(p,2)-$\tiny isomorphism $T_n$ from $Y$ onto
$Y_n$ with $\max\{\|T_n\|_{p,2},\|T_n^{-1}\|_{p,2}\}\leq K.$ Then, 
if $(w_n)$ is any sequence in $[0,1]$ with \newline
 $\sum w_n^{2p/(p-2)}<\infty$, 
$$Z=[(w_n^{2/(p-2)}T_n y):y\in Y]\subset (\sum X_n)_{p,2,(w_n)}$$
is $(p,2)-$\tiny isomorphic to $Y$ with
the norm
$$\max\{\|y\|_p,\|y\|_2 (\sum w_n^{2p/(p-2)})^{(p-2)/2p}\}$$
 and
$(p,2)-$\tiny complemented in $(\sum X_n)_{p,2,(w_n)}$. Moreover the
norms of the operators depend only on $p$ and $K$.
\endproclaim

\demo{Proof} First we will compute the norm of an element in $Z$.
$$\align \|(w_n^{\frac{2}{p-2}}T_n y)\|_{p,2,(w_n)}&=\max\{(\sum w_n^{\frac{2p}{p-2}}\|T_n
y\|_p^p)^{1/p},(\sum w_n^{\frac{4}{p-2}}\|T_n  
y\|_2^2 w_n^2)^{1/2}\} \\
&\sim \max\{(\sum w_n^{\frac{2p}{p-2}}\|
y\|_p^p)^{1/p},(\sum w_n^{\frac{2p}{p-2}}\|
y\|_2^2 )^{1/2}\}\\
&=(\sum w_n^{\frac{2p}{p-2}})^{1/p}\max\{\|y\|_p,
\|y\|_2(\sum w_n^{\frac{2p}{p-2}})^{\frac{p-2}{2p}}\}.
\endalign$$

Thus the map $(\sum w_n^{\frac{2p}{p-2}})^{1/p}y \rightarrow (w_n^{\frac{2}{p-2}}T_n
y)$ defines an
isomorphism from $Y$ in the weigh\-ted norm onto $Z$.
Next we define an operator $P$ by $$P(x_n)=
((\sum w_j^{\frac{2p}{p-2}})^{-1} w_n^{\frac{2}{p-2}}T_n[\sum w_j^{\frac{2(p-1)}{p-2}} T_j^{-1}
P_j x_j]),$$
for every sequence $(x_n)\in (\sum X_n)_{p,2,(w_n)}.$
Clearly $P$ maps into $Z$ and $Pz=z$ for all $z\in Z.$

It remains to check that $P$ is $(p,2)-$\tiny bounded.
$$\|P(x_n)\|\sim (\sum w_n^{\frac{2p}{p-2}})^{\frac{1-p}{p}}\max\{\|y\|_p, 
\|y\|_2(\sum w_n^{\frac{2p}{p-2}})^{\frac{p-2}{2p}}\}$$
where $y=\sum
w_j^{\frac{2(p-1)}{p-2}} T_j^{-1} P_j x_j.$ Since by H\"older's Inequality
with exponent  pairs $(p, p/(p-1))$ and $(2,2),$
$$\|y\|_p\leq\sum
w_j^{\frac{2(p-1)}{p-2}} \|T_j^{-1}\|_p\|P_j\|_p\|x_j\|_p\leq (\sum
w_j^{\frac{2p}{p-2}})^{\frac{p-1}{p}}K^2 (\sum \|x_j\|_p^p)^{1/p}$$
and 
$$\|y\|_2\leq\sum
w_j^{\frac{p}{p-2}}w_j \|T_j^{-1}\|_2\|P_j\|_2\|x_j\|_2\leq (\sum
w_j^{\frac{2p}{p-2}})^{1/2}K^2 (\sum \|x_j\|_2^2w_j^2)^{1/2}$$
Then $$\|P(x_n)\|\leq K^2 \max\{(\sum \|x_n\|_p^p)^{1/p},(\sum
\|x_n\|_2^2w_j^2)^{1/2}
.$$
\qed\enddemo

The next lemma shows that $(p,2)-$\tiny sums of subspaces of $L_p$ behave in a
reasonable way with respect to $(p,2)-$\tiny iso\-mor\-phisms of their summands.

\proclaim{Lemma 2.7} Suppose that for some $\alpha<\omega_1,$
$\{X_\beta:\beta<\alpha\}$ is a family of
subspaces of $L_p(\Omega,\mu)$ for some probability measure $\mu,$
and there is a constant $K$ such that 
for every $\gamma <\beta<\alpha $ there is a
 $(p,2)-$\tiny continuous projection 
$Q_{\gamma,\beta}$ from $X_{\beta}$ onto a subspace
$Y_\gamma$ such that $\|Q_{\gamma,\beta}\|_{p,2}\leq K$ 
and $Y_\gamma$ is $K$ $(p,2)-$\tiny isomorphic to $X_\gamma,$ in the natural
$(p,2)-$\tiny norm.
Let $(w_j)$
be a sequence contained in $(0,1]$ and let $\phi$ be a map from $\N$
into $\alpha $ such that there is an infinite subset $M$ of $\N$ with
$\lim_{n\in M}\phi(n)=\alpha$ and for every $\epsilon>0$,
$$\sum \Sb m\in M \\ w_m<\epsilon \endSb w_m^{2p/(p-2)}=\infty.$$
Then there is a constant $C$ depending on $K$ and $p$ such that 
$(\sum_{n=1}^\infty X_{\phi(n)})_{p,2,(w_n)}$ is
$C$ $(p,2)-$\tiny iso\-mor\-phic
to $(\sum_{n=1}^\infty X_{\alpha_n})_{p,2,(w'_n)},$ where $(w'_n)$
is any sequence in $(0,1]$ satisfying (*) and $\lim_{n\rightarrow
\infty}\alpha_n=\alpha $.
\endproclaim

\demo{Proof} We proceed in the same way that Rosenthal did in the
proof of \cite{R,Theorem 13}. First we single out a special $(p,2)-$\tiny sum.
Let $\psi:\N\times\N \rightarrow \alpha$ be a function such that for any
$\psi(n,j)$ is independent of $j$ and $(\psi(n,1))_{n\in\N}$ is an
enumeration of $\alpha.$
Let $Z=(\sum_{n,j}X_{\psi(n,j)})_{p,2,(w'_{n,j})}$, where
$w'_{n,j}=w'_n$ for all $n,j.$
Clearly if $(w'_n)$ satisfies (*) then we can find disjoint finite
subsets of $\N$, $N_{n,j}$, such that $\alpha_k>\psi(n,j)$ for all
$k\in
N_{n,j}$ and  $w'_j>
(\sum_{k\in N_{n,j}}
w_k^{\frac{2p}{p-2}})^{\frac{p-2}{2p}}
>w'_j/2$
for all $n,j\in \N.$  It follows from Lemma 2.6 
that $X_{\psi(n,j)}$ is $(p,2)$ isomorphic to a
$(p,2)-$\tiny complemented subspace of $(\sum_{n\in N_{n,j}}
X_{\alpha_n})_{p,2,(w'_n)}.$ An application of Lemma 2.5 gives us
that $Z$ is $(p,2)-$\tiny isomorphic to a $(p,2)-$\tiny complemented subspace
 of $(\sum_{n=1}^\infty X_{\alpha_n})_{p,2,(w'_n)}.$ Similarly, $Z$
is $(p,2)-$\tiny iso\-mor\-phic to a $(p,2)-$\tiny complemented subspace
 of $(\sum_{n=1}^\infty X_{\phi(n)})_{p,2,(w_n)}$. Another
application of this argument shows that $Z$ contains
$(p,2)-$\tiny complemented
$(p,2)-$\tiny isomorphs of $(\sum_{n=1}^\infty
X_{\alpha_n})_{p,2,(w'_n)}$ and $(\sum_{n=1}^\infty
X_{\phi(n)})_{p,2,(w_n)}.$

Finally we note that $(\sum_{n=1}^\infty Z)_{p,2,(1)} $ is $(p,2)-$\tiny isomorphic to
$Z$ and thus by the decomposition method (See \cite{R, Proposition 11}.) $Z$ is
$(p,2)-$\tiny isomorphic to \allowlinebreak $(\sum_{n=1}^\infty
X_{\alpha_n})_{p,2,(w'_n)}$ and $(\sum_{n=1}^\infty
X_{\phi(n)})_{p,2,(w_n)}.$
\qed\enddemo

\demo{Proof of Theorem 2.4} It is sufficient to prove this for $j=1$
and then we may conclude by applying Lemma 2.5 and induction.
We follow the pattern of the argument given above for the case
$\alpha
=1.$

It follows from Lemma 2.1 that $\rp{\alpha +\omega}$ is
$(p,2)-$\tiny isomorphic to 
$$(\sum \rp{\alpha + k}{}_0)_{p,2,(1)}\oplus_{p,2}
L_p^0.$$
Let us compare this with $(\sum \rp{\alpha + k})_{p,2,(1)}$. A typical
element of $\rp{\alpha+k}=\Lp k\otimes \rp{\alpha}
$ is of the form 
$$f=\sum_m (\sum_{j=0}^{2^k} a_j^k 1_{I_j^k})\otimes
g^k_m=\sum_{j=0}^{2^k} 1_{I_j^k}\otimes
(\sum_m a_{m,j}^k g^k_m)=\sum_{j=0}^{2^k}
1_{I_j^k}\otimes h_j^k,$$
 where $g^k_m \in \rp{\alpha} $ for all $k,m$
and $h_j^k=\sum_m a_{m,j}^k g^k_m.$ Thus a typical element of
$(\sum \rp{\alpha + k})_{p,2,(1)}$ is $f=(\sum_{j=0}^{2^k}
1_{I_j^k}\otimes h_j^k)_k$ and
$$\|f\|_{p,2,(1)}=\max \{(\sum_{k=1}^\infty \sum_{j=0}^{2^k-1}
\|h_j^k\|^p2^{-k})^{1/p},(\sum_{k=1}^\infty
\sum_{j=0}^{2^k-1}\|h_j^k\|^2
2^{-k})^{1/2}\}.$$
If we define $g_j^k=h_j^k2^{-k/p}$
and $w_j^k= 2^{-k(p-2)/(2p)}$ for
each $j,k,$ then
 $$\|f\|_{p,2}=
\max\{(\sum_{k=1}^\infty
\sum_{j=0}^{2^k-1}\|g_j^k\|^p)^{1/p},(\sum_{k=1}^\infty
\sum_{j=0}^{2^k-1}\|g_j^k\|^2 (w_j^k)^2)^{1/2}\}.$$
Thus $(\sum \rp{\alpha +k})_{p,2,(1)})$
is isometric to $(\sum \rp{\alpha })_{p,(w_k^n)}.$
Now just as in the $\alpha = 1$ case for each $n$ we
let $1_n=(x_k)_k$ where $x_n=1_{[0,1]}$ and $x_k=0$ for $k\neq n.$
Observe that the operator $P(f_n)=((\int f_n)1_n)$ is a $(p,2)-$\tiny norm
1 projection onto $[1_n:n\in \N]$ with kernel
$(\sum \rp{\alpha +k}{}_0)_{p,2,(1)}.$ $[1_n:n\in \N]$ is
$(p,2)-$\tiny isometric 
to $X_{p,2,(1)}$ which is isomorphic to $\ell_2.$ Moreover, $X_p$ is
$(p,2)-$\tiny isomorphic to a $(p,2)-$\tiny complemented subspace of $(\sum
\rp{\alpha +k}{}_0)_{p,2,(1)})$. Indeed, fix $f\in L_{p,0}^2$ and
observe that $L_p^k\otimes f\subset \rp{\alpha+k+2}$ for all $k.$ Thus
$(\sum
\rp{\alpha +k}{}_0)_{p,2,(1)}) \oplus_{p,2} X_{p,2,(1)}$ is $(p,2)-$\tiny
isomorphic to $(\sum
\rp{\alpha +k}{}_0)_{p,2,(1)})$. Hence $(\sum \rp{\alpha+k}{}_0)_{p,2,(1)}$
and $(\sum \rp{\alpha +k})_{p,2,(1)}$ are $(p,2)-$\tiny isomorphic.

It follows that $(\sum \rp{\alpha + k})_I$ is $(p,2)-$\tiny isomorphic to
$(\sum \rp{\alpha })_{p,(w_k^n)}.$
\qed\enddemo

\remark{Remark 2.8} By Lemma 2.7 we can replace
the weights $(w_{k,n})$ in Theorem
2.4 by any other sequence satisfying (*).
\endremark

\proclaim{Corollary 2.9} $\rp{\omega k}$ is $(p,2)-$\tiny isomorphic to the
space 
$$\{(a_{n_1,n_2,\dots,n_k})\in \R^{\N^k}: \|(a_{n_1,\dots,n_k})\|<\infty\}$$
where
 $$ \multline
\|(a_{n_1,\dots,n_k})\|
= \max_{0\leq j\leq k} \biggl\{ \biggl ( \sum_{n_1}\dots \sum_{n_j}
\\
\biggl (
\sum_{n_{j+1}}\dots \sum_{n_k}
a(n_1,\dots,n_k)^2
W^2(\{1,\dots,j\},n_{1},
\dots,n_k)\biggr)^{\frac p2}\biggr)^{\frac1p}\biggr \},
\endmultline$$
and $W^2(\{1,\dots,j\},n_{1},
\dots,n_{k})= \prod_{i=j+1}^k (w^i_{n_{i}})^2.$
\endproclaim

Notice that the norm given in Corollary 2.9 and that in Proposition 1.3 are
similar in form but there are fewer terms and a lack of symmetry in the
formula in Corollary 2.9. We will see much later that this difference is in
fact significant.

\proclaim{Corollary 2.10} For every $\alpha, \omega\leq \alpha < \omega_1$ and
$k\in \N$,
$\rp{\alpha+k}$ is $(p,2)-$\tiny isomorphic to $\rp{\alpha}$.
\endproclaim
\demo{Proof} By Theorem 2.4 $\rp{\alpha}$ is isomorphic to a $(p,2)-$\tiny sum of
spaces $X_\beta$ as in Lemma 2.7. From Lemma 2.7 it follows that $\rp{\alpha}$
is isomorphic to its square. Clearly $\rp{\alpha+k}$ is isomorphic to a
finite direct sum of copies of $\rp{\alpha}$ and thus to $\rp{\alpha}$.
\qed\enddemo

In the next section we will show that the converse of
Corollary 2.10 holds and thus it is actually a
characterization.
Our last goal in this section is to prove that for each $k\in \N$
the space $\rp{\omega k}$ is $(p,2)-$\tiny isomorphic to a
complemented subspace of $\otimes^k X_p$. This will be
accomplished by using the fact that $\rp{\omega
k}$ is $(p,2)-$\tiny isomor\-phic to the $k$-fold $(p,2,(w_n))$-sum of $X_p$,
where $(w_n)$ is any sequence in $[0,1]$ satisfying (*).

Proposition 2.11 is the analogue of \cite{R, Theorem 13}
and the consequences noted in
\cite{R, p 283}.

\proclaim{Proposition 2.11} If $X$ is a subspace of $L_p(\mu)$ for some
probability measure $\mu,$ and if $(w_n)$ is a sequence in $(0,1]$
there are four possible isomorphic
types for $Y=(\sum X)_{p,2,(w_n)}.$ The classes are determined as
follows.
\roster
\item If $\inf w_n >0,$ $Y$ is $(p,2)-$\tiny isomorphic to $(\sum
X)_{p,2,(1)}.$

\item If $\sum w_n^{2p/(p-2)}<\infty,$ $Y$ is isomorphic to $(\sum
X)_{p}.$

\item If $\sum w_n^{2p/(p-2)}=\infty$, $\inf w_n =0$, and there exists
$\epsilon>0$ such that 
$$\sum_{w_n<\epsilon} w_n^{2p/(p-2)}<\infty,$$
then $Y$ is isomorphic to $(\sum
X)_{p,2,(1)} \oplus (\sum X)_{p}.$

\item If $$\sum_{w_n<\epsilon} w_n^{2p/(p-2)}=\infty$$ for every
$\epsilon>0,$ then $Y$ is $(p,2)-$\tiny isomorphic to $(\sum
X)_{p,2,(w_n')},$ \linebreak
where $w_n'=n^{-(p-2)/(2p)}.$
\endroster
\endproclaim

\demo{Proof} Let $(x_n)$ be a sequence in $X.$ It is easy to see
that the third case follows from the first two and that the fourth
case is included in Lemma 2.7. Thus we need only compute the first two
cases.

If $\inf w_n >0,$ then 
$$\sum \|x_n\|_2^2 \geq \sum \|x_n\|_2^2 w_n^2 \geq (\inf w_n)\sum
\|x_n\|_2^2.$$

 If $\sum w_n^{2p/(p-2)}<\infty$ then 
$$\align \sum \|x_n\|_2^2 w_n^2 &\leq (\sum
\|x_n\|_2^p)^{2/p} (\sum w_n^{2p/(p-2)})^{(p-2)/p}\\
&\leq (\sum\|x_n\|_p^p)^{2/p} (\sum w_n^{2p/(p-2)})^{(p-2)/p}\endalign$$
\qed \enddemo

We will next show that $R_p^{\omega k}$ is isomorphic to a complemented
subspace of $\otimes^k X_p.$ The proof is inductive, so we first prove a
proposition which we can use iteratively.

\proclaim{Proposition 2.12} Suppose that $Y$ is a subspace of $L_p(\mu)$
for some probability measure $\mu$ and that $(w_n)$ and $(w_n')$
are sequences in $(0,1].$ For each $n,k\in \N$ let $w_{n,k}=w_n.$
 Then $X_{p,2,(w_n')}\otimes (\sum
Y)_{p,2,(w_{n,k})}$ contains a $(p,2)-$\tiny complemented subspace which is
$(p,2)-$\tiny isomorphic to 
$$(\sum (\sum Y)_{p,2,(w_n)})_{p,2,(w_n')}.$$
\endproclaim

\demo{Proof} For each $k\in \N$ let $Y_k=[(y_{n,j}):y_{n,j}=0
\text{ if }j\neq k]\subset (\sum
Y)_{p,2,(w_{n,k})}.$ Clearly $Y_k$ is $(p,2)-$\tiny isometric to $(\sum
Y)_{p,2,(w_n)}.$ Let $(x_n)$ be the standard basis of
$X_{p,2,(w_n')}$ and consider the space $Z=[x_k\otimes y_k:y_k\in
Y_k\text{ and }k\in \N].$ In order to compute the norm in the
tensor product we must first represent  $(\sum   
Y)_{p,2,(w_{n,k})}$ as a subspace of $L_p(\nu),$  for some
probability measure $\nu.$ We may assume that $X$ contains only
mean 0 functions and thus as in Remark 2.3 
we can consider $\nu$ as an infinite product measure and each
summand of $ (\sum    
Y)_{p,2,(w_{n,k})}$ as a space of functions depending only on the
$(n,k)$-coordinate of the product space. Similarly we realize
$(x_n)$ as a sequence of $L_p$-norm one
independent mean 0 random variables on
some other probability space. Thus a typical element of $Z$ is of the
form $\sum x_k\otimes y_k$ and by Rosenthal's inequality has
norm equivalent to 
$$\align &\hphantom{=}\max\{(\sum \|x_k\otimes y_k\|_p^p)^{1/p},
(\sum \|x_k\otimes y_k\|_2^2)^{1/2}\\
&= \max\{
(\sum \|x_k\|_p^p\| y_k\|_p^p)^{1/p},
(\sum \|x_k\|_2^2 \| y_k\|_2^2)^{1/2}\}\\
&=\max\{(\sum 
\|y_k\|_p^p)^{1/p}, 
(\sum (w_n')^2 \| y_k\|_2^2)^{1/2}\}.\endalign$$

This is precisely the norm in $(\sum Y_k)_{p,2,(w_n')}.$

To see that $Z$ is complemented it is sufficient to observe that
$$X_{p,2,(w_n')}\otimes (\sum   
Y)_{p,2,(w_{n,k})}$$ has an unconditional decomposition with
summands $[x_n]\otimes Y_k$, $n,k \in \N.$ (This follows
from Lemma 1.2.)
\qed\enddemo

\proclaim{Corollary 2.13} Suppose that $(w_n)$ is a sequence in
$(0,1]$ which satisfies (*). Then $$\otimes^k X_{p,2,(w_n)}$$ contains a
$(p,2)-$\tiny complemented subspace $(p,2)-$\tiny isomorphic to
$$\overbrace{(\sum \dots
(\sum \R)_{p,2,(w_n)}\dots )_{p,2,(w_n)}}^{\text{$k$ times}}.$$
Consequently, $\rp{\omega k}$ is $(p,2)-$\tiny isomorphic to a
$(p,2)-$\tiny complemented subspace of $\otimes^k X_p.$
\endproclaim

\demo{Proof} Because $(w_n)$ satisfies (*), $(\sum \R)_{p,2,(w_{n,k})}$ is
$(p,2)-$\tiny isomorphic to \newline
$(\sum \R)_{p,2,(w_n)}.$ Thus the $k=2$ case follows
from Proposition 2.12. Similarly, because Lemma 2.7 implies that 
$$\overbrace{(\sum \dots
(\sum \R)_{p,2,(w_n)}\dots )_{p,2,(w_n)}}^{\text{$k$ times}}$$
is $(p,2)-$\tiny isomorphic to 
$$(\sum\overbrace{(\sum \dots
(\sum \R)_{p,2,(w_n)}\dots )_{p,2,(w_n)}}^{\text{$k-1$ times}})_{p,2,(w_{n,k})
}.$$
the general case follows by induction.
\qed\enddemo

\newpage

\head 3. Isomorphic Classification of $R_p^\alpha$, $\alpha<\omega_1$ \endhead

To complete the isomorphic classification of the spaces $\rp{\alpha}$ we
will develop an invariant based on the presence of nicely placed copies of
$l_2$, $\alpha<\omega_1$.
In a certain gross sense the proof is similar to that of the
isomorphic classification of the spaces $C(\alpha)$, \cite{BP}. We begin with a
few definitions.

\definition{Definition 3.1} A sequence of measurable functions $(f_n)$ on
$(\Omega,\Cal B,\mu)$
will be said to be conditionally independent if there exists a measurable
set $A$ such that supp $(f_n)\subset A$ for all $n \in \N$ and the
restrictions to $A$ are independent for the restricted measure space and
normalized measure, $\mu|_A/\mu(A).$ The set $A$ is said to be a
conditioning set for $(f_n).$
\enddefinition

Conditionally independent sequences occur naturally in the spaces
$\rp{\alpha+k}$. For example if $(x_n)$ is an independent sequence in
$\rp{\alpha}$ then for any set $B$ which is measurable with respect to the
dyadic $\sigma$-algebra $\Cal D_k$, $(1_B\otimes x_n)$ is a conditionally
independent sequence in $\rp{\alpha+k}$.

We will be working with basic sequences in $L_q$, $1<q<2$,
which are equivalent to the basis
of $l_2$ which are well complemented in the sense that the orthogonal
projection is bounded. It is easy to see \cite{F, Lemma 4.4}
that if $[x_n]$
is orthogonally complemented in $L_q(\Omega)$ and $B$ is a measurable set of 
$\Omega_1$ then $[1_B\otimes
x_n]$ is orthogonally complemented in $L_q(\Omega_1\times \Omega).$

\definition{Definition 3.2} Suppose that $(x_n)$ is a sequence in $L_p$, $p>2$,
which is equivalent to the basis of $l_2$.
Let
$$W(x_n)=\lim_{n\rightarrow \infty}\sup
\{\frac{\|x\|_2}{\|x\|_p}:x\in [x_k:k\geq
n]\}.$$
If $(x_n)$
is a sequence in $L_q$, $1<q<2$, which is equivalent to the basis of $l_2$
and has complemented closed span with biorthogonal functionals $(x_n^*)$,
then let $W(x_n)= W(x_n^*).$ 
\enddefinition

This definition for $q<2$ depends on the choice of biorthogonal
functionals. In our application the projection will be the orthogonal
projection and thus the choice will be made automatically. The next lemma
illustrates the effect of conditional independence on $W().$

\proclaim{Lemma 3.3} If $(x_n)$ is a conditionally independent sequence of
mean 0 functions in $L_p$, $p>2$, which is equivalent to
the basis of $l_2$ and $A$ is a conditioning set for $(x_n)$, then
$$\mu(A)^{(p-2)/2p}/K_p\leq W(x_n)\leq\mu(A)^{(p-2)/2p}.$$
\endproclaim

\demo{Proof} 
By H\"older's inequality
$$\|x1_A\|_2^2\leq \|x^2\|_{p/2}\|1_A\|_{p/(p-2)}=\|x\|_p^2 \mu(A)^{(p-2)/p}.$$
This proves the right-hand inequality. For the left-hand we use Rosenthal's
inequality for the measure $\nu=\mu|_A/\mu(A).$
For convenience we assume that $\|x_n\|_{L_p(\nu)}=1.$
$$\align
\|\sum_{n=1}^N x_n\|_p\mu(A)^{-1/p}&=\|\sum_{n=1}^N x_n\|_{L_p(\nu)}\\
&\leq K_p \max\{(\sum_{n=1}^N
\|x_n\|_{L_p(\nu)}^p)^{1/p},(\sum_{n=1}^N\|x_n\|_{L_2(\nu)}^2)^{1/2}\}\\
&= K_p \max\{N^{1/p},(\sum_{n=1}^N\|x_n\|_2^2\mu(A)^{-1})^{1/2}\}.
\endalign$$
Because $(x_n)$ is equivalent to the basis of $l_2$, there is a constant
$c>0$ such that $\|x_n\|_2\geq c$ for all $n.$ Therefore 
$$(\sum_{n=1}^N\|x_n\|_2^2\mu(A)^{-1})^{1/2}\geq
N^{1/2}c/\mu(A)^{1/2}>N^{1/p}$$
for $N$ sufficiently large. For such $N$, we have
$$\|\sum_{n=1}^N x_n\|_p\mu(A)^{-1/p}\leq K_p
(\sum_{n=1}^N\|x_n\|_2^2\mu(A)^{-1})^{1/2}=\|\sum_{n=1}^N x_n \|_2
/\mu(A)^{1/2}.$$
Thus $$\mu(A)^{(p-2)/2p}/K_p\leq \|\sum_{n=1}^N x_n \|_2/\|\sum_{n=1}^N
x_n\|_p.$$
Because the argument applies $(x_{n+s})_{n=1}^N$ for any $s \in \N$,
$\mu(A)^{(p-2)/(2p)}/K_p\leq W(x_n).$
\qed\enddemo

We will also need to apply the function $W()$ to certain $(p,2)-$\tiny sums. These
spaces will in general be isomorphic to complemented subspaces of $L_p$ but
it is convenient to use the same notions in the given norm. In particular
in the next
lemma we will use an isomorph of $D_p$, \cite{F, p. 108}, $(\sum_{n,k}
l_2)_{p,2,(w_{n,k})}$ where $w_{n,k}=(2^{-n})^{(p-2)/2p)}$
for all $k,n$.
For $p>2$, the norm we will use is
$$\multline
\|\sum a^j_{n,k} e^j_{n,k}\|=\\
\max \{(\sum_{n,k} \sum_j
|a^j_{n,k}|^p)^{1/p},(\sum_{n,k} (\sum_j|a^j_{n,k}|^2)^{p/2})^{1/p},(\sum_{n,k}
\sum_j|a^j_{n,k}|^2w_{n,k}^2)^{1/2}\}
\endmultline$$
where $e^j_{n,k}$ denotes the sequence which is 0 except at the
$(j,n,k)$-th place where it is 1. For this space 
$$\|\sum a^j_{n,k} e^j_{n,k}\|_2=(\sum_{n,k}
\sum_j|a^j_{n,k}|^2w_{n,k}^2)^{1/2}$$ and we will compute $W(x_n)$ for
$(x_n)\subset D_{p,(w_{n,k})}$ from
the ratio $\|x\|_2/\|x\|$ for $x\in [x_n:n>k].$ For $q<2$, the norm on
$D_q$ is the dual norm to the norm on $D_p$, where $p$ is the conjugate
index. In the definition of $W()$ for $D_q$ we will take the biorthogonal
functionals from $D_p$ and proceed analogously. Notice that for fixed
 $n,k$,
$(e^j_{n,k}:j\in \N)\subset D_p, p>2$ is
1-equivalent to the basis of $l_2$ and has 1-complemented span and thus so are
the dual functionals. Further, for fixed $j$, $(e^j_{n,k}:n,k \in \N)$ is
1-equivalent to the basis of $X_{p,2,(w_{n,k})}.$

The next lemma is a key step in the proof and is analogous to Schechtman's
result, \cite{S, Proposition 2}, 
that the natural basis of $(\sum l_s)_{l_r}$, $1<q<r<s\leq 2$, is
not equivalent to a sequence of independent
random variables in $L_q.$ We
use several arguments which are given in that paper and leave the reader to
consult \cite{S} for additional detail.

\proclaim{Lemma 3.4} If $T$ is a bounded linear map from $(\sum
l_2)_{q,2,(w_{n,k})}$ into $L_q$, $1<q<2$, and the range of $T$ is
contained in $[x_n]$ where $(x_n)$ is a sequence of independent
mean zero random variables which are basic, then there is a normalized
block basic sequence
$(y^j_{n,k})$ such that $[y^j_{n,k}]$ is a $(q,2)-$\tiny complemented
subspace of $(\sum l_2)_{q,2,(w_{n,k})}$,  $(y^j_{n,k})$ is equivalent 
to the standard basis of $(\sum l_2)_{q,2,(w_{n,k})}$,
for each $n,k$, $W(y^j_{n,k})=w_{n,k}$, and $\|T(y^j_{n,k})\|<2^{-n-k}$
for all $j$.
\endproclaim

\demo{Proof} Let $(e^j_{n,k})$ be the natural basis of $(\sum
l_2)_{q,2,(w_{n,k})}$ where $W(e^j_{n,k}:j\in \N)=w_{n,k}.$ We can assume
that $\liminf_{k\rightarrow \infty} \liminf_{j\rightarrow
 \infty} \|Te^j_{n,k}\|>0$ for
each $n$, otherwise simple averaging of $(e_{n,j}^k)$ 
will produce the required sequence.
(See end of the proof below where we use a similar argument.) By a diagonal
argument and passing to a
subsequence of the index $j$ for each $n,k$ we may assume that 
$(T(e^j_{n,k}))$ is disjointly supported relative to the basic sequence
$(x_n).$ 

For each $n$, $(e^j_{n,k})_{k,j}$  is equivalent to the basis of $l_2$ and
thus $(Te^j_{n,k})_{k,j}$ is $p$-equi-integrable (See \cite{S} for the
definition.). By \cite{S, Lemma
4} and the proof of \cite{S, Proposition 2} we can find  a subsequence of the
$k$'s (which we assume to be the whole sequence for notational convenience)
and for each $k$ a subsequence of the $j$'s (which we again assume is the
whole sequence ) such that $(Te^j_{n,k})_j$ is
$(p,2^{-i})$-equi-distributed and $(Te^j_{n,k})_k$ is
$(p,2^{-t})$-equi-distributed. (The size of $i$ and $t$ are determined by
the estimate of symmetry required in the third display below.)
However for $n$ fixed and large and $t$ and $i$
sufficiently large the behavior in the domain is quite different. Indeed,
$$\|\sum_{j=1}^M e^j_{n,k}\| =  M^{1/2}$$
and
$$\|\sum_{k=1}^M e^j_{n,k}\| = M^{1/q}$$
if $M<w_n^{2q/(q-2)}.$ Thus
$$\|T\|M^{1/2}\geq \|T(\sum_{j=1}^M
e^j_{n,k})\|\sim \|T(\sum_{k=1}^M e^j_{n,k})\|.$$
Therefore 
$$\|T(M^{-1/q}\sum_{k=1}^M e^j_{n,k})\|\leq 2 \|T\|M^{(q-2)/2q}.$$
Also note that if $M$ is essentially equal to $(w_m/w_n)^{2q/(2-q)}$
then the 2
and $p$-norms of the corresponding dual functional are approximately
 equal to $w_m$.
Now select integers $M_{k,m}$ and $n_{k,m}$ such 
that $2\|T\|M_{k,m}^{(q-2)/2q}<2^{-k-m}$
and $M_{k,m}\sim (w_m/w_{n_{k,m}})^{2q/(2-q)}$.
(In particular, for a given $m$ we need to choose $n_{k,m}$ so that 
$(w_m/w_{n_{k,m}})^{2q/(2-q)}>2^{k+m}\|T\|.$ We also assume that
approximations improve as $k\rightarrow \infty.$)

For $k,m\in \N$ find disjoint sets $K_{k,m}$ of $\N$ of cardinality
$M_{k,m}$. (Technically we must actually choose appropriate subsequences of
the index set $j$ so that the required estimates above hold for $n=n_{k,m}$
and $M=M_{k,m}$. We also need to ensure that all of the blocks we construct
are disjointly supported. All of these are possible by simple
diagonalization but are messy notationally, so we leave the details to the
reader.)
Let
$y^j_{m,k}=M_{k,m}^{-1/q}\sum_{s\in K_{k,m}} e^j_{n_{k,m},s}$.
Then 
$$\|\sum_{s\in K_{k,m}}
{e^j_{n_{k,m},s}}^*\|_2/\|\sum_{s\in K_{k,m}}
{e^j_{n_{k,m},s}}^*\|_p\sim M_{k,m}^{1/2}w_{n_{k,m}}/M_{k,m}^{1/p}$$
which converges to $w_m$ as $k$
goes to $\infty.$ Thus $W(y^j_{m,k},j\in N)\sim w_m.$ Therefore
$[y^j_{m,k}:j\in N]$ is $(q,2)-$\tiny isomorphic to $X_{q,2,(w_m)}$ ($m$ is
fixed so the sequence $(w_m)$ is constant.) and by
(duality and)
Proposition 0.2 is $(q,2)-$\tiny complemented with a bound independent of $k$ and
$m.$
It follows from Lemma 2.5 that $[y^j_{k,m}]$ is
$(q,2)-$\tiny complemented.
\qed\enddemo

\remark{Remark 3.5} In the application of Lemma 3.4 and the next one we will
actually be working in some $\Rq{\alpha}$ rather than in sequence norm of 
$(\sum l_2)_{q,2,(w_{n,k})}$. Therefore we need to make some observations
about the $(q,2)-$\tiny norms that occur and the constant $W(y^j_{k,m})$. First
if for each $m,k$ we construct the sequence $(y^j_{k,m})$ in $L_q$ so that
the span is orthogonally complemented in $L_q$, then the estimates can be
made in $L_p=L_q^*$ with the dual functionals. Let us examine the dual situation
in $L_p$ more
closely.

In \cite{F, Theorem 4.8} it is shown that if for each $n\in \N$, 
$X_n$ is an orthogonally complemented subspace of
$L_p[0,1]$ (and mean 0)
and the norms of the projections do not depend on $n$, then the space
constructed by squeezing the supports onto sets of measure $w_n^{2p/(p-2)}$
by the obvious $L_p$-isometry, and then placing the resulting spaces on
independent coordinates produces a $(p,2)-$\tiny complemented
subspace of $L_p( \prod [0,1])$ which is isomorphic to $(\sum
X_n)_{p,2,(w_n)}$ and the image of $X_n$ in the resulting space is still
orthogonally complemented (with the same norm) and the only significant
change is that $l_2$-norms of the elements have been multiplied by a
constant, $w_n$. Second if the spaces $X_n$ are actually the span of
mean 0, independent random variables $(x_{n,k})$ with
$\|x_{n,k}\|_2/\|x_{n,k}\|_p=c,$ then the distance to $l_2$ and the norm of
the orthogonal projection onto $[x_{n,k}:k\in \N]$ is a function of $c$ and
the constants in Rosenthal's inequality.

Finally notice that in the proofs of Lemmas 3.4 and 3.6 only simple averages
 of elements whose dual functionals have the same ratio for $p$ and
$2$-norms are employed. Therefore the resulting dual functionals are also
simple averages. (Of course properly normalized.)
\endremark

We will actually need to consider a sequence of maps rather than just one
and so we will use diagonalization to extend Lemma 3.4 as follows.

\proclaim{Lemma 3.6} If $(T_r)$ is a sequence of bounded linear maps
from $(\sum
l_2)_{q,2,(w_{n,k})}$ into $L_q$, $1<q<2$, and the union of the ranges
of $T_r, r\in\N,$ is
contained in $[x_n]$ where $(x_n)$ is a sequence of independent mean zero
random variables which are basic, then there is a normalized
block basic sequence
$(y^j_{n,k})$ such that $[y^j_{n,k}]$ is a $(q,2)-$\tiny complemented
subspace of $(\sum l_2)_{q,2,(w_{n,k})}$,  $(y^j_{n,k})$ is equivalent
to the standard basis of $(\sum l_2)_{q,2,(w_{n,k})}$,
for each $n,k$, $W(y^j_{n,k})=w_{n,k}$, and $\|T_r(y^j_{n,k})\|<2^{-n-k}$
for all $r<n+k$ and all $j\in \N$.
\endproclaim

\demo{Sketch of Proof} For a given $m$  and $k$
only the finitely many operators $T_r$, such that $r<m+k$,
need be considered. In the proof of Lemma 3.4 
the choice of $M$ depends only on $m$ and $k$ and the norm of $T$. For a
finite number of operators the requirement that $T_r(e^j_{n,k})$ be a
 block subsequence  can be achieved simultaneously by passing to
subsequences and the estimates 
$$\|T_r\|M^{1/2}\geq \|T_r(\sum_{j=1}^M
e^j_{n,k})\|\sim \|T_r(\sum_{k=1}^M e^j_{n,k})\|.$$
can be achieved for all $r<m+k$ as well. Thus the proof above generalizes to
this case.
\qed \enddemo

\remark{Remark 3.7} The previous result could also have been obtained by
applying Lemma 3.4 iteratively and choosing $k$ so that $T_l(y^j_{n,k})$ is
small for
all $l<r$  and then using the fact that a sequence of independent random
variables has an upper $l_q$ estimate to see that
$$\|T_l(M^{-1/q}\sum_{k=1}^M y^j_{n,k})\|
\leq K(\sum_{k=1}^M \|T_l(y^j_{n,k})\|^q)^{1/q}/M^{1/q}\leq
K \max \|T_l(y^j_{n,k})\|.$$
\endremark

We are now ready for the main result of this section. As might be expected
the proof is by transfinite induction.

\proclaim{Theorem 3.8} Let $1<q<2.$ For every
$\alpha < \omega_1$, $R_q^{\alpha+\omega}$ is
not isomorphic to a subspace of $R_q^{\alpha}.$
\endproclaim

Theorem 3.8 is  actually a corollary of the following more informative
result.

\proclaim{Proposition 3.9} There are constants $C,C'$ such that if
$\omega_1>\alpha\geq \omega$ and for each $r\in \N$,
$T_r$ is a bounded
operator from $R_q^{\alpha}$ into $R_q^{\alpha'}$, $\alpha'+\omega\leq
\alpha$, then there is a
a normalized basic sequence $(y_n)$ in $R_q^{\alpha}$ such that
$(y_n)$ is equivalent to the standard basis of $l_2,$
$[y_n]$ is C-$(q,2)-$\tiny complemented, $W(y_n)>C'$, and $\|T_r(y_n)\|<2^{-n}$ for
all $r<n.$
\endproclaim
Before we begin the proof of this let us formalize a consequence of the
statement above that we will use in the induction.

\proclaim{Lemma 3.10} Let $X$ be a subspace of $L_q$ which is complemented by
the orthogonal projection and let $(x_n)$ be an unconditional basis for
$X$  with dual functionals $(x_n^*)\subset L_p$  
such that $(x^*_n)$ is an orthogonal sequence in $L_2$.
Suppose that $(X_j)$ is a sequence of subspaces of $L_q$
such
that for any $k$
if $(T_r)$ is a sequence of operators from $X$ into $(\sum_{j=1}^k
{}' X_j)_I$, then there is a a normalized basic sequence $(y_n)$ in
$X$ such that
$(y_n)$ is equivalent to the standard basis of $l_2,$
$[y_n]$ is C-$(q,2)-$\tiny complemented by the orthogonal projection,
 $W(y_n)>C'$, and $\|T_r(y_n)\|<2^{-n}$ for
all $r<n.$ Let $P_l$ denote the projection from $(\sum_{j=1}^\infty {}'
X_j)_I$ onto $(\sum_{j=1}^l{}' X_j)_I$.
Then for any $\epsilon>0$,
if $(T_r)$ is a sequence of operators from $X$ into
$(\sum_{j=1}^\infty X_j)_I$, then there is a a normalized basic sequence
$(y_n)$ in $X$ such that
$(y_n)$ is equivalent to the standard basis of $l_2,$
$[y_n]$ is $K(C+\epsilon)$-$(q,2)-$\tiny complemented, $W(y_n)>(C'-\epsilon)/K$,
where $K$  is the constant in the upper $l_2$-estimate for $(x^*_n)$,
and $\|P_l T_r(y_n)\|<2^{-n}$
for all $l,r<n.$
\endproclaim

\demo{Proof} The proof is a gliding hump argument.  Indeed, for each $l$ we
use the hypothesis for the sequence of operators
$(P_j T_r)_{r=1,}^\infty{}_{j\leq l}$
to get
a sequence $(y^l_n)$ as above. Let $y_1=y_1^1,$ and assume that $y_j$ has
been chosen for $j<m.$ For $n$ large $y^m_n$ satisfies $\|P_l
T_r(y^m_n)\|<2^{-m}$ for all $l,r<m.$ We know that $[y^m_n]$ is
complemented. Let $(y^{m*}_n)$ be the dual functionals in $X^*$.
By choosing $n$
large enough we can assume that the support of $y^{m*}_n$ relative to the
dual basis $(x_n^*)$ of $X^*$ is disjoint from the supports of $y^*_j$ for
$j<m.$ By construction $\|y^{m*}_n\|_2/\|y^{m*}_n\|_p>C'$,  for all $n,m$.
Therefore if we choose $y_m=y^m_{n(m)},$ so that $(y^{m*}_{n(m)})$ is an
orthogonal sequence then $(y_m)$ will be equivalent to the basis of $l_2$
and have orthogonally complemented closed span.
\qed\enddemo

\remark{Remark 3.11} The proof of Lemma 3.10 does not really require that $X$
have an unconditional basis, but that it has an unconditional F.D.D. so
that the copies of $l_2$ are taken from different summands
In this way selecting elements from  the different $l_2$ bases will still
produce a copy of $l_2$. Moreover in the application of Lemma 3.10 below, the
elements $(y_n)$ will actually be independent and mean 0, so that the
constant $K$ in the lemma can be replaced by 1.
\endremark

\remark{Remark 3.12} In \cite{BRS} it is shown that there is a certain isomorph
of $L_p$, $X^p_{\Cal D},$ so that for every $\alpha<\omega_1$,
$\rp{\alpha}$ is (isometric to a)
contractively complemented in $X^p_{\Cal D}.$ Also the
projection is the orthogonal projection. Moreover, if
$\alpha<\beta<\omega_1$, $\rp{\alpha}$ is contractively complemented by the
orthogonal projection in $\rp{\beta}$.
\endremark

\demo{Proof of Proposition 3.9} We
 will use induction on $\alpha$. Given $\epsilon>0$, 
the constants $C,C'$  can be chosen to be $(1+\epsilon) C_p$ and
$1-\epsilon$, respectively, where $C_p$ is the $(q,2)-$\tiny norm of the
orthogonal projection onto the span of the Rademachers in $L_p.$
At each stage of the
induction we actually will produced a badly preserved copy of
$(\sum l_2)_{p,2,(w_{n,k})}.$

If $\alpha=\omega$, $\rp{\alpha'}$ is finite
dimensional and the result is obvious. Note also that we may assume that
$\alpha$ is a limit ordinal since for $n\in \N$,
$\Rq{\alpha+n}$ contains a contractively
complemented subspace isometric to $\Rq{\alpha}$ and
$\alpha'+\omega\leq\alpha+n$ implies that $\alpha'+\omega\leq\alpha.$
Similarly we may asssume that $\alpha=\alpha'+\omega,$ because if
$\alpha>\alpha'+\omega$ then $\Rq{\alpha}$ contains a
contractively complemented subspace isometric to $\Rq{\alpha'+\omega}$ and
the inductive hypothesis for $\alpha'+\omega$
gives the result.

Suppose $\alpha=\omega\cdot 2.$
Because $\Rq{\omega}$ is isomorphic to $X_q$, we can assume that the the
union of the ranges
of $(T_r)$ is in the span of a sequence of mean zero independent random
variables.
$\Rq{\omega\cdot 2}=(\sum \Rq{\omega+n})_I$ contains  a well
$(q,2)-$\tiny complemented
($(q,2)-$\tiny isomorphic) copy of
$(\sum l_2)_{q,2,(w_{n,k})}$ since $\Rq{\omega}$ contains well
$(q,2)-$\tiny complemented copies of $X_{q,2,\delta},$ for any $0<\delta\leq 1.$
In particular, in order to avoid any problems with constants, for each
$n,k$ we can take a
sequence of independent mean zero Bernoulli random variables and find
an (almost) isometric  copy of their span in $\Rq{\omega+m(n,k)}$ with the
weight $W(y_{n,k})=w_{n,k}.$
Applying Lemma 3.6 to this copy of $(\sum l_2)_{q,2,(w_{n,k})}$ yields the
required sequence $(y_n)$. Note that the constants which occur here do not
depend on passing through the isomorphisms, because as noted in Remark 3.5
we actually get a
sequence $(y_n)$ of independent mean 0 random variables which have closed
span
complemented by the orthogonal projection from $L_q.$ The $(q,2)-$\tiny norm of
this projection is the same as the $(q,2)-$\tiny norm of the projection onto a
similar block in $X_q$ and the isomorphism to $l_2$ is similarly
 controlled. At one point in the argument we need the slightly
 stronger statement that the projections are in fact the orthogonal
projections
in order to apply Lemma 3.10. Technically our inductive hypothesis should be
strengthened to include this information.

Now assume that the result holds for all $\beta<\alpha.$
$\Rq{\alpha}$ is
equal to $(\sum_{\beta<\alpha} \Rq{\beta})_I$. Let $\beta_n=\alpha'+n$.
We only make use of the
summands $\Rq{\alpha'+n}$ for $n\in \N$ in $\Rq{\alpha'+\omega}$. Thus for
our purposes we can consider $(\sum_{n}' \Rq{\beta_n})_I$ in place of
$\Rq{\alpha}.$
Also we can replace $\Rq{\alpha'}$ by
the (codimension one) isomorphic
space $(\sum' \Rq{\beta'_l})_I$, where $(\beta'_l)$ is an
enumeration of $\{\beta<\alpha'\}.$
For each
$n\in N$ let $P_n$ be the orthogonal projection from $(\sum_l'
\Rq{{\beta_l'}})_I$ onto $(\sum_{s=1}^n{\vphantom{\sum}}'
\Rq{{\beta_s'}})_I$.
(This is in fact just
conditional expectation.) Consider the sequence of operators
$(T'_t)_t=(P_n T_r)_{n,r}$ in some order.

The proof divides into two cases depending on the nature of $\alpha'.$
First we consider the case $\alpha'=\gamma+\omega,$ for some $\gamma\geq
\omega.$
Observe that for any $n$, $\Rq{\beta'_n}$ is isomorphic to a subspace of
$\Rq{\gamma}$
and therefore, for each $n$, $(\sum_{s=1}^n{\vphantom{\sum}}'
\Rq{{\beta_s'}})_I$ is also isomorphic to a subspace of $\Rq{\gamma}$.
Therefore the inductive hypothesis applies
to $\beta_n$ for all $n$ and thus also for maps
from $\Rq{\alpha'+n}$ into $(\sum_{s=1}^n{\vphantom{\sum}}' \Rq{{\beta'_s}})_I$.
($\Rq{\alpha+n}$ contains a well complemented isometric copy of
$\Rq{\beta_s+\omega}$ for all $s\in \N$.) 

Let $\phi:\NxN $ into $\N$ be an injection such that $\phi(n,k)>n$.
For each $n,k$ let $X_{n,k}$ be an isometric (in $L_q$ norm)
copy of $\Rq{\alpha'}$ in $\Rq{\beta_{\phi(n,k)}}$ which is contractively
complemented but is supported on a set of measure
$w_{n,k}^{2p/(p-2)}$,i.e., the
norm in $L_2$ has been multiplied by $w_{n,k}$.  Therefore for each $n\in \N$
there is a normalized sequence $(z^j_{n,k})_j$ in $X_{n,k}$ such that 
$(z^j_{n,k})_j$ is equivalent to the standard basis of $l_2,$
$[z^j_{n,k}:j\in \N]$ is C-$(q,2)-$\tiny complemented,
$w_{n,k}\geq W((z^j_{n,k})_j)>C'w_{n,k}$, and $T'_r(z^j_{n,k})<2^{-j}$ for
all $r<j.$ By a diagonalization argument we can find $J_{n,k}\subset \N$
for all $n,k$ such that $(T(z^j_{n,k}))_{j\in J_{n,k},n,k}$ is a block of
the independent sum $(\sum' \Rq{\beta'_l})_I$. It follows that
$[z^j_{n,k}:j\in J_{n,k},n,k]$ is $(q,2)-$\tiny isomorphic to 
$(\sum l_2)_{q,2,(w_{n,k})}.$ Moreover, because of the diagonalization
step, the range of $T_r$ is contained in the span of a sequence of
independent random variables. Therefore by Lemma 3.6 there is a 
 a normalized
block basic sequence
$(y^j_{n,k})$ such that $[y^j_{n,k}]$ is a $(q,2)-$\tiny complemented
subspace of $[z^j_{n,k}:j\in J_{n,k},n,k]$,  $(y^j_{n,k})$ is equivalent
to the standard basis of $(\sum l_2)_{q,2,(w_{n,k})}$,
for each $n,k$, $W(y^j_{n,k})=w_{n,k}$, and $\|T_r(y^j_{n,k})\|<2^{-n-k}$
for all $r<n+k$ and all $j\in \N$. The subsequence $(y^1_{1,k})$ meets our
requirements.

If $\alpha'$ is not of the form $\gamma+\omega,$ then we need to apply
Lemma 3.10  to overcome a small difficulty and then proceed similarly to the
previous case. As before $\Rq{\alpha'}$ is isomorphic to $(\sum_n'
\Rq{\beta'_n})_I$ where $(\beta'_n)$ is an enumeration of the ordinals less
than $\alpha'.$ However these are not all isomorphic to subspaces of
$\Rq{\gamma}$ for some $\gamma<\alpha'.$ However if we consider the maps
$(T_r)$ from $X_{n,k}$ as above, we have a maps from an isometric copy of
$\Rq{\alpha'}$ with its $l_2$-norm multiplied by $w_{n,k}$ into
$(\sum_n'\Rq{\beta'_n})_I$. By the inductive hypothesis a sequence of maps
from $\Rq{\alpha'}$ into $(\sum_{n\leq m}'\Rq{\beta'_n})_I$ will satisfy
the hypothesis of Lemma 3.10. Therefore by Lemma 3.10 we can find for each $n,k$
 a normalized basic sequence $(z^j_{n,k})$ in $X_{n,k}$ such that
$(z^j_{n,k})$ is equivalent to the standard basis of $l_2,$
$[z^j_{n,k}]$ is C-$(q,2)-$\tiny complemented,
$w_{n,k}\geq W(z^j_{n,k})>C'w_{n,k}$, and $T'_r(z^j_{n,k})<2^{-j}$ for
all $r<j$. At this point the proof continues as in the previous case.
\qed\enddemo

\proclaim{Corollary 3.13} Suppose that $\omega \leq \alpha
<\alpha'<\omega_1.$ Then $\rp{\alpha}$ is isomorphic to $\rp{\alpha'}$ if
and only if $\alpha+\omega >\alpha'.$
\endproclaim

\demo{Proof} One direction is Corollary 2.10 and the other is immediate
from Proposition 3.9.
\qed\enddemo

\remark{Remark 3.14} By modifying Schechtman's proof that $\otimes^{k-1}
X_q$ does not contain
$$(\sum (\sum \dots (\sum \ell_{r_1})_{r_2}\dots
)_{r_{k-1}})_{r_k},$$
to use $D_p$ as in the argument above, we can also show
that for each $k\in \N$,
$\Rq{\omega k}$ is not isomorphic to a subspace of $\otimes^{k-1} X_q.$ The
point is that the first part of Schechtman's argument is in fact to use the
induction hypothesis to to locate a sequence of independent random
variables which contain a copy of $(\sum \ell_{r_k-1})_{r_k})$. Our
argument needs a sequence of independent random variables which contain
 a copy of $D_p$.
\endremark

\newpage

\head 4. Isomorphism from $X_p\otimes X_p$ into $(p,2)-$\tiny sums
\endhead

In this section we begin to investigate isomorphisms from $X_p\otimes X_p$ into
subspaces of $L_p$. Because, for $p>2$, $X_p\otimes X_p$ is isomorphic to a
subspace of $(\sum l_2)_p$ and $(\sum l_2)_p$ is isomorphic to a (complemented)
subspace of $X_p\otimes X_p$, one cannot say too much. However there are
some restrictions which are related to the $l_2$ structure which reveal
themselves. In later sections we will consider complemented embeddings and
obtain stronger results.

In this  section we will be working with $X_p$ as a sequence space. Thus if
$(e_n)$ is the usual basis of $X_{p,(w_n)}$, then we will write the norm as
$$\|\sum a_n e_n\|=\max\{\|\sum a_n e_n\|_2,\|\sum a_n e_n\|_p\} $$
where 
$$\|\sum a_n e_n\|_2
=(\sum a_n^2 w_n^2)^{1/2}\text{ and }\|\sum a_n e_n\|_p=(\sum a_n^p
)^{1/p}.$$
Our first result is about $X_p$ itself. (See \cite{A2} for related results.)
 Notice that it says that any operator from $X_p$ into $L_p$ acts like an
$L_2$ bounded operator on a large part of the basis.

\proclaim{Lemma 4.1}
Suppose that $(z_m)$ is a normalized standard basis of $X_p$ and $T$ is an
isomorphism of $X_p$ into $L_p.$ If $B=\{m:\|T z_m\|_2 > \|T\| \|z_m\|_2\}$
then 
$$\sum_{m\in B} \|z_m\|_2^{2p/(p-2)}\leq 1.$$
\endproclaim

\demo{Proof} Let $F$ be a finite subset of $B$ and let $y=\sum_{m\in F}
\|z_m\|^{2/(p-2)}_2 z_m.$ Then 
$$\align 
\|y\|_2 &=\bigg (\sum_{m\in F} \|z\|_2^{2p/(p-2)}\bigg )^{1/2}\\
\intertext{and}
\|y\|_p &=\bigg (\sum_{m\in F} \|z\|_2^{2p/(p-2)}\bigg )^{1/p}
\endalign $$

Therefore 
$$\align
\|T\|\bigg ( \sum_{m\in F} \|z_m\|_2^{2p/(p-2)}\bigg )^{1/2} &<
 \bigg (\sum_{m\in F} \|z_m\|_2^{4/(p-2)} \|Tz_m\|_2^2 \bigg )^{1/2}\\
&=
\bigg (\int \|T\sum_{m\in F} \|z_m\|_2^{2/(p-2)}r_m(t)z_m\|_2^2 dt\bigg )^{1/2}
\\
&\leq \bigg (\int \|T\sum_{m\in F} \|z_m\|_2^{2/(p-2)}r_m(t)z_m\|_p^p
dt\bigg )^{1/p} \\
&\leq \|T\| \max \{\|y\|_p,\|y\|_2\}.
\endalign
$$

If $\sum_{m\in F} \|z_m\|_2^{2p/(p-2)} >1$, then $\|y\|_2>\|y\|_p$ and we
have a contradiction.
\qed \enddemo

In the next result we obtain a ``type 2 inequality'' for the restriction of
an operator on
$X_p\otimes X_p$ to a large subspace. In $X_p\otimes X_p$ as a sequence
space there are four norms which appear in the expression in Proposition~1.3.
Below we need only one of them. If $(e_n)$ and $(e_n')$ are copies of the
natural basis of $X_{p,(w_n)}$, we will use the notation
$$\|\sum_{n,m} a_{n,m} e_n\otimes e_m'\|_2=
(\sum_{n,m} a_{n,m}^2 w_n^2 w_m^2)^{1/2}$$
for this $ell_2$-norm.

\proclaim{Proposition 4.2} If $T$ is an isomorphism from $X_p \otimes X_p$
into $L_p$ and $(e_n)$ and $(e_n')$ are copies of the natural basis of
$X_{p,(w_n)}$, then there are complemented subspaces, $Y$ and $Z$, of $X_p
$ each isomorphic to $X_p$, with standard $X_p$ bases, $(y_m)$ and $(z_m)$,
respectively, which are subsequences of $(e_n)$ and $(e_n')$, 
respectively, such that $\|T(y_m\otimes z_k)\|_2 \leq \|T\| \|y_m \otimes
z_k\|_2,$ for all $m$ and $k$. Consequently, for all $\sum a_{m,k
} y_m\otimes z_k \in Y \otimes Z$,
$$\bigg(\int\|T \sum a_{m,k}r_{m,k}(t)y_m \otimes z_k\|_2^2 dt\bigg)^{1/2} \leq
\|T\|\|\sum a_{m,k}y_m \otimes z_k\|_2$$
where $(r_{m,k})$ is a doubly indexed set of Rademacher functions.
\endproclaim
\demo{Proof} Let $(u_{m,n})$ and $(v_{m,n})$ be bases of $X_p$ such that
$\|u_{m,n}\|_2=\|v_{m,n}\|_2=w_m$ for all $m,n$ and $w_m \downarrow 0.$ By
the previous lemma for each $m$ and $n$ we have that if
$$\align 
B_{m,n}&=\{(k,j):\|T\|\|u_{m,n}\otimes v_{k,j}\|_2<\|T u_{m,n}\otimes
v_{k,j}\|_2\}\\
\intertext{and}
C_{m,n}&=\{(k,j):\|T\|\|u_{k,j}\otimes v_{m,n}\|_2<\|T u_{k,j}\otimes
v_{m,n}\|_2\}\endalign$$
then
$$\sum_{(k,j)\in B_{m,n}} \|u_{m,n}\otimes v_{k,j}\|_2^{2p/(p-2)} \leq
1\quad\text{and}\quad 
\sum_{(k,j)\in C_{m,n}} \|u_{k,j}\otimes v_{m,n}\|_2^{2p/(p-2)} \leq 1.$$
Also note that for all $k,j,m,n$, $(k,j)\notin B_{m,n}$ if and
only if $(m,n)\notin C_{k,j}.$

To get the subspaces $Y$ and $Z$ we will inductively choose disjoint finite
subsets $F_s$ and $G_s$ of $\Bbb N\times \Bbb N$ and infinite subsets $M_s,
N_s$
of $\Bbb N \times \Bbb N$ such that for all $s \in \Bbb N$ and for all $(m,n)\in F_s,(k,j) \in G_s$
$$ \align M_{s+1}\subset M_s &\text{ and }N_{s+1}\subset N_s,\\
\|T u_{m,n}\otimes v_{k,j}\|_2 &\leq \|T\| \|u_{m,n}\otimes v_{k,j}\|_2
,\\
M_s\supset G_t \text{ for }t\leq s&\text{ and }N_s\supset F_t \text{ for }t\leq
s,\\
|N_s \cap \{(k,r):r \in \Bbb N\}|=\infty
&\text{ and }|M_s \cap \{(k,r):r \in \Bbb N\}|=\infty
\qquad\text{for all }k.
\endalign$$
Let $F_1 =\{(1,1)\}, M_1=\NxN\setminus B_{1,1}, G_1=\{(1,j_1)\}$ for some
$(1,j_1)\in M_1$
and $N_1=\NxN \setminus C_{1,j_1}.$ Note that $(1,1)\notin
C_{1,j_1}$ because $(1,j_1) \notin B_{1,1}$, and thus $F_1\subset
N_1.$

Now suppose that we have defined $F_s, G_s, M_s$ and $N_s.$ Let
$F_{s+1}\subset N_s \setminus \cup_{t\leq s} F_t$ such that $F_{s+1}$ is finite
and $F_{s+1} \cap
\{k\}\times \Bbb N \neq \emptyset$ for $k=1,2,\dots , s+1.$ Define
$M_{s+1}= M_s \setminus \cup_{(m,n)\in F_{s+1}} B_{m,n}.$ Observe
that if $t\leq s$ then $G_t\cap \cup_{(m,n)\in F_{s+1}} B_{m,n}=
\emptyset$ because $F_{s+1}\cap \cup_{(k,j)\in G_t} C_{k,j}=
\emptyset.$ Choose a finite
set $G_{s+1} \subset M_{s+1}\setminus \cup_{t\leq s} G_t$ such that $
G_{s+1} \cap
\{k\}\times \Bbb N \neq \emptyset$ for $k=1,2,\dots , s+1.$ Define
$N_{s+1}=N_s \setminus \cup_{(k,j)\in G_{s+1}} C_{k,j}.$
Note that for $t\leq s+1$ and $(k,j)\in G_{s+1},$ $F_t \cap
C_{k,j}=\emptyset$ because $(k,j)\notin B_{m,n}$ for all
$(m,n)\in \cup_{t\leq s+1} F_t.$

This completes the induction step. Because the sets $F_s$ are disjoint and $
F_s \cap \{k\}\times \Bbb N \neq \emptyset$ for all $s \geq k,$ $\cup F_s
\cap \{k\}\times \Bbb N$ is infinite for all $k.$ Similarly the same is
true for $\cup G_s.$ Thus
$$Y=[u_{m,n}:m,n\in \cup F_s]\text{ and }Z=[v_{m,n}:m,n\in \cup G_s]$$
are isomorphic to $X_p.$
We have that $\|T u_{m,n}\otimes v_{k,j}\|_2 \leq \|T\|
\|u_{m,n}\otimes v_{k,j}\|_2$ for all $(m,n)\in \cup F_s, (k,j)\in \cup
G_s $ because $(m,n)\in F_s, (k,j)\in G_t, s\leq t,$ implies that
$G_t\subset M_s \subset \NxN \setminus B_{m,n}$ and if $t<s$ then $F_s
\subset N_t \subset \NxN \setminus C_{k,j}.$

The last statement of the conclusion follows from the fact that
$(r_{m,k}y_m\otimes z_k)$ is an orthogonal sequence and hence
$$ \align
\int \|T \sum a_{m,k} r_{m,k} y_m\otimes z_k\|^2_2 &=
\sum |a_{m,k}|^2\|Ty_m\otimes z_k\|^2_2\\
&\leq \|T\|^2\sum
|a_{m,k}|^2\|y_m\otimes z_k\|^2_2\\
&= \|T\|^2\|\sum a_{m,k}
y_m\otimes z_k\|^2_2.\endalign$$
\qed
\enddemo
\remark{Remark 4.3} If it were possible to pass to  subsequences of
the bases $(y_m), (z_k)$ for which each still spanned $X_p$ and the image
$(Ty_m\otimes z_k)$ was unconditional then the
average over signs could be removed from the conclusion of
Proposition 4.2. Unfortunately we do not know if this is
possible.
\endremark

In view of Proposition 4.2 we can usually
assume that we have passed to the
subspaces $Y$ and $Z$ and thus for the given bases of $X_p$,
$$\|T(x_i\otimes y_j)\|_2\leq \|T\|
\|x_i\otimes y_j\|_2,\tag{T2}$$
for all $i,j\in \N.$

The next result shows that diagonal blocks of the basis of $X_p\otimes X_p$
must be mapped by a bounded operator from  $X_p\otimes X_p$ into a
$(p,2)-$\tiny sum so that the norm $\|\cdot\|_2$ is well controlled.
In the computation we will use some of the other norms from the formula in
Proposition 1.3, so we introduce special notation for them.
For each $i \in \N$, let $$\Cal R_i \sum_{n,m} a_{n,m} e_n\otimes
e_m'=\sum_m a_{i,m}e_m'$$ and $$\Cal C_i \sum_{n,m} a_{n,m} e_n\otimes
e_m'=\sum_n a_{n,i}e_n. $$ Then define
$$\|\sum_{n,m} a_{n,m} e_n\otimes e_m'\|_R=(\sum_i \|\Cal R_i \sum_{n,m}
a_{n,m} e_n\otimes e_m'\|_2^p)^{1/p}=(\sum_i (\sum_m a_{i,m}^2
w_m^2)^{p/2})^{1/p}$$ and
$$\|\sum_{n,m} a_{n,m} e_n\otimes e_m'\|_C=(\sum_i \|\Cal C_i \sum_{n,m}
a_{n,m} e_n\otimes e_m'\|_2^p)^{1/p}=(\sum_i (\sum_n
a_{n,i}^2 w_n^2)^{p/2})^{1/p}.$$ Finally, define 
$$\|\sum_{n,m} a_{n,m} e_n\otimes e_m'\|_p=(\sum_{n,m} a_{n,m}^p)^{1/p}.$$
Thus if $z\in X_p\otimes X_p,$ we have
$$\|z\|=\max \{\|z\|_2,\|z\|_R,\|z\|_C,\|z\|_p\}.$$

\proclaim{Lemma 4.4} If $(Y_n)$ is a sequence of subspaces of
$L_p[0,1]$ with $C-$\tiny uncondi\-tional bases
and $T:X_{p,w}\otimes X_{p,w}\rightarrow (\sum Y_n)_{p,2}$ is an
isomorphism and $(z_n)$ is a sequence of norm 1 elements in
$X_{p,w}\otimes X_{p,w}$
such that 
\item{a.} $(z_n)$ is block diagonal, i.e., there exist increasing
sequences of integers $(m_j)$ and $(p_j)$ such that
$z_j=(Q_{m_{j+1}}-Q_{m_j})\otimes (Q_{p_{j+1}}-Q_{p_j})z_j$ for
all $j$,
\item{b.} $(Tz_n)$ is a block of the basis of $(\sum Y_n)_{p,2}$,
\item{c.} there exists an increasing sequence of integers $(M_n)$
such that 
$$\align
\|(I-P_{M_n})Tz_{n-1}\|&=0 \\
\intertext{and}
\|(I-P_{M_n})Tz_n\|_2&>\|T\| C \|z_n\|_2\qquad\text{for all
$n$}
\endalign$$
Then 
$$\bigg( \sum_1^N \|z_n\|_2^2 \bigg )^{1/2}\leq N^{1/p}\qquad\text{
for all $N$}$$
\endproclaim
\demo{Proof} We estimate the norms of $\sum z_n$ and $\sum T
z_n.$ Clearly $\|\sum z_n\|_2 =$\newline $ (\sum \|z_n\|_2^2)^{1/2}$ and 
$\|\sum z_n\|_p \leq (\sum \|z_n\|^p)^{1/p}$. Moreover, because
the sequence $(z_n)$ is block diagonal $(\sum \|z_n\|^p)^{1/p}$
dominates the mixed norms, 
$$\|\sum_n z_n\|_R=(\sum_i \|\Cal R_i \sum_n
z_n\|_{2}^p)^{1/p}\text{ and } \|\sum_n z_n\|_C=(\sum_i \|\Cal C_i \sum_n
z_n\|_{2}^p)^{1/p}.$$
 On the other hand if $(r_n)$ is a sequence
of Rademacher functions and $F\subset \N$ with $|F|=N,$
$$\align 
\|T\|\max\{N^{1/p},(\sum_{n\in F} \|z_n\|^2)^{1/2}\}&\geq
\|T\|\|\sum_{n\in F} z_n\| \\
&\geq \int \|\sum_{n\in F} r_n Tz_n\|_2 \\
&= (\sum_{n\in F} \|T z_n\|_2^2)^{1/2}\\ &\geq (\sum_{n\in F}
\|(I-P_{M_n})Tz_n\|_2^2)^{1/2}\\
 &>C\|T\| (\sum_{n\in F} \|z_n\|^2)^{1/2}.
\endalign $$
Therefore, $N^{1/p}\geq (\sum_{n\in F} \|z_n\|^2)^{1/2}.$
\qed\enddemo

Let $(w_n)$ be a sequence in $(0,1]$ which decreases to $0$ and let
$w_{n,k}=w_n$ for all $n,k\in \N.$ Throughout the next few sections
$(w_{n,k})$ will denote a doubly indexed sequence of this type.
For the space $X_{p,2,(w_{n,k})}$ only some subsequences of the basis are
again bases of $X_{p,2,(w_n')}$ for some $(w_n')$ satisfying (*). The next
definition contains a large enough class of such subsequences that we can
restrict to this class for various gliding hump arguments.

\definition{Definition 4.5} A subset $S$ of $\NxN$ is said to be {\it
rich} if there exists $M\subset \Bbb N $, $M$ infinite, such that
for every $m\in M$, $\{(m,k)\in S:k\in \Bbb N\}$ is infinite.
\enddefinition

The rich sets are a fairly nice class from a combinatorial standpoint. For
example if $S$ is a rich set and $A\subset S$ then either $A$ or $S
\setminus A$ must contain a rich set. Also if $A$ contains a rich set, 
$K$, then there is a maximal rich set $K'$ with $K\subset K' \subset A.$
Indeed $K'=\{(m,n):|\{m\}\times \N \cap A| =\infty\}.$ Another
important point for us is that rich sets can be constructed by an induction
procedure which imitates the usual method of showing that the rationals are
countable.

In $X_p \otimes X_p$ each row and column is isomorphic to $X_p$, but 
$\|\cdot\|_2$ is affected by the choice of row or column. Previous results
in this section show how the norm $\|\cdot\|_2$ in $X_p$ (Lemma 4.1) or
in $X_p \otimes X_p$ (Proposition 4.2) is ``felt'' by an
operator at least on the basis vectors. In the next result we see that it
also ``felt'' on each of the row and column spaces. This phenomenon is more
subtle since it is caused by the other rows and columns.

\proclaim{Proposition 4.6} If $(Y_n)$ is a sequence of subspaces of
$L_p[0,1]$ with $C-$\tiny uncondi\-tional bases
and $T:X_{p,(w_{n,k})}\otimes X_{p,(w_{n,k})}
\rightarrow (\sum Y_n)_{p,2}$ is an isomorphism,
then for every $\epsilon >0$ there exists a rich subset $K$ of
$\NxN$ and for each $\kappa \in K$, there is an integer $M_\kappa$ such that 
$$\|(I-P_{M_\kappa})Tz\|_2 \leq \|T\| C  w_\kappa (1+\epsilon) \|z\|$$
for all $z \in [x_\kappa] \otimes X_{p,(w_{n,k})},$
where $P_{M_\kappa}$ is the restriction operator from
$(\sum_{n=1}^\infty Y_n)$ onto $(\sum_{n=1}^{M_\kappa} Y_n).$
\endproclaim

\demo{Proof:}First we fix $m,n$ and suppose that
there is no such $M$ for some $\epsilon$ and $[x_{m,n}]\otimes
X_{p,(w_{n,k})}.$
 Then let
$(\delta_i)$ be a sequence of positive numbers tending to 0. We may
inductively choose a  sequence of unit vectors 
 $(z_{m,n,k})_{k=1}^\infty$ in
$[x_{m,n}]\otimes X_{p,w}$
and an increasing sequence of positive integers $(M_k)$ such that
$$\align 
\|(I-P_{M_k})Tz_{m,n,k-1}\|&<\delta_k \qquad\text{for all $k$}\\
\intertext{and}
\|(I-P_{M_k})Tz_{m,n,k}\|_2&>\|T\| C  w_{m,n}.
\endalign$$
By passing to a subsequence if necessary we may assume that there is a $z_0$
such that for each $k$, $z_{m,n,k}=z_0+z_{m,n,k}'$, $(z_{m,n,k}')$ is a
(perturbation of a)
block of the basis of
$[x_{m,n,k}] \otimes X_{p,(w_{n,k})},$  and $(Tz'_{m,n,k})$ is
 (a perturbation of) a block of the basis of
$(\sum Y_m)_{p,2}$ which is disjoint from $Tz_0.$
Observe that for sufficiently large $k$
$$\|(I-P_{M_k})Tz_{m,n,k}'\|_2 >(1+\epsilon)\|T\| C  w_{m,n}\|z_{m,n,k}'\|.$$
Thus we may assume that $z_0=0.$

If for some $m$ there are infinitely many $n$ for which there is no
$M_{m,n}$, then by a diagonalization argument we can find a sequence
$(z_{m,n,k(n)})_{n\in N}$ and an increasing sequence of
integers $(M_n)$ such that
$$\align 
\|(I-P_{M_n})Tz_{m,n,k(n)}\|&<\delta_n \qquad\text{for all $n$}\\
\|(I-P_{M_n})Tz_{m,n,k(n)}\|_2&>\|T\| K  w_{m,n}.
\endalign$$
We may also assume that if $z_{m,n,k(n)}=x_{m,n} \otimes \zeta_n$
then $(\zeta_n)$ is a block of the basis of $X_{p,(w_{n,k})}.$

Because $w_{m,n}$ is the same for all $n$,
Lemma 4.4 and a perturbation argument shows that this is impossible. Thus
for any $m$ there are only finitely many $n$ for which the required integer
$M_{m,n}$ does not exist, and $K$ may be obtained by discarding these
finitely many $n$ for each $m.$
\qed \enddemo
 
We will now see that the estimate on the norm $\|\cdot\|_2$ actually reduces
the tail to an $\ell_p$ sum.

\proclaim{Corollary 4.7} Suppose that $(w_{n,k})$ is a sequence of
positive numbers as above,
$T$ is an isomorphism from $X_{p,(w_{n,k})}\otimes X_{p,(w_{n,k})}$ into
$(\sum Y_n)_{p,2}$ as in Proposition 4.6. 
Let $S_M$ be the map from $(\sum Y_n)_{p,2}$ into
$(\sum_{n=M+1}^\infty Y_n)_p$ and let $P_M$ be the natural projection from
$(\sum Y_n)_{p,2}$ onto $(\sum_{n=1}^M Y_n)_{p,2}$.
Then  there exists a rich set $K\subset
\NxN$ and for each $\kappa \in K$ an integer $M_\kappa$ such that
$(S_{M_\kappa}+P_{M_\kappa})T$ is an isomorphism from $[e_\kappa]
\otimes X_{p,w}$
into $(\sum_{n=1}^M
Y_n)_{p,2}\oplus (\sum_{n=M_\kappa+1}^\infty Y_n)_p$, where $(e_{n,k})$ is
the natural basis of $X_{p,(w_{n,k})}.$
\endproclaim

\demo{Proof} First we find a rich set $K'$ as in Proposition 4.6 with
$\epsilon<1.$ Because $w_n \downarrow 0$ there is a rich set $K\subset K'$
such that $4 C \|T^{-1}\|\|T\|w_\kappa<1$ for all $\kappa \in K.$
From Proposition 4.6 we have that  for $z\in [e_\kappa]\otimes
X_{p,(w_{n,k})}$,
$$\align
\|Tz\|_2^2 &= \|(I-P_{M_\kappa})Tz\|_2^2 +\|P_{M_\kappa} Tz\|^2_2\\
&\leq \|T\|^2C^2 4 \sup w_\kappa^2 \|z\|_2^2+\|P_{M_\kappa} Tz\|^2\\
\intertext{(by Proposition 4.6 and $\|\cdot\|_2\leq \|\cdot\|$ in $(\sum
Y_n)_{p,2}$)}
&\leq \|z\|_2^2/(\|T^{-1}\|^2 4)+\|P_{M_\kappa} Tz\|^2
\tag{4.7.1}
\endalign$$
since by the choice of $K'$, $\|T\|^2 C^2 4\sup w_\kappa^2
\leq \|T^{-1}\|^{-2}/4.$
For any $z \in [e_\kappa]\otimes X_{p,(w_{n,k})}$,
$$\align
\|T^{-1}\|^{-1}\|z\|
&\leq \|Tz\|\\&=\max\{\|Tz\|_2,(\sum
\|(P_{n+1}-P_n)Tz\|_p^p)^{1/p}\}\\
&\leq \max\{(\|T\|^2C^2 4\sup w_\kappa^2
\|z\|_2^2+\|P_{M_\kappa} Tz\|^2)^{1/2},\\
&\qquad
\|P_{M_\kappa} Tz\|+(\sum_{n=M_\kappa+1}^\infty
\|(P_{n+1}-P_n)Tz\|_p^p)^{1/p}\}
\endalign$$
If $$(\|T\|^2C^2 4\sup w^2 \|z\|_2^2+\|P_{M_\kappa} Tz\|^2)^{1/2}>\|Tz\|\geq
\|T^{-1}\|^{-1}\|z||$$
then by (4.7.1)
$$\|T^{-1}\|^{-2}\|z\|^2\leq \|z\|_2^2/(\|T^{-1}\|^2 4)+\|P_{M_\kappa}
Tz\|^2$$ and consequently
$$\|T^{-1}\|^{-2}\|z\|^2(3/4)\leq \|P_{M_\kappa} Tz\|^2.$$
Thus $$\|T\|\|T^{-1}\|
\|P_{M_\kappa} Tz\|\sqrt{12}/3\geq \|Tz\|$$ in this case.

It follows that $\|((S_{M_\kappa}+P_{M_\kappa})T)^{-1}\|\leq \|T\|\|T^{-1}\|
\sqrt{12}/3.$ \qed\enddemo

\remark{Remark 4.8}
For a general isomorphism from $X_p\otimes X_p$ into $(\sum Y_n)_{p,2}$, we
cannot hope to show that there is a bound $\infty>M\geq M_\kappa$. To see
this notice that there is a simple isomorphic embedding of $X_{p,(w_n)} \otimes
X_{p,(w_n)}$ into $(\sum X_{p,(w_n)})_{p,2,(w_n)}\oplus (\sum
X_{p,(w_n)})_{p,2,(w_n)}.$ Let $(e_n)$ and $(e_n')$ be copies of the
natural basis of $X_{p,(w_n)}.$ Let $(d_{k,j})$ and $(d'_{k,j})$ be two
copies of the natural basis of $(\sum X_{p,(w_n)})_{p,2,(w_n)}$. Define
$T(e_n\otimes e_m')=d_{n,m}\oplus d_{m,n}'$ and extend linearly. It is easy
to see from Proposition 1.3 and the definition of the $(p,2)-$\tiny sum
( or Corollary 2.9) that $T$ is an isomorphism. Obviously no uniform bound
on $(M_\kappa)_{\kappa \in K}$ can be found for this operator and any rich
set $K$.
\endremark

\newpage

\head 5. Selection of bases in $X_p\otimes X_p$ \endhead

Because of the multi-index nature of the natural basis of $X_{p,(w_{n,k})}
\otimes X_{p,(w_{n,k})}$
and the technical complexities of gliding hump type arguments with respect
to a multi-index, 
we introduce a method of producing subsequences of the basis which
still span a copy of $X_{p,(w_{n,k})}\otimes X_{p,(w_{n,k})}$
 but which can be used without
directly worrying about the nature of the underlying index set. We will use
a fairly general setup in this section which may be applicable to other
bases with complicated natural orderings.

\definition{Definition 5.1} Let $X$ be a Banach space with basis $(x_i)$.
A set $\Cal S$ of infinite subsets of
$\N$ is {\it
$K$-admissible for $(x_i)$} if 
\item{(1)} For each $n\in \N$ and $S \in \Cal S$,
$\{n, n+1, \dots \}\cap S \in \Cal S$.
\item{(2)} $\N \in \Cal S.$
\item{(3)} For each $S\in \Cal S,$
 $(x_i)_{i\in S}$ is a basis
for a subspace of $X$ which is $K$-isomorphic to $X$.
\enddefinition

\definition{Definition 5.2} Suppose that $X$ is a Banach space with basis
$(x_i),$ $\tau$ is a topology on $X$, and $\Cal S$ is a admissible
subset of $2^\N$ for $(x_i)$. We will say that $(x_i)$ has SP
 (selection property) with respect to $\tau$ and $\Cal S$
if there is a winning strategy for the second player in the two
player game described below.
\item{(0)} $S_0'=\N.$
\item{(1)} On each turn $n$, $n=1,2,\dots,$ the first player must define
 a multi-function
$ F_n$ on $X$ with range $F_n \subset
\{0,1,2,\dots,N_n\}$, $N_n$ 
finite,  and $F_n^{-1}(j)$ is $\tau$-open for each $j$
and a set $S_n\in \Cal S$ with $S_n \subset S_{n-1}'$ and
$\{i_1,i_2,\dots,i_{n-1}\}\subset S_n$.
\item{(2)} On the turn $n$ the second player must choose an
integer $i_n \in S_n$, $i_n>i_{n-1},$ and
a set $S_n' \in \Cal S$ with $S_n'\subset S_n$
and $\{i_1,i_2,\dots,i_n\}\subset S_n'$.

Player 2 wins if (and only if) 
$\{x_{i_k}:k\in \N\}$ is in $\Cal S$ and 
$F_n(x_{i_j})$ is constant for all $j\geq
n.$
\enddefinition

In the game the function $F_n$ given by Player 1 defines a $\tau$-open
cover of $X$ and Player 2 is forced to choose one of the sets, $O$, 
from the open
cover and select all further elements for the subsequence from $O$.
Usually it will not be necessary to define the multi-functions $F_n$ on
anything more than the $\tau$-closure of the basis $(x_i).$

In $X_p$ the most useful basis is indexed by $\N\times
\N.$ To treat such cases 
we will use an order like that used to
enumerate the
rationals. To this end let $\phi$ be a bijection from $\N$ onto
$\N \times
\N$ such that $\phi(j)^1+\phi(j)^2\leq \phi(i)^1+\phi(i)^2$
if $j \leq
i$. (Here we use superscripts to denote the coordinates of
elements in $\N\times \N,$ e.g., $(n,k)^1=n.$)
We order the basis of $X_{p,(w_{n,k})}$ as
$(x_{\phi(i)})_{i=1}^\infty.$
The admissible set $\Cal S$ in this case is $\phi^{-1}(\{M \subset \N
\otimes \N:  \text{ there is an infinite }N' \subset \N
\text{ such that for each }n \in N', (n,k)\in M \text{ for
infinitely many }k\text{ and }(n,k)\notin M \text{ for all }n\notin
N'\text{ and }k\in \N\}).$

If $M$ is an infinite subset of $\Bbb N$ and $i\in M$ then $o(i)$
will be the ordinal of the element in $M$ as an ordered set with the
inherited order. If there is some ambiguity about which set is under
consideration, we will add the set as another parameter as in $o(i,M).$

\proclaim{Lemma 5.3} The standard basis of $X_{p,(w_{n,k})}$, where
$(w_{n,k})$ satisfies (*), has the SP with respect to the weak
topology and the class $\Cal S$ defined above.
\endproclaim

\demo{Proof} 
Let $(x_{n,k})$ be the standard basis of $X_{p,(w_{n,k})}$.
 We will use induction to choose a subsequence $(x_{n,m})_{n\in \Bbb N, m
\in M_n}$ with $\{(n,m):n \in \N, m\in M_n\}\in \phi(\Cal S).$

For notational convenience we will assume that 
$0\in F_m^{-1}(0),$ for each $m.$ If $A\in \phi(\Cal S)$, $A^1$ will denote
$\{n:(n,k)\in A\text{ for some }k\},$ the projection of $A$ into the first
coordinate.

Consider $F_1, S_1.$ 
We have assumed that $0\in F_1^{-1}(0)$.
Because $\text{w-}\lim_k x_{n,k}=0$, for each $n\in \phi(S_1)^1$
there exists $M_{n,1}\subset (n\times\N)\cap \phi(S_1)$,
$M_{n,1}$ infinite, such that $F_1(x_{n,k})=0$ for all $n$ and all $k\in
M_{n,1}.$ Let 
$i_1=\phi^{-1}((n_1,k))$ where $o(n_1,\phi(S_1)^1)=1=\phi(1)^1$ and 
$o(k,M_{n_1,1})=1.$ Define
$S_1'=\phi^{-1}(\cup_{n\in \phi(S_1)^1} M_{n,1}).$

Let Player 1 choose the multi-function
$F_2$ and $S_2\subset S_1',$ with $S_2 \in \Cal S.$
Because $F_2(0)=0$
and $\text{w-}\lim_{k:(n,k)
\in M_{n,1}} x_{n,k}=0$, for each $n\in \phi(S_2)^1$ we can
find  an infinite subset  $M_{n,2}$ of $M_{n,1}\cap \phi(S_2)$ such that
$F_2(x_{n,k})=0$ and $i_1<\phi^{-1}(n,k)$
for all $k \in  M_{n,2}$. Let $i_2 =
\phi^{-1}(n_2,k)$ where $o(n_2,\phi(S_2)^1)=\phi(2)^1$ and
$o(k,M_{n_2,2})=\phi(2)^2$ and let
$S_2'= \phi^{-1}(\cup_{n\in \phi(S_2)^1} M_{n,2}).$

Suppose we have chosen $i_1 < i_2 < \dots < i_{m-1},$ and
$M_{n,m-1} =\phi(S_{m-1}')\cap (n\times \N),$ for all $n\in
\phi(S_m)^1.$ Let Player 1 choose $F_m$ and $S_m.$
Because $F_m(0)=0$
and $\text{w-}\lim_{k:(n,k)
\in M_{n,m-1}} x_{n,k}
=0$, for each $n\in \phi(S_m)^1$ we can
find  an infinite subset  $M_{n,m}$ of $M_{n,m-1}\cap \phi(S_m)$ such that
$F_m(x_{n,k})=0$ and
$i_{m-1}<\phi^{-1}(n,k)$ for all $k \in 
M_{n,m}$. Let $i_m =
\phi^{-1}(n_m,k)$ where $o(n_m,\phi(S_m)^1)=\phi(m)^1$ and
$o(k,M_{n_m,m})=\phi(m)^2$ and let
$S_m'= \phi^{-1}(\cup_{n\in \phi(S_m)^1} M_{n,m}).$

In this way the subsequence $(x_{\phi(i_j)})$ satisfies
$F_m(x_{\phi(i_j)})=0$ for all $j\geq m.$ Moreover since
for any $j\geq m$, $o(\phi(i_m)^1,\phi(S_j)^1)=\phi(m)^1$
 and $o(\phi(i_m)^2,\{\phi(i_j)^2:\phi(i_j)^1=
\phi(m)^1\})=\phi(i_m)^2=o(\phi(m)^2,\N),$ $\{i_j:j\in \N \}$ is in
$\Cal S.$
\qed \enddemo

In the next result we pass to a tensor product. We will not define a class
of admissible sets for the product index set but instead use the game in
each factor to accomplish our goals.
Below we use an unspecified tensor product of two Banach spaces. The
only property that we require of the tensor product is that
$\|x \otimes y\|\leq
\|x\|\|y\|$ for all $x,y.$

\proclaim{Proposition 5.4} Suppose that $(x_i)$ and $(y_j)$ are shrinking
unconditional bases with SP with respect to the weak topology
and admissible sets $\Cal S(X)$ and $\Cal S(Y)$, respectively,
 and that  $T$ is a bounded operator from $[x_i
\otimes y_j: i,j \in \Bbb N]$ into a space $Z$ with normalized shrinking
basis $(z_k).$ Then
given $\epsilon > 0$ there are $I \in \Cal S(X)$ and $J \in \Cal S(Y)$
and finite subsets $N_{i,j}$ of $\Bbb N$
such that 
\item{1)} $\sum_{i\in I,j\in J}\|T(x_i \otimes y_j)|_{\Bbb N\setminus 
N_{i,j}}\| < \epsilon $
\item{2)} $N(i,j)\cap N(i',j')=\emptyset$ if $i'\neq i$ and
$j\neq j'$; $i=i'$, $j\neq j'$ and $\max(o(j),o(j'))>o(i)$;
or $i\neq i'$, $j=j'$ and $\max(o(i),o(i'))\geq o(j).$ 
\endproclaim
\demo{Proof:} First let $\epsilon(i,j)>0$ such that
$\sum_{i,j}\epsilon(i,j)<\epsilon,$ and $\epsilon(i,j)$ is
decreasing, i.e., $\epsilon(i,j)\leq \epsilon(i',j')$ if $i\leq
i', j\leq j'.$ We will use two interweaving
games to choose the
sets $I$ and $J$ and subscript the associated functions by $X$ or $Y$ to
keep the notation straight. We begin the games with only
trivial conditions: $F_{X,1}$ and $F_{Y,1}$ are constant functions
and $S_{X,1}=S_{Y,1}=\N.$
Let $i_1,j_1$ be the first elements chosen in each of the games and
$S_{X,1}',S_{Y,1}'$ be the resulting elements of $\Cal S(X)$ and
$\Cal S(Y)$, respectively.
Choose a finite subset $N_{1,1}$ of $\N$ such that 
$\|T(x_{i_1}\otimes y_{j_1})|_{\Bbb N\setminus N_{1,1}}\|<\epsilon(1,1)$.

Let $\eta_1=\{1,2,\dots,\max N_{1,1}\}.$
Define $F_{Y,2}(x_{i_1},y)=0$ if $\sum_{k\in \eta_1}
 |z^*_k(Tx_{i_1}\otimes y)|
<\epsilon(1,2)/2$ and $F_{Y,2}(y)=1$ if $\sum_{k\in \eta_1}
|z^*_k(Tx_{i_1}\otimes y)|
>\epsilon(1,2)/3$. Let $S_{X,2}=S_{X,1}'$
 and $S_{Y,2}=S_{Y,1}'.$

Player 2 must now select $j_2$ such that $F_{Y,2}(y_{j_2})=0$
because $(T(x_{i_1}\otimes y_j))$
converges to 0 weakly. Choose a finite subset $N_{1,2}$ of $\N\setminus
\eta_1$
such that 
$\|T(x_{i_1}\otimes y_{j_2})|_{\Bbb N\setminus
N_{1,2}}\|<\epsilon(1,2)$.

Let $\eta_2=\{1,2,\dots,\max( N(1,1)\cup N(1,2))\},$
$$F_{X,2}(x)=0 \text{ if }\max_{s=1,2}
\sum_{k\in \eta_2 }|z_k^*(T(x\otimes y_{j_s}))|
<\epsilon(2,2)/2$$ and
$$F_{X,2}(x)=1\text{ if }\max_{s=1,2}
\sum_{k\in \eta_2 }|z_k^*(T(x\otimes y_{j_s}))|
>\epsilon(2,2)/3.$$
Let $S_{X,2}=S_{X,1}'.$

Player 2 in the $X$ game must choose $i_2$ such that $F_{X,2}(x_{i_2})=0$.
For $k=1,2$ choose a finite set $N(2,k) \subset \N \setminus \eta_2
,$ such that 
$\|T(x_{i_2}\otimes y_{j_k})|_{\Bbb N\setminus
N_{2,k}}\|<\epsilon(2,k)$.

Let $\eta_3=\{1,2,\dots,\max \cup_{s\leq 2, t\leq 2} N(s,t)\},$
$$F_{Y,3}(y)=0\text{ if }\max_{n'=1,2} 
\sum_{k\in \eta_3}
|z_k^*(T(x_{i_{n'}}\otimes y))|
<\epsilon(2,3)/2$$ and
$$F_{Y,3}(y)=1\text{ if }\max_{n'=1,2} 
\sum_{k\in \eta_3}
|z_k^*(T(x_{i_{n'}}\otimes y))|
>\epsilon(2,3)/3.$$
Let $S_{Y,3}=S_{Y,2}'.$

Player 2 in the $Y$ game must choose $j_3$ such that $F_{Y,3}(y_{j_3})=0$.
For $k=1,2$ choose a finite set $N(k,3) \subset \N \setminus \eta_3
$ such that 
$\|T(x_{i_k}\otimes y_{j_3})|_{\Bbb N\setminus
N_{k,3}}\|<\epsilon(k,3)$.

This completes the first few steps of the induction. Assume that $i_1,
\dots, i_r$ and $j_1,\dots j_{r+1}$ are known and the corresponding sets
$N(n,m),n=1,2,\dots r,m=1,2,\dots , r+1$ have been chosen.

We continue  with the $r+1$ turn of the $X$ game. 
Let $$\eta_{2r}=\{1,2,\dots,\max \cup_{s\leq r,t\leq r+1}
N(s,t)\}$$ and
$$F_{X,r+1}(x)=0 \qquad\text{ if } \max_{s\leq r+1}
\sum_{k\in \eta_{2r}}
|z_k^*(T(x\otimes y_{j_s}))|
<\epsilon(r+1,r+1)/2$$ and
$$F_{X,r+1}(x)=1\text{ if }\max_{s\leq r+1}
\sum_{k\in \eta_{2r}}
|z_k^*(T(x\otimes y_{j_s}))|
>\epsilon(r+1,r+1)/3.$$ Let $S_{X,r+1}=S_{X,r}'.$

Player 2 in the $X$ game must choose $i_{r+1}$ such that
$F_{X,r+1}(x_{i_{r+1}})=0$.
For $k=1,2,\dots,r+1$ choose
a finite set $N(r+1,k) \subset \N \setminus \eta_{2r}
$ such that 
$\|T(x_{i_2}\otimes y_{j_k})|_{\Bbb N\setminus
N_{r+1,k}}\|<\epsilon(r+1,k)$.

Let $\eta_{2r+1}=\{1,2,\dots,\max \cup_{s\leq r+1,t\leq r+1}
N(s,t)\}$ and
$$F_{Y,r+2}(y)=0 \text{ if } \max_{n'\leq r+1} 
\sum_{k\in \eta_{2r+1}}
|z_k^*(T(x_{i_{n'}}\otimes y))|
<\epsilon(r+1,r+2)/2$$ and
$$F_{Y,r+2}(y)=1 \text{ if } \max_{n'\leq r+1} 
\sum_{k\in \eta_{2r+1}}
|z_k^*(T(x_{i_{n'}}\otimes y))|
>\epsilon(r+1,r+2)/3$$
Let $S_{Y,r+2}=S_{Y,r+1}'.$

Player 2 in the $Y$ game must choose $j_{r+2}$ such that
 $F_{Y,r+2}(y_{j_{r+2}})=0$.
For $k=1,2,\dots, r+1$ choose a finite set
$N(k,r+2) \subset \N \setminus \eta_{2r+1}
$ such that 
$$\|T(x_{i_k}\otimes y_{j_{r+2}})|_{\Bbb N\setminus
N(k,r+2)}
\|<\epsilon(k,r+2).$$

It is easy to see that we have now completed the next step of the induction
and the result follows.
\qed
\enddemo

\newpage

\head 6. $X_p\otimes X_p$-preserving operators on  $X_p\otimes X_p$ \endhead

The purpose of the section is to prove a criterion which guarantees that an
operator on $X_p\otimes X_p$ preserves isomorphically a copy of the space.
Similar results exist for spaces with an unconditional basis which have
many subsequences equivalent to the original basis. The tensor product
makes the combinatorics more difficult here.
In Proposition 5.4 we were unable to get completely disjoint blocks. 
In order to use the SP for the basis of each factor of
$X_p\otimes X_p$ to get really disjoint blocks it is necessary to
have some quantitative information. The next two lemmas give us estimates
of how many vectors we will need in order to insure that we can get
reasonable disjointness for a step of a gliding hump argument.

In the next lemma we make use of the concept of a lower $\ell_r$-estimate. We say
that a basic sequence $(z_n)$ has a lower $\ell_r$-estimate with constant
$C$ if (and only if)
$$\|\sum a_n z_n\|\geq C (\sum a_n^r)^{1/r}$$ for every
sequence of real numbers $(a_n).$

\proclaim{Lemma 6.1} If $T$ is an operator from $X_{p,(w_n)}$ into a space
$Z$ with a basis $(z_n)$ and normalized
biorthogonal functionals $(z_n^*)$ such that
for some $r\geq 2,$ $(z_n)$ satisfies a lower $\ell_r$ estimate with constant C,
then for any $k\in \N$ and $\|T\|>\epsilon>0$,
$$\sum_{n\in F} w_n^{2p/(p-2)}
\leq \|T\|^r 2k^{r+1}/(\epsilon^r C^r \min\{w_i^r:1\leq i \leq k\}),$$
where $(x_n)$ is the standard basis of $X_{p,(w_n)}$ and
$$F=\{n: |z_i^*(Tx_n)|>\frac{\epsilon w_n}{k} \text{ or }
|z_n^*(Tx_i)|>\frac{\epsilon w_i}{k}\text{ for some }i\leq k\}.$$
\endproclaim
\demo{Proof} For $1\leq i \leq k$ let
$F_i=\{n:|z_i^*(Tx_n)|>\frac{\epsilon w_n}{k}\}$ and let $G_i=
\{n:|z_n^*(Tx_i)|>\frac{\epsilon w_i}{k}\}.$ Then for $n\in F$,
$n\in F_i$ or $G_i$ for at
least one $i$.
For $n\in F_i$ let $a_n=w_n^{(p+2)/(p-2)}.$ Then
$$\align
\|\sum_{n\in F_i} a_n x_n\|&= \max \{(\sum_{n\in F_i}a_n^2
w_n^2)^{1/2},(\sum_{n\in F_i}a_n^p)^{1/p}\}\\
&= \max \{(\sum_{n\in F_i}
w_n^{4p/(p-2)})^{1/2},(\sum_{n\in F_i}w_n^{p(p+2)/(p-2)})^{1/p}\}\\
&\leq  \max \{(\sum_{n\in F_i}
w_n^{2p/(p-2)})^{1/2},(\sum_{n\in F_i}w_n^{2p/(p-2)})^{1/p}\}.
\endalign$$
Let $W_i=\sum_{n\in F_i}w_n^{2p/(p-2)}.$
For each $i$ and choice of signs $\sigma_n$,
$$|z_i^*(T\sum_{n\in F_i}\sigma_n a_n x_n)|\leq
\|T\|\max\{W_i^{1/2},W_i^{1/p}\}.$$ Thus 
$$\epsilon
W_i/k=\sum_{n\in F_i}\epsilon w_n^{2p/(p-2)}/k
 = \sum_{n\in F_i}\epsilon w_n
|a_n|/k \leq
\|T\|\max\{W_i^{1/2},W_i^{1/p}\}.$$
If $W_i>1$ then
$W_i^{1/2}>W_i^{1/p}$ and hence $W_i\leq \|T\|^2 k^2/\epsilon^2;$ if
$W_i\leq 1$, 
$$W_i\leq \|T\|^{p/(p-1)}k^{p/(p-1)}/\epsilon^{p/(p-1)}\leq
\|T\|^2 k^2/\epsilon^2$$
 also.

For each $i\leq k$, $$\|Tx_i\|\geq C(\sum_{n\in G_i}|z_n^*(Tx_i)|^r)^{1/r}
\geq C\epsilon w_i k^{-1} |G_i|^{1/r}
\geq C\epsilon w_i k^{-1}
(\sum_{n\in G_i}w_n^{2p/(p-2)})^{1/r}.$$
Thus $\sum_{n\in G_i}w_n^{2p/(p-2)}\leq k^r\|T\|^r/(\epsilon^r C^r w_i^r).$
Let $W_i'=\sum_{n\in G_i} w_n^{2p/(p-2}.$

Then $$\align
\sum_{n\in F} w_n^{2p/(p-2)}&\leq\sum_{i=1}^k W_i+W_i'\\
&\leq
\sum_{i=1}^k(\frac{\|T\|^2 k^2}{\epsilon^2}+\frac{\|T\|^r k^r}{\epsilon^r
C^r w_i^r})\\
&\leq
\frac{\|T\|^r 2k^{r+1}}{\epsilon^r C^r \min\{w_i^r:1\leq i \leq k\}}.
\endalign$$
\qed\enddemo

\medskip

Notice that if a finite sequence of positive
numbers $(\epsilon_i)$ is specified in advance and some control
on $\min \{w_i^r:1\leq i \leq k\}$ is given, it is possible to predict the
number of elements required to produce an approximately blocked
image. To be more precise we have the following.

\proclaim{Lemma 6.2} Let $(\epsilon_i)$ be a sequence of positive
numbers, $C,D,w_0 \in \R^+$ and $K\in \N.$ There exists an integer
$N_0,$ such that
if $T$ is an operator from $X_{p,(w_n)_{n=1}^{N_0}}$ into a space
$Z$ with a basis $(z_n)$ and biorthogonal functionals $(z_n^*)$
such that
for some $r\geq 2,$ $(z_n)$ satisfies a lower $\ell_r$ estimate
with constant $C,$ $w_n\geq w_0$ for all $n\leq N,$ and
$\|T\|\leq D,$
then  there exist $\{n_j:1\leq j \leq K\}$ such that for $1\leq
i < j\leq K,$
$$|z_{n_i}^*(Tx_{n_j})|<\epsilon_j w_{n_j}/(j-1)$$
and 
$$|z_{n_j}^*(Tx_{n_i})|<\epsilon_i w_{n_i}/(j-1).$$
\endproclaim

\demo{Proof} We use induction on $K.$
Notice that if $K=1$, the requirements are vacuous, i.e., $N_0=1$
works. If true for
$K$,  let $N(K)$ denote the required integer. We need a sufficiently
large $N(K+1)$ to apply Lemma 6.1 with $\epsilon=\epsilon_{K+1}$ (We may
assume that $\epsilon_{K+1}=\min\{\epsilon_j:1\leq j \leq
K+1\}.$), $k=K$, $z_i=z_{n_i}$, and the basis of $X_{p,(w_n)}$
reordered so that $x_{n_i}$ is the $i$th basis vector and the
remaining  basis vectors follow these. Then, because
$\sum_{n=K+1}^{N(K+1)} w_n^{2p/(p-2)} \geq \sum_{n=K+1}^{N(K+1)}
w_0^{2p/(p-2)},$ if
$(N(K+1)-(K+1)) w_0^{2p/(p-2)} >D^r 2K^{r+1}/(\epsilon_{K+1}
w_0 C)^r,$ the inductive hypothesis gives us $\{n_i:1\leq i \leq k\}$
and Lemma 6.2
 will produced $n_{k+1}\in \{K+1,\dots,N(K+1)\}.$ \qed\enddemo

With these lemmas we can now show that operators on $X_p\otimes X_p$ with a
significant diagonal are isomorphisms on a subspace isomorphic to the whole
space. The proof makes use of the basic technique of \cite{CL}. The problem
here is to overcome the technical difficulties in getting a large block
basic sequence.

\medskip

\proclaim{Proposition 6.3} Suppose that $T$ is a bounded operator
on $X_{p,(w_{n,k})}\otimes X_{p,(w_{n,k})}$ such that there
exists $\epsilon >0$ with $|(x_{n,k}^*\otimes
y_{m,j}^*)(Tx_{n,k}\otimes y_{m,j})|\geq \epsilon$ for all
$n,k,m,j.$ Then there are rich subsets $K,J$ of $\N\times \N$
such that $T|_{[x_{n,k}\otimes y_{m,j}:(n,k)\in K, (m,j)\in J]}$
is an isomorphism. Moreover, $T([x_{n,k}\otimes y_{m,j}:(n,k)\in K,
(m,j)\in J])$ is complemented.
\endproclaim
\demo{Proof} First observe that if $J,K$ are rich subsets of $\N\times\N$
and $$Z={[x_{n,k}\otimes y_{m,j}:(n,k)\in K, (m,j)\in J]},$$
then we can restrict our attention to $Z.$
Indeed, because the basis of $X_{p,(w_{n,k})}\otimes X_{p,(w_{n,k})}$
is unconditional, we can compose
$T|_Z$ with the basis
projection $P$ onto $Z$
and a diagonal operator $D$ on $Z$ 
 such that $D(x_{n,k}\otimes
y_{m,j})=((x_{n,k}^*\otimes
y_{m,j}^*)(Tx_{n,k}\otimes y_{m,j}))^{-1}x_{n,k}\otimes
y_{m,j}$ for all $n,k,m,j.$ The result will be obtained by
showing that $K,J$ can be chosen so that the resulting
composition is a perturbation of the identity on $Z.$
As in \cite{CL} this will imply that $T(Z)$ is complemented since for any $z\in
Z$, $(TDP)Tz=T(DPTz)\approx Tz.$

By Proposition 5.4 for any $\delta>0$ we may find rich sets $K',J'$ and finite
subsets of $\N^4$, $N_{k,j},k\in K', j\in J'$ such that
\item{(1)} $\sum_{k\in K',j\in J'}\|T(x_k \otimes y_j)|_{
\N^4\setminus
N_{k,j}}\| < \delta$
\item{(2)} $N(k,j)\cap N(k',j')=\emptyset$ if $k'\neq k$ and
$j\neq j'$; if $k=k'$ and $\max(o(j),o(j'))>o(k)$;
or if $j=j'$ and $\max(o(k),o(k'))\geq o(j).$

If $\delta$ is sufficiently small, it follows from standard arguments that
if $\eta_{n}=\{j:o(j)\leq o(n)\} $ and $\eta'_{m}=\{k:o(k)<o(m)\}
$ for all $n\in K',m\in J',$ then $T(Z)$ is isomorphic to the UFDD
$$\sum_{n\in K'}
[T(x_n\otimes y_j)|_{\cup\{N(n,j):j\in \eta_n}:j\in \eta_n]\oplus
\sum_{m\in J'}
[T(x_k\otimes y_m)|_{\cup\{N(k,m):k\in \eta'_m}:k\in \eta'_m].$$
In particular, if $\delta <\epsilon$, then $(i,j)\in N(i,j)$ for all $i,j.$
We need to refine the index sets $J',K'$ further to get an actual
(perturbation of a ) basis rather than a UFDD. To save
the notational burden of carrying a set of error terms through the computation
we will assume that $\|T(x_k \otimes y_j)|_{
\N^4\setminus
N_{k,j}}\|=0$ for all $k\in K', j\in J'.$

Let $(\epsilon_i)$ be a sequence of positive numbers such that
$\sum \epsilon_i<\epsilon'/4^p,$ where $\epsilon'<\min(1,\epsilon).$
We will use the estimate in Lemma 6.1
repeatedly in the proof, but its exact nature is not used. The important
thing is the dependence only on $r,C,\epsilon,k,\min\{w_i:i\leq k\},$ and
 $\|T\|.$ In this proof $r=p$, $C=1$ and $\|T\|$ are fixed,
so we will denote the
estimate by $H(\epsilon,w,k),$ where $w=\min\{w_i:i\leq k\}.$ 

We will essentially follow the two game argument of Proposition 5.4
but this time make more careful choices. We will again make use
of the admissible class $\Cal S$ for $X_{p,(w_{n,k})}$. As before the
selection of $i_1,j_1$ is not controlled. The choice $j_2$ is more critical
since we must make the contribution to the support of
$Tx_{\phi(i)}\otimes y_{\phi(j_1)}$
small for all $i$ in some set in the class $\Cal S(X).$

For each $j\in S_{Y,1}'$, $j>j_1,$
let 
$$\multline
K_j=\{i:i\in S_{X,2}',
|(x_{\phi(i)}^*\otimes y_{\phi(j_1)}^*)(Tx_{\phi(i)}\otimes
y_{\phi(j)})|<\epsilon_1 w_{\phi(j)},\\
\text{ and }
|(x_{\phi(i)}^*\otimes y_{\phi(j)}^*)(Tx_{\phi(i)}\otimes
y_{\phi(j_1)})|<\epsilon_1 w_{\phi(j_1)}\}.
\endmultline$$
We need to show that $\{j:K_j\supset K_j'\backepsilon
i_1\in  K_j'\text{ and }K_j'\in \Cal S(X)\}\cup\{j_1\}$
contains a set $S_{Y,2}\in \Cal S(Y)$, with $j_1 \in S_{Y,2}$ for
which $\phi(S_{Y,2})\cap(\{\phi(j_1)^1\}\times \N)$ is infinite.
We will first show that  for almost all $j$ there exists $K_j'\subset K_j$
such that $K_j' \in \Cal S(X).$ For each $i$, by applying Lemma 6.1  to the
operator $T_i:[x_{\phi(i)}\otimes y_{\phi(j)}:j\in \N]\rightarrow
[x_{\phi(i)}\otimes y_{\phi(j)}:j\in \{j_1\}+(\N\setminus \{j_1\})]$ defined by
$$T_i(z)=x_{\phi(i)}^*\otimes
y_{\phi(j_1)}^*(Tz)x_{\phi(i)}\otimes y_{\phi(j_1)}
+\sum_{j\neq j_1} x_{\phi(i)}^*\otimes y_{\phi(j)}^*(Tz)
x_{\phi(i)}\otimes y_{\phi(j)},$$
we have that $\sum_{j:i \notin K_j}
w^{2p/(p-2)}_{\phi(j)}<H(\epsilon_1, \{w_{\phi(j_1)}\},1).$ 
If for some $j$, $K_j$ contains no subset which is
in $\Cal S(X)$, then $\phi(K_j) \cap (\{m\}\times\N)$ is finite for
all but finitely many $m.$  Lemma 6.1 implies that
$\sum w_{\phi(j)}^{2p/(p-2)}$
over such $j$ must be no more than $H(\epsilon_1, \{w_{\phi(j_1)}\},1).$ 
Indeed,
if the sum exceeded this bound then there would be a finite set
$F$ of such $j$ such that $\sum_{j\in F} w_{\phi(j)}^{2p/(p-2)}>
H(\epsilon_1, \{w_{\phi(j_1)}\},1).$
Then since no subset of $K_j$ is in $\Cal S(X),$
there would be infinitely many $m$ for which $\phi(K_j) \cap
(\{m\}\times\N)$ is finite for all of those $j\in F$. Thus for any $i\in
S_{X,2}\setminus \cup_{j\in F} K_j$
$T_i$ would violate the conclusion of Lemma 6.1.
In particular 
for those $j$ such that $\phi(j)\in (\{\phi(j_1)^1\}\times \N)$,
we have that $w_{\phi(j)}=w_{\phi(j_1)}$ and thus only finitely many of these
$K_j$ will fail to contain a subset $K_j'$ as above.
 Similarly,
if
$F\subset \{j:K_j\not\supset K_j'\backepsilon i_1\in K_j',
 K_j'\in \Cal S(X)\text{ but } \exists K \subset K_j,
K \in S(Y)\}$ and $\sum_{j\in F}
w_{\phi(j)}^{2p/(p-2)}>
H(\epsilon_1, \{w_{\phi(j_1)}\},1)$ then $\phi(K_j) \cap
(\{\phi(i_1)^1\}\times\N)$ is finite
for all of those $j\in F$, thus
for any $i\in
S_{X,2}'\setminus \cup_{j\in F} K_j$ with $\phi(i)^1=\phi(i_1)^1$,
$T_i$ would violate the conclusion of Lemma 6.1.
It follows that the
required set $S_{Y,2}$ exists. 

Player 2 selects $j_2$ from $S_{Y,2}$ and selects $S_{Y,2}'\subset S_{Y,2}$
which is in $\Cal S(Y)$ and contains $\{j_1,j_2\}.$ By our choice of
$S_{Y,2}$, $K_{j_2}$ contains $K_{j_2}'$ such that $i_1\in K_{j_2}'\in
\Cal S(Y).$ Observe that this implies that
$$\align
|(x_{\phi(i_1)}^*\otimes y_{\phi(j_1)}^*)(Tx_{\phi(i_1)}\otimes
y_{\phi(j_2)})|&<\epsilon_1 w_{\phi(j_2)},\\
\intertext{ and }
|(x_{\phi(i_1)}^*\otimes y_{\phi(j_2)}^*)(Tx_{\phi(i_1)}\otimes
y_{\phi(j_1)})|&<\epsilon_1 w_{\phi(j_1)}.
\endalign$$
The set $S= K_{j_2}'$ is
our candidate for $S_{X,2}$ but we must refine it a little.

For each $i\in S_{X,1}'$
let 
$$\multline N_i=\{j:j\in S_{Y,2}',
|(x_{\phi(i_1)}^*\otimes y_{\phi(j)}^*)(Tx_{\phi(i)}\otimes
y_{\phi(j)})|<\epsilon_1 w_{\phi(i)}\\
\text{ and } 
|(x_{\phi(i)}^*\otimes y_{\phi(j)}^*)(Tx_{\phi(i_1)}\otimes
y_{\phi(j)})|<\epsilon_1 w_{\phi(i_1)}\}.
\endmultline$$
We need to show that $\{i:N_i\supset N_i'\backepsilon j_1,j_2\in N_i',
N_i'\in \Cal S(Y)\}\cup\{i_1\}$
contains a set $S_{X,2}\in \Cal S(X)$, $i_1 \in S_{X,2}$, for
which $\phi(S_{X,2})\cap(\{\phi(i_1)^1\}\times \N)$ is infinite.
First observe that for $i$ sufficiently large $j_1,j_2\in N_i.$ Next
we will show that  for almost all $i$ there exists $N_i'\subset N_i$
such that $N_i' \in \Cal S(Y).$
By applying Lemma 6.1, for each $j$, to $T_j
:[x_{\phi(i)}\otimes y_{\phi(j)}:j\in \N]\rightarrow
[x_{\phi(i_1)}\otimes y_{\phi(j)}:i\in \{i_1\}+(\N\setminus \{i_1\})]$
defined by
$$T_i(z)=x_{\phi(i_1)}^*\otimes
y_{\phi(j)}^*(Tz)x_{\phi(i_1)}\otimes y_{\phi(j)}
+\sum_{i\neq i_1} x_{\phi(i)}^*\otimes y_{\phi(j)}^*(Tz)
x_{\phi(i)}\otimes y_{\phi(j)},$$
we have that $\sum_{i:j \notin N_i}
w^{2p/(p-2)}_{\phi(i)}<H(\epsilon_1, \{w_{\phi(i_1)}\},1).$
If for some $i$, $N_i$ contains no subset which is
in $\Cal S(Y)$, then $\phi(N_i) \cap (\{m\}\times\N)$ is finite for
all but finitely many $m.$  Lemma 6.1 implies that
$\sum w_{\phi(i)}^{2p/(p-2)}$
over such $i$ must be no more than $H(\epsilon_1, \{w_{\phi(i_1)}\},1).$
Indeed
if the sum exceeded this bound then there would be a finite set
$F$ of such $i$ such that $\sum_{i\in F}
w_{\phi(i)}^{2p/(p-2)}>H(\epsilon_1, \{w_{\phi(i_1)}\},1).$
Then there would be infinitely many $m$ for which
$\phi(N_i) \cap
(\{m\}\times\N)$ is finite for all of those $i\in F$, thus for any
$j\in
S_{Y,2}'\setminus \cup_{i\in F} N_i$,
$T_j$ would violate the conclusion of Lemma 6.1.
In particular
for those $i$ such that $\phi(i)\in (\{\phi(i_1)^1\}\times \N)$,
$w_{\phi(i)}=w_{\phi(i_1)}$ and thus only finitely many of these
$N_i$ will fail to contain a subset $N_i'$ as above. Similarly,
if
$F\subset \{i:N_i\not\supset N_i'\backepsilon j_2\in N_i',
\{j_1\}\cup N_i'\in \Cal S(Y)\text{ but } \exists N \subset N_i,
N \in S(Y)\}$ and $\sum_{i\in F}
w_{\phi(i)}^{2p/(p-2)}>2H(\epsilon_1, \{w_{\phi(i_1)}\},1),$
then $\phi(N_i) \cap
(\{\phi(j)^1\}\times\N)$ is finite  for $j=j_1\text{ or }j_2$
for all of those $i\in F$, thus for at least one of $j'=j_1,j_2$
there exists $F'\subset F$ such that  $\sum_{i\in F'}
w_{\phi(i)}^{2p/(p-2)}>H(\epsilon_1, \{w_{\phi(i_1)}\},1).$
and
$j\in
S_{Y,2}'\setminus \cup_{i\in F} N_i$ with $\phi(j)^1=\phi(j')$.
This
would violate the conclusion of Lemma 6.1 for $T_j$.
It follows that the
required set $S_{X,2}$ exists.

Player 2 selects $i_2$ from $S_{X,2}$ and selects $S_{X,2}'\subset S_{X,2}$
which is in $\Cal S(X)$ and contains $\{i_1,i_2\}.$ By our choice of
$S_{X,2}$, $N_{i_2}$ contains $N_{i_2}'$ such that $\{j_1,j_2\}\subset 
N_{i_2}'\in
\Cal S(Y).$ 
Note that for all $j\in N_{i,2}'$, and $s\neq t\in \{1,2\}$,
$$|(x_{\phi(i_s)}^*\otimes y_{\phi(j)}^*)(Tx_{\phi(i_t)}\otimes
y_{\phi(j)})|<\epsilon_{(s\vee t)-1}w_{\phi(i_t)}.$$
The set $S= N_{i_2}'$ is our candidate for $S_{Y,3}$
however we must now refine it further as we did above to produce $S_{X,2}.$
This completes the initial phase of the induction.

Suppose that we have played $r$ turns of each game to get $I=\{i_1,\dots,i_r\}$
, $J=\{j_1,\dots,j_r\},$ $J\subset S_{Y,r}'\in \Cal S(Y)$ and $I\subset
S_{X_r}'\in \Cal S(X).$

For each $j\in S_{Y,r}'\setminus J$
let 
$$\multline K_j=\{i:i\in S_{X,r}',
|(x_{\phi(i)}^*\otimes y_{\phi(j')}^*)(Tx_{\phi(i)}\otimes
y_{\phi(j)})|<\epsilon_r w_{\phi(j)}/r\\
\text{ and } 
|(x_{\phi(i)}^*\otimes y_{\phi(j)}^*)(Tx_{\phi(i)}\otimes
y_{\phi(j')})|<\epsilon_r w_{\phi(j')}/r\text{ for all }j'\in J\}
.\endmultline$$
We need to show that $\{j:K_j\supset K_j'\backepsilon
I \subset K_j'\text{ and }K_j'\in \Cal S(X)\}\cup J$
contains a set $S_{Y,r+1}\in \Cal S(Y)$, $J\subset S_{Y,r+1}$ for
which $\phi(S_{Y,r+1})\cap(\{\phi(j)^1\}\times \N)$ is infinite for all $j
\in J.$
We will show that  for almost all $j$ there exists $K_j'\subset K_j$
such that $K_j' \in \Cal S(X).$ 
For each $i$ let
operator $T_i:[x_{\phi(i)}\otimes y_{\phi(j)}:j\in \N]\rightarrow
[x_{\phi(i)}\otimes y_{\phi(j)}:j\in J+(\N\setminus J)]$, where the
notation $J+(\N\setminus J)$ means that the elements in $J$ precede the
others, defined by
$$T_i(z)=\sum_{j\in J} x_{\phi(i)}^*\otimes
y_{\phi(j)}^*(Tz)x_{\phi(i)}\otimes y_{\phi(j)}
+\sum_{j\notin J } x_{\phi(i)}^*\otimes y_{\phi(j)}^*(Tz)
x_{\phi(i)}\otimes y_{\phi(j)}.$$
By Lemma 6.1 for each $i$,
$\sum_{j:i \notin K_j}
w^{2p/(p-2)}_{\phi(j)}<
H(\epsilon_r, \min\{w_{\phi(j')}:j'\in J\},r)$  
If for some $j$, $K_j$ contains no subset which is
in $\Cal S(X)$, then $\phi(K_j) \cap (\{m\}\times\N)$ is finite for
all but finitely many $m.$  Lemma 6.1 implies that
$\sum w_{\phi(j)}^{2p/(p-2)}$
over such $j$ must be no more than
$H(\epsilon_r, \min\{w_{\phi(j')}:j'\in J\},r).$ 
Indeed
if the sum exceeded this bound then there would be a finite set
$F$ of such $j$ such that
$\sum_{j\in F} w_{\phi(j)}^{2p/(p-2)}
>H(\epsilon_r, \min\{w_{\phi(j')}:j'\in J\},r).$
Then there would be infinitely many $m$ for which $\phi(K_j) \cap
(\{m\}\times\N)$ is finite for all of those $j\in F$, thus any $i\in
S_{X,r}'\setminus \cup_{j\in F} K_j$
would violate the conclusion of Lemma 6.1.
In particular 
for those $j$ such that $\phi(j)\in (\{\phi(j')^1\}\times \N)$, for some
$j'\in J,$
$w_{\phi(j)}=w_{\phi(j')}$ and thus only finitely many of these
$K_j$ will fail to contain a subset $K_j'$ as above.
 Similarly,
if
$F\subset \{j:K_j\not\supset K_j'\backepsilon I\subset K_j',
 K_j'\in \Cal S(X)\text{ but } \exists K \subset K_j,
K \in S(Y)\}$ and
$\sum_{j\in F}
w_{\phi(j)}^{2p/(p-2)}>r H(\epsilon_r, \min\{w_{\phi(j')}:j'\in J\},r)$
then for each $i\in I$,let  $F_i=\{j\in F:\phi(K_j) \cap
(\{\phi(i)^1\}\times\N) \text{ is finite}\}.$
For at least one $i' \in I$,
$\sum_{j\in F_i}
w_{\phi(j)}^{2p/(p-2)}>H(\epsilon_r, \min\{w_{\phi(j')}:j'\in J\},r)$
thus any $i\in
S_{X,2}'\setminus \cup_{j\in F_i} K_j$ with $\phi(i)^1=\phi(i')^1$
would violate the conclusion of Lemma 6.1.
It follows that the
required set $S_{Y,r+1}$ exists. 

Player r+1 selects $j_{r+1}$ from $S_{Y,r+1}$
and selects $S_{Y,r+1}'\subset S_{Y,r+1}$
which is in $\Cal S(Y)$ and contains $J\cup\{j_{r+1}\}.$ By our choice of
$S_{Y,r+1}$, $K_{j_{r+1}}$ contains $K_{j_{r+1}}'$
such that $I\subset K_{j_{r+1}}'\in
\Cal S(Y).$
Note that for all $i\in K_{j,r+1}'$, and $s\neq t\in \{1,2,\dots,r+1\}$,
$$|(x_{\phi(i)}^*\otimes y_{\phi(j_s)}^*)(Tx_{\phi(i)}\otimes
y_{\phi(j_t)})|<\epsilon_{((s\vee t))-1}w_{\phi(j_t)}.$$
The set $S= K_{j_{r+1}}'$ is
our candidate for $S_{X,r+1}$ but we must again refine it a little.

For each $i\in S\setminus I$
let 
$$\multline N_i=\{k:k\in S_{Y,r+1}',
|x_{\phi(i')}^*\otimes y_{\phi(k)}^*(Tx_{\phi(i)}\otimes
y_{\phi(k)})|<\epsilon_r w_{\phi(i)}/r\\
\text{ and }.
|(x_{\phi(i)}^*\otimes y_{\phi(j)}^*)(Tx_{\phi(i')}\otimes
y_{\phi(j)})|<\epsilon_r w_{\phi(i')}/r\text{ for all }i'\in I\}
\endmultline$$
We need to show that $\{i:N_i\supset N_i'\backepsilon J\subset N_i',
 N_i'\in \Cal S(Y)\}\cup I$
contains a set $S_{X,r+1}\in \Cal S(X)$, $I\subset S_{X,r+1}$, for
which $\phi(S_{X,r+1})\cap(\{\phi(i)^1\}\times \N)$ is infinite for all
$i\in I.$
First observe that for $i$ sufficiently large $J\subset N_i.$ Next
we will show that  for almost all $i$ there exists $N_i'\subset N_i$
such that $N_i' \in \Cal S(Y).$
By applying Lemma 6.1 for each $j$, as above,
$\sum_{i:j \notin N_i}
w^{2p/(p-2)}_{\phi(i)}<H(\epsilon_r, \min\{w_{\phi(i')}:i'\in I\},r).$
If for some $i$, $N_i$ contains no subset which is
in $\Cal S(Y)$, then $\phi(N_i) \cap (\{m\}\times\N)$ is finite for
all but finitely many $m.$  Lemma 6.1 implies that
$\sum w_{\phi(i)}^{2p/(p-2)}$
over such $i$ must be no more than
$H(\epsilon_r, \min\{w_{\phi(i')}:i'\in I\},r).$  
Indeed
if the sum exceeded this bound then there would be a finite set
$F$ of such $i$ such that $\sum_{i\in F}
w_{\phi(i)}^{2p/(p-2)}>H(\epsilon_r, \min\{w_{\phi(i')}:i'\in I\},r).$
Then there would be infinitely many $m$ for which
$\phi(N_i) \cap
(\{m\}\times\N)$ is finite for all of those $i\in F$, thus any
$j\in
S_{Y,r+1}'\setminus \cup_{i\in F} N_i$
would violate the conclusion of Lemma 6.1.
In particular
for those $i$ such that $\phi(i)\in (\{\phi(i')^1\}\times \N)$, $i'\in I$,
$w_{\phi(i)}=w_{\phi(i')}$ and thus only finitely many of these
$N_i$ will fail to contain a subset $N_i'$ as above. Similarly,
if
$F\subset \{i:N_i\not\supset N_i'\backepsilon J\cup\{j_{r+1}\}\subset N_i',
 N_i'\in \Cal S(Y)\text{ but } \exists N \subset N_i,
N \in S(Y)\}$ and $\sum_{i\in F}
w_{\phi(i)}^{2p/(p-2)}>r H(\epsilon_r, \min\{w_{\phi(i')}:i'\in I\},r).$
then $\phi(N_i) \cap
(\{\phi(j)^1\}\times\N)$ is finite  for some $j\in J\cup\{j_{r+1}\}$
for all of those $i\in F$, thus for at least one of $j'=j_1,j_2,\dots,
j_{r+1}$
there exists $F'\subset F$ such that  $\sum_{i\in F'}
w_{\phi(i)}^{2p/(p-2)}>H(\epsilon_r, \min\{w_{\phi(i')}:i'\in I\},r)$ 
and
$j\in
S_{Y,r+1}'\setminus \cup_{i\in F} N_i$ with $\phi(j)^1=\phi(j')$.
This
would violate the conclusion of Lemma 6.1.
It follows that the
required set $S_{X,r+1}$ exists.

Player r+1 selects $i_{r+1}$ from $S_{X,r+1}$
and selects $S_{X,r+1}'\subset S_{X,r+1}$
which is in $\Cal S(X)$ and contains $\{i_1,i_2,\dots,i_{r+1}\}.$
By our choice of
$S_{X,r+1}$, $N_{i_{r+1}}$ contains $N_{i_{r+1}}'$ such that
$\{j_1,\dots,j_{r+1}\}\subset N_{i_{r+1}}'\in
\Cal S(Y).$ The set $S= N_{i_{r+1}}'$
is our candidate for $S_{Y,r+2}$.

This completes the induction step.

Therefore there are rich subsets $I=\{i_1,i_2,\dots\},J=\{j_1,j_2,\dots\}$
such that 
\item{(1)} $|(x_{\phi(i_s)}^*\otimes y_{\phi(j)}^*)
(Tx_{\phi(i_t)}\otimes y_{\phi(j)})|
<\epsilon_{(s\vee t)-1} w_{\phi(i_t)}/((s\vee t)-1)$ for all $
s\neq t \in \N$ and $j\in J.$ 
\item{(2)} $|(x_{\phi(i)}^*\otimes y_{\phi(j_s)}^*)
(Tx_{\phi(i)}\otimes y_{\phi(j_t)})|
<\epsilon_{(s\vee t)-1} w_{\phi(j_t)}/((s\vee t) -1)$ for all $i\in I$ and $
s\neq t \in \N$
\item{(3)} $(x_{\phi(i)}^*\otimes y_{\phi(j)}^*)
(Tx_{\phi(i')}\otimes y_{\phi(j')})=0
$ for all $i'\in I$ and $j'\in J$ such that $i\neq i'$ and $j\neq j'$;
$i=i'$ and $o(j\vee j')>o(i),j\neq j'$;
or $i\neq i',o(i\vee i')\geq o(j)$ and $j=j'.$ 

It remains to show that the restriction of $T$ to $Z=[x_{\phi(i)}\otimes
y_{\phi(j)}:i\in I, j\in J]$ is an isomorphism. As noted above we may
assume by composing with the projection onto $Z$
that $T$ is actually a map from $Z$ to itself and that $\epsilon
=1.$ Moreover, we can assume that $I=J=\N$ by replacing $i$ by $o(i,I)$ and
$j$ by $o(j,J).$ 
Thus  with the notation appropriately revised we have
\item{(0)} $( x_{\phi(i)}^*\otimes y_{\phi(j)}^*)
(Tx_{\phi(i)}\otimes y_{\phi(j)})=1$ for all $i,j$,
\item{(1)} $|(x_{\phi(i)}^*\otimes y_{\phi(j)}^*)
(Tx_{\phi(i')}\otimes y_{\phi(j)})|
<\epsilon_{i\vee i'-1} w_{\phi(i')}/(i\vee i'-1)$ for all $
i\neq i' \in \N$ and $j\in J.$ 
\item{(2)} $|(x_{\phi(i)}^*\otimes y_{\phi(j)}^*)
(Tx_{\phi(i)}\otimes y_{\phi(j')})|
<\epsilon_{j\vee j'-1} w_{\phi(j')}/(j\vee j'-1)$ for all $i\in I$ and $
j\neq j'\in \N$
\item{(3)} $(x_{\phi(i)}^*\otimes y_{\phi(j)}^*)
(Tx_{\phi(i')}\otimes y_{\phi(j')})=0
$ for all $i'\in I$ and $j'\in J$ such that $i\neq i'$ and $j\neq j'$;
$i=i'$ and $j\vee j'>i,j\neq j'$;
or $i\neq i',i\vee i'\geq j$ and $j=j'.$ 

We will now estimate $\|Tz-z\|$ for $z\in Z$ with finite support relative
to the basis.

Let $z=\sum_{i\in I}\sum_{j\in J} 
a_{i,j} x_{\phi(i)}\otimes y_{\phi(j)})$.

First we estimate the $\ell_2$-norm.
$$\align
\|Tz&-z\|_2^2\\
&= \sum_{i\in \N}\sum_{j\in \N}
\left|\sum\Sb i',j'\in \N\\(i',j')\neq (i,j)\endSb a_{i',j'}
(x_{\phi(i)}^*\otimes y_{\phi(j)}^*)
(Tx_{\phi(i')}\otimes y_{\phi(j')})\right|^2
w_{\phi(i)}^2 w_{\phi(j)}^2\\
&=
\sum_{j=1}^\infty \sum_{i< j}
\left|\sum\Sb i',j'\in \N\\(i',j')\neq (i,j)\endSb a_{i',j'}
(x_{\phi(i)}^*\otimes y_{\phi(j)}^*)
(Tx_{\phi(i')}\otimes y_{\phi(j')})\right|^2
w_{\phi(i)}^2 w_{\phi(j)}^2\\
&+
\sum_{i=1}^\infty \sum_{j\leq i}
\left|\sum\Sb i',j'\in \N\\(i',j')\neq (i,j)\endSb a_{i',j'}
(x_{\phi(i)}^*\otimes y_{\phi(j)}^*)
(Tx_{\phi(i')}\otimes y_{\phi(j')})\right|^2
w_{\phi(i)}^2 w_{\phi(j)}^2\\
&=
\sum_{j=1}^\infty \sum_{i<j}
\left|\sum_{1\leq i'\neq i <j}a_{i',j}
(x_{\phi(i)}^*\otimes y_{\phi(j)}^*)
(Tx_{\phi(i')}\otimes y_{\phi(j)})\right|^2
w_{\phi(i)}^2 w_{\phi(j)}^2\\
&+
\sum_{i=1}^\infty \sum_{j\leq i}
\left|\sum_{1\leq j'\neq j \leq i}a_{i',j'}
(x_{\phi(i)}^*\otimes y_{\phi(j)}^*)
(Tx_{\phi(i)}\otimes y_{\phi(j')})\right|^2
w_{\phi(i)}^2 w_{\phi(j)}^2\\
\intertext{ by (3)}
&\leq
\sum_{j=1}^\infty \sum_{i<j}
\left(\sum_{1\leq i'< i }|a_{i',j}|
\frac{\epsilon_{i-1}w_{\phi(i')}}{i-1}+\sum_{i':i<i'< j}
|a_{i',j}|\frac{\epsilon_{i'-1}w_{\phi(i')}}{i'-1}\right)^2
w_{\phi(i)}^2 w_{\phi(j)}^2\\
&+\sum_{i=1}^\infty \sum_{j\leq i}
\left(\sum_{1\leq j'< j }|a_{i,j'}|
\frac{\epsilon_{j-1}w_{\phi(j')}}{j-1}j+\sum_{j<j'\leq i}
|a_{i,j'}|\frac{\epsilon_{j'-1}w_{\phi(j')}}{j'-1}\right)^2
w_{\phi(i)}^2 w_{\phi(j)}^2\\
\intertext{ by (1) and (2)}
&\leq
\sum_{j=1}^\infty \sum_{i<j}
\left(\sum_{1\leq i'< i }|a_{i',j}|^2
\frac{\epsilon_{i-1}w_{\phi(i')}^2}{i-1}+\sum_{i':i<i'< j}
|a_{i',j}|^2\frac{\epsilon_{i'-1}w_{\phi(i')}^2}{i'-1}\right)
w_{\phi(i)}^2 w_{\phi(j)}^2\\
&+\sum_{i=1}^\infty \sum_{j\leq i}
\left(\sum_{1\leq j'< j }|a_{i,j'}|^2
\frac{\epsilon_{j-1}w_{\phi(j')}^2}{j-1}+\sum_{j<j'\leq i}
|a_{i,j'}|^2\frac{\epsilon_{j'-1}w_{\phi(j')}^2}{j'-1}\right)
w_{\phi(i)}^2 w_{\phi(j)}^2\\
\intertext{ by the convexity of $t^2$}
&\leq
\sum_{j=1}^\infty \sum_{i'=1}^{j-1} |a_{i',j}|^2w_{\phi(i')}^2
w_{\phi(j)}^2(\sum_{i:j>i>i'}  \frac{w_{\phi(i)}^2\epsilon_{i-1}}{i-1}
+\sum_{i:i<i'}\frac{\epsilon_{i'-1}w_{\phi(i)}^2}{i'-1})\\
&+\sum_{i=1}^\infty \sum_{j'=1}^{i} |a_{i,j'}|^2w_{\phi(i)}^2
w_{\phi(j')}^2(\sum_{j:i\geq j>j'} \frac{w_{\phi(j)}^2\epsilon_{j-1}}{j-1} 
+\sum_{j:j<j'}\frac{\epsilon_{j'-1}w_{\phi(j)}^2}{j'-1})
\\
&\leq \|z\|_2^2\epsilon'/4^p
\endalign$$

Next  we estimate the $\ell_p$-norm.
$$\align
|Tz-z|_p^p&=
\sum_{i\in \N}\sum_{j\in \N}
\left|\sum\Sb i',j'\in \N\\(i',j')\neq (i,j)\endSb a_{i',j'}
(x_{\phi(i)}^*\otimes y_{\phi(j)}^*)
(Tx_{\phi(i')}\otimes y_{\phi(j')})\right|^p \\
&=
\sum_{j=1}^\infty \sum_{i<j}
\left|\sum_{1\leq i'\neq i<j }a_{i',j}
(x_{\phi(i)}^*\otimes y_{\phi(j)}^*)
(Tx_{\phi(i')}\otimes y_{\phi(j)})\right|^p\\
&+
\sum_{i=1}^\infty \sum_{j\leq i}
\left|\sum_{1\leq j'\neq j \leq i}a_{i',j'}
(x_{\phi(i)}^*\otimes y_{\phi(j)}^*)
(Tx_{\phi(i)}\otimes y_{\phi(j')})\right|^p\\
\intertext{(by (3))}
&\leq
\sum_{j=1}^\infty \sum_{i<j}
\left(\sum_{1\leq i'< i }|a_{i',j}|
\frac{\epsilon_{i-1}w_{\phi(i')}}{i-1}+\sum_{i':i<i'< j}
|a_{i',j}|\frac{\epsilon_{i'-1}w_{\phi(i')}}{i'-1}\right)^p\\
&+\sum_{i=1}^\infty \sum_{j\leq i}
\left(\sum_{1\leq j'< j }|a_{i,j'}|
\frac{\epsilon_{j-1}w_{\phi(j')}}{j-1}+\sum_{j<j'\leq i}
|a_{i,j'}|\frac{\epsilon_{j'-1}w_{\phi(j')}}{j'-1}\right)^p\\
\intertext{(by (1) and (2))}
&\leq
\sum_{j=1}^\infty \sum_{i<j}
\left(\sum_{1\leq i'< i }|a_{i',j}|^p
\frac{\epsilon_{i}w_{\phi(i')}}{i-1}+\sum_{i':i<i'< j}
|a_{i',j}|^p\frac{\epsilon_{i'-1}w_{\phi(i')}}{i'-1}\right)\\
&+\sum_{i=1}^\infty \sum_{j\leq i}
\left(\sum_{1\leq j'< j }|a_{i,j'}|^p
\frac{\epsilon_{j-1}w_{\phi(j')}}{j-1}+\sum_{j<j'\leq i}
|a_{i,j'}|^p\frac{\epsilon_{j'-1}w_{\phi(j')}}{j'-1}\right)\\
\intertext{(by the convexity of $t^p$)}
&\leq
\sum_{j=1}^\infty \sum_{i':1\leq i'<j} |a_{i',j}|^p
(\sum_{i:i>i'}  \frac{w_{\phi(i)}\epsilon_{i-1}}{i-1}
+\sum_{i:i<i'}\frac{\epsilon_{i'-1}w_{\phi(i)}}{i'-1})\\
&+\sum_{i=1}^\infty \sum_{j':1\leq j'<i} |a_{i,j'}|^p 
(\sum_{j:j>j'} \frac{w_{\phi(j)}\epsilon_{j-1}}{j-1}
+\sum_{j:j<j'}\frac{\epsilon_{j'-1}w_{\phi(j)}}{j'-1})
\\
&\leq |z|_p^p\epsilon'/4^p
\endalign$$

The estimates for the two row and column norms are similar so we will only do
one.
$$\align
\|T&z-z\|_{R}^p\\
&= \sum_{i\in \N}\left(\sum_{j\in \N}
\left|\sum\Sb i',j'\in \N\\(i',j')\neq (i,j)\endSb a_{i',j'}
(x_{\phi(i)}^*\otimes y_{\phi(j)}^*)
(Tx_{\phi(i')}\otimes y_{\phi(j')})\right|^2
w_{\phi(j)}^2\right)^{\frac{p}{2}}\\
&\leq
2^{p-1}\left(\sum_{i=1}^\infty\left( \sum_{j> i}
\left|\sum\Sb i',j'\in \N\\(i',j')\neq (i,j)\endSb a_{i',j'}
(x_{\phi(i)}^*\otimes y_{\phi(j)}^*)
(Tx_{\phi(i')}\otimes y_{\phi(j')})\right|^2
 w_{\phi(j)}^2\right)^{\frac{p}{2}}\right.\\
&+\left.
\sum_{i=1}^\infty \left(\sum_{j\leq i}
\left|\sum\Sb i',j'\in \N\\(i',j')\neq (i,j)\endSb a_{i',j'}
(x_{\phi(i)}^*\otimes y_{\phi(j)}^*)
(Tx_{\phi(i')}\otimes y_{\phi(j')})\right|^2
 w_{\phi(j)}^2\right)^{\frac{p}{2}}\right)
\\
&=
2^{p-1}\left(\sum_{i=1}^\infty \left(\sum_{j>i}
\left|\sum_{1\leq i'\neq i < j }a_{i',j}
(x_{\phi(i)}^*\otimes y_{\phi(j)}^*)
(Tx_{\phi(i')}\otimes y_{\phi(j)})\right|^2
 w_{\phi(j)}^2\right)^{\frac{p}{2}}\right.\\
&+\left.
\sum_{i=1}^\infty \left(\sum_{j\leq i}
\left|\sum_{1\leq j'\neq j \leq i}a_{i',j'}
(x_{\phi(i)}^*\otimes y_{\phi(j)}^*)
(Tx_{\phi(i)}\otimes y_{\phi(j')})\right|^2
w_{\phi(j)}^2\right)^{\frac{p}{2}}\right)\\
\intertext{(by (3))}
&\leq
2^{p-1}\left(\sum_{i=1}^\infty \left(\sum_{j>i}
\left(\sum_{i'=1}^{i-1}|a_{i',j}|
\frac{\epsilon_{i-1}w_{\phi(i')}}{(i-1)}+\sum_{i'=i+1}^{j-1}
|a_{i',j}|\frac{\epsilon_{i'-1}w_{\phi(i')}}{(i'-1)}\right)^2
w_{\phi(j)}^2\right)^{\frac{p}{2}}\right.\\
&+\left.\sum_{i=1}^\infty\left( \sum_{j\leq i}
\left(\sum_{j'=1}^{j-1}|a_{i,j'}|
\frac{\epsilon_{j-1}w_{\phi(j')}}{j-1}+\sum_{j'=j+1}^{i}
|a_{i,j'}|\frac{\epsilon_{j'-1}w_{\phi(j')}}{(j'-1)}\right)^2
w_{\phi(j)}^2\right)^{\frac{p}{2}}\right)\\
\intertext{(by (1) and (2))}
&\leq
2^{p-1}\left(\sum_{i=1}^\infty\left( \sum_{j>i}
\left(\sum_{i'=1}^{i-1}|a_{i',j}|^2
\frac{\epsilon_{i-1}w_{\phi(i')}^2}{(i-1)}+\sum_{i'=i+1}^{j-1}
|a_{i',j}|^2\frac{\epsilon_{i'-1}w_{\phi(i')}^2}{i'-1}\right)
w_{\phi(j)}^2\right)^{\frac{p}{2}}\right.\\
&+\left.\sum_{i=1}^\infty \left(\sum_{j\leq i}
\left(\sum_{j'=1}^{j-1}|a_{i,j'}|^2
\frac{\epsilon_{j-1}w_{\phi(j')}^2}{j-1}+\sum_{j<j'\leq i}
|a_{i,j'}|^2\frac{\epsilon_{j'-1}w_{\phi(j')}^2}{j'-1}\right)
w_{\phi(j)}^2\right)^{\frac{p}{2}}\right)
\\
\intertext{( by the convexity of $t^2$)}
&\leq
2^{p-1}\left(\sum_{i=1}^\infty 
\left(\sum_{i':1\leq i'<i} \left(\sum_{j>i}|a_{i',j}|^2
w_{\phi(j)}^2\right)^{\frac{p}{2}}
\frac{w_{\phi(i')}^2\epsilon_{i-1}}{i-1}\right. \right .\\
&+\left. \left.
\sum_{i'>i}\left(\sum_{j>i'}|a_{i',j}|^2
w_{\phi(j)}^2\right)^{\frac{p}{2}}
\frac{\epsilon_{i'-1}w_{\phi(i)}^2}{i'-1}\right
)\right.\\
&+\left.\sum_{i=1}^\infty \left(\sum_{j':1\leq j'\leq i} |a_{i,j'}|^2
w_{\phi(j')}^2(\sum_{j:j>j'} \frac{w_{\phi(j)}^2\epsilon_{j-1}}{j-1}
+\sum_{j:j<j'}\frac{\epsilon_{j'-1}w_{\phi(j)}^2}{j'-1})\right)^{\frac{p}{2}}
\right)
\\
\intertext{(by the convexity of $t^{\frac{p}{2}}$)}
&\leq
2^{p-1}\left(\sum_{i=1}^\infty
\left(\sum_{j>i}|a_{i',j}|^2
w_{\phi(j)}^2\right)^{\frac{p}{2}}\frac{\epsilon'}{4^p}\right.\\
&+\left.\sum_{i=1}^\infty \left(\sum_{j':1\leq j'\leq i}
|a_{i,j'}|^2
w_{\phi(j')}^2\right)^{\frac{p}{2}}
(\frac{\epsilon'}{4^p})^{\frac{p}{2}}\right)\\
&\leq \|z\|_{R}^p\frac{2^p\epsilon'}{4^p}
\endalign$$

Because $\epsilon'< 1$, we have that there is an $\epsilon'<1$
such that $\|Tz-z\|<\epsilon'\|z\|$ for all
$z\in Z.$ Thus $T|_Z$ is an isomorphism.
\qed
\enddemo

Notice that in fact we can make the constant
$\epsilon'$ in the proof above as close to 0 as we wish.

Proposition 6.3 shows that we can focus on the coordinates of an operator
on $X_p\otimes X_p$ in attempting to show that an operator is large. We
will use this in Section 8 to show that certain operators do not exist. For
further development of the isomorphic theory it is also necessary to know
that there are many natural operators on a space. We will finish this
section by looking at some of the natural classes of operators and some
interesting subsets of the natural bases of $X_p\otimes X_p.$

\definition{Definition 6.4} If $T$ is an operator from $X_{p,(w_n)}
\otimes X_{p,(w_n')}$ into $X_{p,(u_n)} \otimes X_{p,(u_n')}$, say
that $T$ is $(p,R,C,2)-$\tiny bounded if it is bounded with respect to each
of the four norms, i.e., there is a constant $K$ such that for all $z \in
X_{p,(w_n)} \otimes X_{p,(w_n')}$,
$\|Tz\|_p\leq K \|z\|_p,\|Tz\|_R\leq K \|z\|_R,\|Tz\|_C\leq K \|z\|_C$ and
$\|Tz\|_2\leq K \|z\|_2.$
\enddefinition

Observe that if $T_1$ is a $(p,2)-$bounded operator from $X_{p,(w_n)}$ into 
$X_{p,(u_n)}$ and $T_2$ is a $(p,2)-$bounded operator from $X_{p,(w_n')}$
into $X_{p,(u_n')}$, then $T_1\otimes T_2$ is $(p,R,C,2)-$\tiny bounded.
Also any basis projection relative to the natural basis of $X_{p,(w_n)}
\otimes X_{p,(w_n')}$ is $(p,R,C,2)-$\tiny bounded. Thus there are many
such operators.

\proclaim{Lemma 6.5} If $(F_i)$ and $(G_i)$ are sequences of finite subsets
of $\N$ such that $\max F_i <\min F_{i+1}$ and $\max F_i <\min F_{i+1}$
then 
$$[e_n\otimes e_m':(n,m)\in F_i\times G_i\text{ for some $i$}],$$
where $(e_n)$ and $(e_m')$ are natural bases for $X_{p,(w_n)}$ and
$X_{p,(w_n')}$, respectively, is isomorphic to a complemented subspace of
$X_p.$
\endproclaim

\demo{Proof} Observe that the block subspaces $Z_i=[e_n\otimes
e_m':(n,m)\in F_i\times G_i]$ are $(p,2)-$\tiny complemented in
$X_{p,(w_n)} \otimes X_{p,(w_n')}$ and that $Z=\sum_{i\in \N} Z_i$ is a natural
$(p,2)-$\tiny sum decomposition. Therefore, $Z$ is $(p,2)-$\tiny isomorphic
to a $(p,2)-$\tiny complemented subspace of $\rp{\omega}$ and thus of
$X_p.$
\enddemo

It follows from Lemma 6.5 and the fact that if $Y$ is a complemented
subspace of $X_p$ then $Y\oplus X_p$ is isomorphic to $X_p$, \cite{JO}, that 
$X_{p,(w_n)} \otimes X_{p,(w_n')}$ is isomorphic to $[e_n\otimes e_m':n\neq
m].$  Thus the diagonal is a negligible subspace of $X_p\otimes X_p.$ We
will next show that the upper (lower) triangle of $X_p\otimes X_p$ is
isomorphic to the whole space. Notice that this question depends on the
representation of $X_p$. There are some representations where the argument
follows from a Cantor-Bernstein result \cite{W2}, \cite{Woj1}.

\proclaim{Lemma 6.6} Suppose that $(w_n)$ and $(w_n')$ are
sequences in $(0,1]$ which
satisfy (*)  and that $w_{n,k}=w_n$ and $w_{n,k}'=w_n'$ for all $n\in \N.$ Then
$X_{p,(w_{n,k})}\otimes X_{p,(w_{n,k}')}$ is isomorphic to $[e_{\phi(i)}\otimes
e_{\phi(j)}:i<j],$ where $\phi:\N\rightarrow \NxN$ is a bijection as in
Section 5.
\endproclaim
\demo{Proof} It is sufficient to show that $(e_{\phi(i)}\otimes
e_{\phi(j)})_{i< j}$ contains a subsequence which is equivalent to the
natural basis of $X_{p,(w_{n,k})}\otimes X_{p,(w_{n,k}')}$. Observe that the
reflection along the diagonal mapping, $R$, defined by $R(e_{\phi(i)}\otimes
e_{\phi(j)})=e_{\phi(j)}\otimes e_{\phi(i)})$ extends 
linearly to an isomorphism. Because every weight is repeated
infinitely often, we  can split the columns with weight $w_n$ into two
infinite sets with indices $K_n$ and $L_n$, i.e., for $k \in K_n\cup L_n$,
$\|e_{\phi(k)}\|_2=w_n$, $K_n\cap L_n=\emptyset$ and $|K_n|=|L_n|=\infty.$
Let $\psi_n:K_n\cup L_n\rightarrow K_n$ be injective and satisfy
$\psi_n(k)> k$ for all $k$ and let $\zeta_n:K_n\cup L_n\rightarrow L_n$
do likewise. Define  $\psi:\N \rightarrow \cup_n K_n$ by $\psi(i)=\psi_n(i)$
if $i\in K_n\cup L_n$ and $\zeta:\N \rightarrow \cup_n L_n$ by
$\zeta(i)=\zeta_n(i)$ if  $i\in K_n\cup L_n$.
Define 
$$S(e_{\phi(i)}\otimes e_{\phi(j)})=\cases e_{\phi(i)}\otimes
e_{\phi(\psi(j))}, &\text{if $i\leq j$}\\
e_{\phi(j)}\otimes e_{\phi(\zeta(i))},
 &\text{if $j<i$},\endcases$$ 
and extend linearly. It is easy to check that $S$ is an isomorphism.
\qed\enddemo

\proclaim{Lemma 6.7} Suppose that $(w_n)$ and $(w_n')$ are
sequences in $(0,1]$ which
satisfy (*)  and that $w_{n,k}=w_n$ and $w_{n,k}'=w_n'$ for all
$n,k\in \N.$  Let $(u_n)$ and
$(u_n')$ be two sequences in $(0,1]$ which satisfy (*). Let $(e_{n,k})$,
$(e_{n,k}')$, $(d_n)$ and $(d_n')$ be the natural bases for $X_{p,(w_{n,k})}$,
$X_{p,(w_{n,k}')}$, $X_{p,(u_n)}$ and $X_{p,(u_n')}$, respectively.
Then there is a $(p,R,C,2)-$\tiny complemented subspace of $[d_n\otimes
d_m':n<m]$ which is $(p,R,C,2)-$\tiny isomorphic to $[e_{\phi(i)}\otimes
e_{\phi(j)}:i<j].$
\endproclaim
\demo{Proof} By \cite{R} there is a block basic sequence $(D_n)$
 of the basis, $(d_n)$, of $X_{p,(u_n)}$ 
which is $(p,2)-$\tiny equivalent to the basis $(e_{n,k})$ with
$(p,2)-$\tiny complemented closed span. Similarly, there is a block basic
sequence $(D_n')$  of the basis, $(d_n')$, of $X_{p,(u_n')}$
which is $(p,2)-$\tiny equivalent to the basis $(e_{n,k}')$ with
$(p,2)-$\tiny complemented closed span. Let $P$ and $P'$ be the
corresponding projections. Then $P\otimes P'$ is a $(p,R,C,2)-$\tiny
bounded projection from $X_{p,(u_n)}\otimes X_{p,(u_n')}$ onto the subspace
$[D_n\otimes D_m'].$ If we restrict this map to $[d_n\otimes d_m':n<m]$ and
compose with a basis projection to eliminate any partial support of blocks
$D_n\otimes D_n'$, we get the required projection and subspace.
\qed\enddemo

\proclaim{Proposition 6.8} Let $(u_n)$ and
$(u_n')$ be two sequences in $(0,1]$ which satisfy (*) and let $(d_n)$ and
$(d_n')$ be the natural bases for $X_{p,(u_n)}$
and $X_{p,(u_n')}$, respectively. Then
$[d_n\otimes d_m':n<m]$ is isomorphic to $X_p\otimes X_p.$
\endproclaim
\demo{Proof} It follows from Lemma 6.7 with $w_n=u_n$ and $w_n'=u_n'$
and an argument analogous to that of Rosenthal \cite{R, Theorem 13} 
that $[d_n\otimes d_m':n<m]$ is
isomorphic to its square and to $[e_{\phi(i)}\otimes e_{\phi(j)}':i<j],$
with the same notation as in the previous lemma.
By Lemma 6.6 $[e_{\phi(i)}\otimes e_{\phi(j)}:i<j]$ is isomorphic to
$X_p\otimes X_p.$
\qed\enddemo

\newpage

\head 7. Isomorphisms of $X_p\otimes X_p$ onto complemented subspaces of 
$(p,2)-$\tiny sums \endhead

Our goal is to show that $X_p\otimes X_p$ is not isomorphic to a
complemented subspace of any $R_p^\alpha,\alpha<\omega_1.$ Because
$X_p\otimes X_p$ is isomorphic to a subspace of $R_p^{\omega2}$ the
projection must play a fundamental role. We have already seen that
$\rp{\alpha}$ is isomorphic to a $(p,2)-$\tiny sum of spaces
$\rp{\alpha_n}$ In this section we will develop a
number of technical results which describe the restrictions on a
complemented isomorphic embedding of $X_p\otimes X_p$ into a
$(p,2)-$\tiny sum of subspaces of $L_p.$

Below we will be working with a $(p,2)-$\tiny sum of spaces $Y_j$ and
we will denote the natural projection onto $[Y_j:j\leq k]$ by
$P_k.$ We will use the sequence space norm rather than the norm as an
independent sum. In applications to $\rp{\alpha}$
we will need to change to the independent sum
and this will introduce a constant $C$, the constant in Rosenthal's
inequality, from the equivalence between the norm
on $(\sum Y_j)_{p,2,(w_n)}$ and
and embedding $Y_n$ into $L_p$ on independent coordinates
and using the $L_p$-norm. This change is only an annoying technicality.

The proof of the next lemma and the one following are somewhat easier if we
assume that each space $Y_n$ has an unconditional basis and that we have done a
quasi-blocking of the image of the basis of $X_p\otimes X_p$
relative to the unconditional basis of $(\sum Y_n)_{p,2}$ as
in the conclusion of Proposition 5.4. In that case we can arrange things
so that some of the error estimates can be replaced by zero
and get a stronger conclusion. We
summarize this as Lemma 7.4 below, but the reader may find the proofs of the
first two lemmas are easier to understand on the first reading if he
considers this easier case.

The first lemma is closely related to the results of Section 4 except that
we assume that we have a projection onto the range of the operator.
Throughout this section and the next we will assume that the operator from
$X_p\otimes X_p$ satisfies (T2). There is no loss of generality since by
Proposition 4.2 we can always restrict the operator to a natural
complemented subspace which is isomorphic to the whole space and so that
the restricted operator satisfies (T2).

\proclaim{Lemma 7.1}
Let $T$ be an isomorphism from $X_{p,(w_n)} \otimes X_{p,(w_n)}$ into $(\sum
Y_j)_{p,2,(w_n')}$ satisfying (T2)
and  let $P$ be a projection onto the range of $T$.
Let $\rho$, $\rho'$ and $\delta $ be positive constants with
$\rho>\rho'$
and let $(x_i)$ and $(y_j)$ be two copies of
the usual basis for $X_{p,(w_n)}$. Suppose
that $j$ is such that $w_j<\frac{\delta^{1/2}
(3\rho+\rho') }{6\sqrt{2} \|T\|\|T^{-1}\|\|P\|} $, $(N_k)_{k=1}^K$,
and $(N_k')_{k=1}^K$
are strictly increasing sequences of integers with $N_k\leq
N_k'<N_{k+1}$ for all $k$, $(F_k)_{k=1}^K$ is a
disjoint sequence of finite subsets of $\N$ such that
for all $i\in F_k$
$$\align 
|x_i^*\otimes
y_j^*(T^{-1}P(P_{N_k}-P_{N_{k-1}'})T(x_i\otimes y_j))| &\geq \rho,\tag{7.1.1}\\
|x_i^*\otimes
y_j^*(T^{-1}P(P_{N_k}-P_{N_{k-1}'})\sum \Sb s\neq i \\ s\in F_k \endSb
w_s^{2/(p-2)}T(x_s\otimes y_j))| 
&\leq (\rho-\rho') w_i^{2/(p-2)}/4\tag{7.1.2}\\
\intertext{
and for all $k$}
1\geq \sum_{i\in F_k} w_i^{2p/(p-2)} &\geq \delta \tag{7.1.h3}
\endalign$$
Then $K\leq \max\{2\delta^{-1},(\|T\|^2\|T^{-1}\| \|P\|   8/(\delta^{1/2}
(3\rho+\rho')))
^{2p/(p-2)}\}.$
\endproclaim

\demo{Proof}
First by replacing $\delta$ by $\delta/2$ we may assume that for all
$i\in F_k$,
$$x_i^*\otimes
y_j^*(T^{-1}P(P_{N_k}-P_{N_{k-1}'})T(x_i\otimes y_j)) \geq \rho,$$
that is, they all have the same sign which we assume to be positive.
Because $(r_i Tx_i\otimes y_j)_{i\in F_k}$, where $(r_i)$ is the
sequence of Rademacher, is orthogonal in $L_2(\Omega \times
[0,1])$, there is a choice of signs $(\epsilon_i)$ such that
$$\|\sum_{i\in F_k} \epsilon_i  w_i^{2/(p-2)}T x_i\otimes y_j\|_2
\leq \frac{3}{2} (\sum_{i\in F_k} w_i^{4/(p-2)} \|T x_i\otimes
y_j\|_2^2)^{1/2}.\tag{7.1.4}$$
Let $z_k=\sum_{i\in F_k} \epsilon_i  w_i^{2/(p-2)} x_i\otimes
y_j$ for each $k.$
Since $(x_i\otimes y_j)_i$ is equivalent to the usual basis of
$X_{p,(w_i)}$ it follows from Proposition 0.2 that $[z_k]$ is
complemented in $[x_i\otimes y_j:i\in \N]$ with projection $Q_1.$
This in turn induces an operator $Q$ from $(\sum
Y_j)_{p,2,(w_n)}$ onto $[z_k]$, namely, $Q=Q_1T^{-1}P.$ Since
we will need to do some computations with $Q$, let us be more
explicit about its evaluation.

Let
$$z_k^*=(\sum_{i\in F_k} w_i^{2p/(p-2)})^{-1}\sum_{i\in F_k} \epsilon_i 
w_i^{2(p-1)/(p-2)} x_i^*\otimes
y_j^*.$$
Then $Q(z)=\sum_{k=1}^K z_k^*(T^{-1}Pz)z_k.$

Let $R_k=(P_{N_k}-P_{N_{k-1}'})$ for all $k.$
Now we want to restrict $Q$ to the subspace
$Z=[R_kT(z_k)]$ and pass to its
``diagonal''. Observe that
$$\align
(\sum_{i\in F_k} w_i^{2p/(p-2)})&|z_k^*(T^{-1}PR_kT(z_k))|\\
&=
|\sum_{i\in F_k} [w_i^{2p/(p-2)}(x_i^*\otimes
y_j^*)(T^{-1}PR_kT(x_i\otimes y_j))\\
&\quad +w_i^{2(p-1)/(p-2)}
(x_i^*\otimes
y_j^*)(T^{-1}PR_k\sum \Sb s\neq i \\ s \in F_k \endSb w_s^{2/(p-2)}
T(x_s\otimes y_j))]|\\
&\geq 
\sum_{i\in F_k} w_i^{2p/(p-2)}\rho- w_i^{2(p-1)/(p-2)}w_i^{2/(p-2)}
(\rho-\rho')/4 \\
\intertext{(by 7.1.2)}
&= \sum_{i\in F_k} w_i^{2p/(p-2)}(3\rho+\rho')/4. \tag{7.1.5}
\endalign$$
Therefore $\|R_kTz_k\|\geq (3\rho+\rho')/(4\|T\|\|P\|).$
Thus we have that $(R_kT(z_k))$ is a sequence of non-zero blocks in
an unconditional sum  with constant 1 and thus is an
unconditional basic sequence. Applying Tong's Theorem \cite{L-T,
Proposition 1.c.8} to $Q$ restricted
to the subspace
$Z$, we get that the operator $Q_D,$ defined by
$$Q_D(\sum a_k R_kTz_k)=\sum a_k
z_k^*(T^{-1}PR_kTz_k)z_k,$$
 is bounded with the norm
at most $\|T\|\|T^{-1}\|\|P\|.$

Now let us estimate the norm more directly. First 
because $R_k$ is contractive in both norms.
$$\align\|\sum a_k R_kTz_k\| &\leq  \max \{(\sum
|a_k|^p\|R_kTz_k\|_p^p)^{1/p},(\sum
|a_k|^2\|R_kTz_k\|_2^2)^{1/2} \}\\
&\leq
 \max \{(\sum
|a_k|^p\|Tz_k\|_p^p)^{1/p},\\
&\qquad\qquad (\sum
|a_k|^2\frac{9}{4} \sum_{i\in F_k} w_i^{4/(p-2)}\|Tx_i\otimes
y_j\|_2^2)^{1/2}\}\\
\intertext{(by 7.1.4)}
&\leq  \max \{(\sum
|a_k|^p \sum_{i\in F_k} w_i^{2p/(p-2)})^{1/p}\|T\|,\\
&\qquad\qquad (\sum
|a_k|^2\frac{9}{4} \sum_{i\in F_k} w_i^{2p/(p-2)}w_j^2
\|T\|^2)^{1/2}\}\\
\intertext{(by (T2))}
&\leq  \|T\|\max \{(\sum
|a_k|^p)^{1/p},(\frac{3}{2})w_j(\sum
|a_k|^2)^{1/2}\} \tag{7.1.6}
\endalign$$
since $\sum_{i\in F_k} w_i^{2p/(p-2)}\leq 1.$

Next we estimate the norm of the image under $Q_D$.
We consider only the case $a_k=1$ for all $k.$
Let $z_k'=T^{-1}PR_kTz_k$ for
all $k$.
$$\align
\|\sum_k z_k^*(z_k')z_k\| &\geq 
\max\{(\sum_k
|z_k^*(z_k')|^p\sum_{i\in F_k}
w_i^{2p/(p-2)})^{1/p},\\
&\qquad \qquad (\sum_k
|z_k^*(z_k')|^2\sum_{i\in F_k}
w_i^{2p/(p-2)})^{1/2}\}\\
&\geq  ((3\rho+\rho')/4) \max \{(\sum_k
\sum_{i\in F_k}
w_i^{2p/(p-2)})^{1/p},\\
&\qquad\qquad(\sum_k
\sum_{i\in F_k}
w_i^{2p/(p-2)})^{1/2}\}\\
\intertext{(by 7.1.5)}
&\geq  ((\rho+3\rho')/4) \max \{(\delta/2 )^{1/p}
K^{1/p},(\delta/2 )^{1/2}
K^{1/2}\}\tag{7.1.7}
\endalign$$
Combining the estimates 7.1.6 and 7.1.7 yields,
$$ \|Q_D\| \|T\|\max \{K^{1/p},(\frac{3}{2})w_jK^{1/2}\} \geq
((\rho-\rho')/4)  \max \{(\delta K/2)^{1/p},(\delta
K/2)^{1/2}\}.$$
If $K \delta\geq 2$, then the inequality becomes
$$\|T^{-1}\|\|P\| \|T\|^2\max \{K^{1/p},(\frac{3}{2})w_j K^{1/2}\} \geq
((3\rho+\rho')/4)(K\delta/2
)^{1/2},$$
and thus either $$K^{1/p}<(\frac{3}{2})w_j K^{1/2}$$
from which it would follow that 
$$\|T^{-1}\|\|P\|\|T\|^2(6\sqrt{2}) w_j>(3\rho+\rho')\delta^{1/2},$$
or
$$\sqrt{2}\|T^{-1}\|\|P\| \|T\|^2 \geq ((3\rho+\rho')/4)(\delta
)^{1/2}K^{1/2-1/p}.$$

Because we have assumed that $w_j$ is small, the former is
impossible and the conclusion follows.
\qed\enddemo

The next lemma is essentially
a reformulation of Lemma 7.1 for the special case of the
weights $(w_{n,k}).$

\proclaim{Lemma 7.2}
Suppose that $(x_i)$ and $(y_j)$ are bases of $X_{p,(w_{n,k})}$.
Let $T$ be an isomorphism from $X_{p,(w_{n,k})} \otimes X_{p,(w_{n,k})}$
into $(\sum
Y_j)_{p,2,(w_n)}$ satisfying (T2) and
let $P$ be a projection onto the range of $T$.
Let $\rho$, $\rho_1$ and $\delta $ be positive constants,
$\rho >\rho_1$.
Then there exist integers $M_m$, $m\in \N$ such that
if $w_j<\frac{\delta^{1/2}
(4\rho-\rho_1) }{6\sqrt{2} \|T\|\|T^{-1}\|\|P\|} $, 
$(N_k)_{k=1}^K$ and $(N_k')_{k=1}^K$
are strictly increasing sequences of integers with $N_k\leq N_k'\leq
N_{k+1}$ for all $k$, $(m(k))_{k=1}^K\subset \N$,
$(H_k)_{k=1}^K$ is a
sequence of finite subsets of $\N$ such that $|H_k| \geq
M_{m(k)}$ for all $k$, $H_k\cap H_{k'}=\emptyset$ if $m(k)=m(k'), k\neq k',$ and
for all $i\in H_k$
$$|x_{m(k),i}^*\otimes y_j^*((P_{N_{k}}-P_{N_{k-1}'})x_{m(k),i}\otimes y_j)|
\geq \rho
\tag{7.2.1}$$
then $K<\max\{2\delta^{-1},(\|T\|^2\|T^{-1}\| \|P\|  8/(\delta^{1/2}
(4\rho-\rho_1)))
^{2p/(p-2)}\}.$
\endproclaim

\demo{Proof}
Fix $m\in
\N.$
For $\epsilon_i=\epsilon=\rho_1 w_m^{4p/(p-2)}/16$,
$C=1,$ $D=\|T\| \|T^{-1}\| \|P\|,$ $ w_0=w_{m,j}=w_m,$
$r=p$ and $K=\delta(w_m^{2p/(p-2)})^{-1}$, we obtain an integer
$M_m$ from Lemma 6.2.

We want to construct blocks as in Lemma 7.1.
We are already given the integers $N_k$ and sets $H_k$ so
we only need to  refine the $H_k$'s to get the sets $F_k$,
define the $z_k$'s as in Lemma 7.1 and check the hypotheses of Lemma 7.1.

First we apply Lemma 6.2. We use the operator
$QT^{-1}P(P_{N_{k+1}}-P_{N_k'})T$
and the sequence $(x_{m(k+1),n}\otimes y_j)_{n\in H_{k+1}}$
to obtain a subset $F_{k+1}$ of
$H_{k+1}$ of cardinality $K_{k+1}=\delta(w_{m(k+1)}^{2p/(p-2)})^{-1}$ such that
$$\multline
|x_{m(k+1),i}^*\otimes y_j^*(T^{-1}P(P_{N_{k+1}}
-P_{N_k'})Tx_{m(k+1),i'}\otimes y_j)|\\
<\rho_1
w_{m(k+1)}^{4p/(p-2)+1}/16,
\endmultline\tag{7.2.2}$$
for all $i\neq i',i,i'\in F_{k+1}.$
(Here the parameters
$w_{n_i}$ which occur in Lemma 6.2 are not dependent on the index {i}.)

This completes the inductive definition of the sets $F_k.$
It remains to verify that our choice of $F_{k+1}$ verifies the hypothesis
of Lemma 7.1.

First we have that  for $n\in F_{k+1}$, $n\in H_{k+1}$ and thus by 7.2.1
$$
x_{m(k+1),n}^*\otimes y_j^*(T^{-1}P(P_{N_{k+1}}-P_{N_k'})
Tx_{m(k+1),n}\otimes y_j)\geq \rho$$

By 7.2.2
$$\align
&|x_{m(k+1),n}^*\otimes y_j^*(T^{-1}P((P_{N_{k+1}}-P_{N_k'})
\sum \Sb s\neq n\\ s\in F_{k+1}
\endSb w_{m(k+1)}^{2/(p-2)}T(x_{m(k+1),s}\otimes y_j)|\\
&\qquad \leq | F_{k+1}| \rho_1
w_{m(k+1)}^{4p/(p-2)+1}/16\\
&\qquad \leq \rho_1
w_{m(k+1)}^{2p/(p-2)}/16\\
&\qquad <\rho_1
w_{m(k+1)}^{2/(p-2)}/4.\tag{7.2.3}
\endalign$$

This shows that the hypotheses of Lemma 7.1 are satisfied with
$\rho,$ $\rho'=\rho-\rho_1,$ and $\delta.$
Thus $K\leq \max\{2\delta^{-1},(\|T\|^2\|T^{-1}\| \|P\|  8/(\delta^{1/2}
(4\rho-\rho_1))
^{2p/(p-2)}\}.$
\qed\enddemo

In order to show that there is no isomorphism from $X_p\otimes X_p$ into
$\rp{\alpha}$ with complemented range, we will show that any isomorphism
from $X_p\otimes X_p$ into a $(p,2)-$\tiny sum actually must have a large
part in finitely many of the summands. The previous lemmas give us tools to
use in quantitative gliding hump arguments. Lemma 7.3 is the first in a
series of lemmas which estimate how much is mapped outside a finite number of
summands.

\proclaim{Lemma 7.3} 
Suppose that $T$, $P$, and $P_N$ are as in Lemma 7.1, 
$0<\epsilon <1$, and $\rho_1>0$.
Then there exists an integer
$N_0$ and integers $M_m$, 
such that if
$$w_j<(4\epsilon-\rho_1)/(12\|T\|\|T^{-1}\|\|P\|)$$
then for each $m\in \N$, and $N\geq N_0$,
$$|\{t\in \N:
|x_{m,t}^*\otimes y_j^*(T^{-1}P(I-P_N)Tx_{m,t}\otimes
y_j)|>\epsilon\}|\leq M_m $$.
\endproclaim
\demo{Proof} For each $m$
the integers $M_m$ are chosen by the criteria established in 
Lemma 7.2 with $\rho=\epsilon$, $\delta=1/2$, and $\rho_1.$

Suppose no $N_0$ works for these $M_m.$
Then for each $N$ there exist an $m\in \N$, $N'\geq N$, and
a set $L_m\subset \N$ of cardinality at least $M_m$ such that 
$$|x_{m,s}^*\otimes y_j^*(T^{-1}P(I-P_{N'})Tx_{m,s}\otimes
y_j)|\geq \epsilon\tag{7.3.1}$$
for all $s\in L_m.$ We will apply this inductively to
construct a sequence of
pairs of
integers $(N_k,N_k')$ and finite subsets  of $\N$, $(F_k)$, as in Lemma 7.2 with
$\rho=\epsilon$ and $\delta=1/2$. Suppose we have found
$N_1,N_1',\dots,N_k$ and $F_1,\dots,F_k.$
By assumption for $N=N_k$ there exists an integer
$m_{k+1}$, $N_k'\geq N_k$ and an infinite set $L_{m_{k+1}}$ as above. 

By a simple perturbation argument we may assume that each of the elements
$Tx_{m_{k+1},i }\otimes y_j, i\in F_{k+1}$
is nonzero in a finite number of the summands $Y_n$
and thus there is an $N_{k+1}$
such that
$$P_{N_{k+1}}Tx_{m_{k+1},i }\otimes y_j=Tx_{m_{k+1},i}\otimes y_j\tag{7.3.2}$$
for all
$i\in F_{k+1}.$ This completes the induction step of the construction.

To complete the proof we need to show that if we continue in this way we can
find $K$ blocks as in Lemma 7.2 and thus reach a contradiction for $K$ large
enough. However note that
we have by 7.3.1 and 7.3.2
$$\align
|x_{m_{k+1},n}^*&\otimes y_j^*(T^{-1}P(P_{N_{k+1}}-P_{N_k'})
Tx_{m_{k+1},n}\otimes y_j)|\\
&= |x_{m_{k+1},n}^*\otimes y_j^*(T^{-1}P(I-P_{N_k'})
Tx_{m_{k+1},n}\otimes y_j)|\geq \epsilon.\endalign$$

This shows that with $\delta =1/2$ the hypotheses of Lemma 7.2 are satisfied.
Therefore we have a contradiction for large $k$
and the choice of $M_m$ works for some $N_0.$
\qed\enddemo

As we noted before Lemma 7.1 the proof simplifies and the results
strengthened if
we assume that we have an unconditional basis in the range space.

\proclaim{Lemma 7.4} Suppose that there is a constant $D$ such that for all
$n$, $Y_n$ is
a subspace of $L_p$ with a $D$-unconditional basis and that $(x_i)$ and
$(y_j)$ are bases of $X_{p,(w_{n,k})}$.
Suppose that $T$, $P_N$, and $P$ are as in Lemma 7.1,
$0<\epsilon <1$,  and $j$ is such that
$w_j<4\epsilon^2\delta/(9\|T^{-1}\|^2\|P\|^2\|T\|^2).$
Further assume that there are finite sets $N_{i,j}$ for all $i,j\in \N$
such that 
\roster
\item $T(x_i\otimes y_j)$ is supported in $N(i,j)$
\item $N(i,j)\cap N(i',j')=\emptyset$ if $i'\neq i$ and
$j\neq j'$; $i=i'$, $j\neq j'$ and \newline $\max(o(j),o(j'))>o(i)$;
or $i\neq i'$, $j=j'$ and $\max(o(i),o(i'))\geq o(j).$
\endroster
Then there exists an integer
$N_0$ 
such that for all $N\geq N_0$, $\sum_{n\in F} w_{m,n}^{2p/(p-2)}<\infty$
where $F=\{n:|x_n^*\otimes y_j^*(T^{-1}P(I-P_N)Tx_n\otimes
y_j)|\geq\epsilon\}.$
\endproclaim

\demo{Proof} We use arguments like that in the proofs of Lemmas
7.1, 7.2 and 7.3 but we use unconditionality to avoid the use of Lemma 6.2. 

Suppose that no such $N_0$ exists and fix $\delta, 1\geq \delta>1/2$
Then we can find a strictly increasing sequence of integers $(N_k)$ and
finite disjoint subsets $(F_k)$ of $\N$ such that
$$ 1 \geq \sum_{s\in F_k} w_s^{2p/(p-2)}> \delta \text{ for each $k$,}$$
$$N(s,j)\cap N(n,j)=\emptyset \text{ if $s,n\in \cup_k F_k$ and $s\neq k,$}$$
and for all $n\in F_k,$
$$|x_n^*\otimes y_j^*(T^{-1}P(I-P_{N_k})Tx_n\otimes
y_j)|\geq \epsilon$$
and
$$P_{N_{k+1}}x_n\otimes y_j=x_n\otimes y_j.$$

Let $x_n'=(I-P_{N_k})Tx_n\otimes y_j$ for all $n\in F_k$ and all $k.$ 
$(x_n')_{n\in F_k,k\in \N}$ is an unconditional basic sequence in
$(\sum Y_n)_{p,2}.$
Define an operator $S$ from $[x_n':n\in F_k,k\in \N]$ into $[x_n\otimes
y_j:n\in F_k,k\in \N]$ by $QT^{-1}P$ where $Q$ is the basis projection onto
$[x_n\otimes
y_j:n\in F_k,k\in \N]$.
By Tong's Theorem  the diagonal operator $S_D$ defined by
$S_D(\sum_n a_n x_n')=\sum_n a_n (x_n^*\otimes
y_j^*(T^{-1}Px_n'))x_n\otimes y_j$ is bounded. As in the proof of Lemma 7.1
choose signs $\epsilon_n$ such that $\|\sum_{n\in F_k} \epsilon_n
w_n^{2/(p-2)} Tx_n\otimes y_j\|_2\leq (3/2)\|T\|(\sum_{n\in F_k} w_n^{2p/(p-2)}
w_j^2)^{1/2}.$

Let $$z_k'=\sum_{n\in F_k} \epsilon_n w_n^{2/(p-2)} x_n',$$
$$z_k=\sum_{n\in F_k}
\epsilon_n w_n^{2/(p-2)}
x_n\otimes y_j,$$
and $$z_k^*=(\sum_{n\in F_k} w_n^{2p/(p-2)})^{-1}
\sum_{n\in F_k} \epsilon_n w_n^{2(p-1)/(p-2)}x_n^*\otimes
y_j^*.$$
Let $Q'$ be the usual projection onto $[z_k]$. Then 
$$Q'S_D z_k'=
z_k^*(S_D z_k')z_k= (\sum_{n\in F_k} w_n^{2p/(p-2)})^{-1}
(\sum_{n\in F_k} w_n^{2p/(p-2)} x_n^*\otimes
y_j^*(T^{-1}Px_n'))z_k.$$
Notice that  $|x_n^*\otimes
y_j^*(T^{-1}Px_n')|\geq \epsilon $ for all $n.$
Thus
$$\align
\|S_D\|\max\{(\sum_{i=1}^K \|z_i'\|_p^p)^{1/p},(\sum_{i=1}^K
\|z_i'\|_2^2)^{1/2}\} &\geq \|Q'S_D(\sum_{i=1}^K z_i')\|\\
&\geq \|\sum_{i=1}^K (\epsilon) z_k\|\\
&\geq \epsilon\max\{(K\delta)^{1/p},(K\delta)^{1/2}\}.
\endalign$$
Because 
$$\|z_i'\|_2 \leq (3/2)\|T\|(\sum_{n\in F_k} w_n^{2p/(p-2)}
w_j^2)^{1/2},$$
$$(\sum_{i=1}^K
\|z_i'\|_2^2)^{1/2}\leq (3/2)\|T\|K^{1/2}w_j^{1/2}.$$
For $K>\delta^{-1}$, if
$(3/2)\|T\|K^{1/2}w_j^{1/2}>\|T\|K^{1/p}$ then we have that
$$\|S_D\|C(3/2)\|T\|K^{1/2}w_j^{1/2}\geq \epsilon(K\delta)^{1/2}.$$
But $w_j<4\epsilon^2\delta/(9\|S_D\|^2\|T\|^2)$, so we must have
$(3/2)K^{1/2}w_j^{1/2}\leq K^{1/p}$ and therefore
$$\|S_D\|C \|T\|K^{1/p}\geq \epsilon(K\delta)^{1/2}$$
or equivalently,
$$K\leq (\|S_D\|C\|T\|/(\epsilon\delta^{1/2}))^{2p/(p-2)}.$$
Thus the construction can only be made finitely many times and the claimed
$N_0$ exists.
\qed\enddemo

After our detour into the case of unconditional basis, we continue
enlarging the sets of indices for which basis vectors are mapped into a
finite number of summands.

\proclaim{Lemma 7.5}Suppose that
$(x_{i,h})$ and
$(y_{j,l})$ are standard bases of $X_{p,(w_{n,k})}$.
Let $T$ be an isomorphism from $X_{p,(w_{n,k})} \otimes X_{p,(w_{n,k})}$
into $(\sum
Y_j)_{p,2,(w_n)}$ satisfying (T2) and let $P$ be a
projection onto the range of $T$.
Let $P_N$ denote the projection onto the first $N$ summands
of $(\sum Y_n)_{p,2}.$ Let $\epsilon>0$ and let
$M_j$ and $M_i$ be the integers given by Lemma
7.2 for
$\rho=3\epsilon/4$, $\delta=1/2$, and $\rho_1=\epsilon/4.$
Then for each $i,j$ such that
$\max\{w_i,w_j\}<(100 \|T\|\|T^{-1}\|\|P\|)^{-1}$,
there exist two
infinite sets $H,L$ and an integer $N_0$ such that 
for all $N>N_0$,
if $L'\subset L$ and  $|L'|\geq M_j$ then
$$|\{h:|x_{i,h}^*\otimes
y_{j,l}^*(T^{-1}P(I-P_N)Tx_{i,h}\otimes
y_{j,l})|\geq\epsilon\quad \forall l\in L'\}|<\infty ,$$
and
if $H'\subset H$ and $|H'|\geq M_i$ then
$$|\{l:|x_{i,h}^*\otimes
y_{j,l}^*(T^{-1}P(I-P_N)Tx_{i,h}\otimes
y_{j,l})|\geq\epsilon\quad\forall h\in H' \}|<\infty.$$
Moreover, for all $h\in H, l\in L$,
$$|x_{i,h}^*\otimes
y_{j,l}^*(T^{-1}P(I-P_{N_0})Tx_{i,h}\otimes
y_{j,l})|<\epsilon.$$
\endproclaim

\demo{Proof} The proof is by induction.
Just to get the induction started let $h_1=1$ and  by Lemma 7.3 there is
an $N_1\in \N$ such that for
any $N\geq N_1$,
$|x_{i,1}^*\otimes
y_{j,l}^*(T^{-1}P(I-P_{N})Tx_{i,1}\otimes
y_{j,l})|<\epsilon$ for all but $M_j$ many $l$.
For $N=N_1$ there are only finitely many $l$ such that the inequality above
fails, so let $L_1$ be the set of $l$ such that
$|x_{i,1}^*\otimes
y_{j,l}^*(T^{-1}P(I-P_{N_1})Tx_{i,h}\otimes
y_{j,l})|<\epsilon$.

Consider the following (non-mutually exclusive) possibilities.
\roster
\item For every infinite $H\subset \N$ and $N\geq N_1$
there are infinitely many $l\in L_1$ 
and infinite 
subsets $H_l$ of $H$ such that for all $h\in H_l$, 
$$|x_{i,h}^*\otimes
y_{j,l}^*(T^{-1}P(I-P_{N})Tx_{i,h}\otimes
y_{j,l})|<\epsilon.$$

or

\item There is an $N\geq N_1$, an infinite set $H$
and  there are $M_j$ integers $l$
such that
$$|x_{i,h}^*\otimes
y_{j,l}^*(T^{-1}P(I-P_{N})Tx_{i,h}\otimes
y_{j,l})|\geq\epsilon$$ for all $h\in H.$
\endroster

If the second possibility occurs, we can find $N_1'\geq N_1$, a subset
$L_1'$ of $L_1$ with cardinality $M_j$, and an infinite
subset $H_1$ of $\N$ such that for all
$h\in H_1, l\in L_1'$
$|x_{i,h}^*\otimes
y_{j,l}^*(T^{-1}P(I-P_{N_1'})Tx_{i,h}\otimes
y_{j,l})|\geq\epsilon$. By applying Lemma 7.3 at most $M_j$ times we can
find an integer $N_2>N_1'$ such that for any $N\geq N_2$ and $l\in L_1'$,
$$|\{h:|x_{i,h}^*\otimes
y_{j,l}^*(T^{-1}P(I-P_{N})Tx_{i,h}\otimes
y_{j,l})|\geq\epsilon/4\}|<M_i'.$$
(We are assuming that the integer $M_i'$ 
obtained from Lemma 7.3 is not necessarily the same as
that we have obtained from Lemma 7.2.)
By eliminating a finite number of $h \in H_1$ we may assume that
$$|x_{i,h}^*\otimes
y_{j,l}^*(T^{-1}P(I-P_{N_2})Tx_{i,h}\otimes
y_{j,l})|<\epsilon/4$$
 for all $h\in H_1\cup \{h_1\}.$ Let $l_1$ be any
element of $L_1'$. Notice that for any $h\in H_1$, we have
$$|x_{i,h}^*\otimes
y_{j,l}^*(T^{-1}P(P_{N_2}-P_{N_1'})Tx_{i,h}\otimes
y_{j,l})|\geq\epsilon-\epsilon/4=3\epsilon/4$$ and thus we have a block as
in Lemma 7.2.

If the second possibility does not occur, then we choose $l_1\in L_1$ and
an infinite set $H_1$ such that for all $h\in H_1$,
$|x_{i,h}^*\otimes
y_{j,l_1}^*(T^{-1}P(I-P_{N_1})Tx_{i,h}\otimes
y_{j,l_1})|<\epsilon$. Let $N_2=N_1.$

Now we choose $h_2,$ by interchanging the roles of $l$ and $h$ in the
argument above.

If there is an $N\geq N_2$, an infinite set $L\subset L_1$
 and  there are $M_i$ integers $h$
such that
$|x_{i,h}^*\otimes
y_{j,l}^*(T^{-1}P(I-P_{N})Tx_{i,h}\otimes
y_{j,l})|\geq\epsilon$ for all $l\in L$, then let  
$N_2'>N_2$, $H_2'\subset H_1$ with cardinality $M_i$, and  let $L_2$ be an
infinite subset of $L_1\setminus L_1'$ such that 
for all $l\in L_2$ and $h\in H_2'$, $|x_{i,h}^*\otimes
y_{j,l}^*(T^{-1}P(I-P_{N_1'})Tx_{i,h}\otimes
y_{j,l})|\geq\epsilon$. By applying Lemma 7.3 at most $M_i$ times we can
find an integer $N_3>N_2'$ such that for any $N\geq N_3$ and $h\in H_2'$,
$$|\{l:|x_{i,h}^*\otimes
y_{j,l}^*(T^{-1}P(I-P_{N})Tx_{i,h}\otimes
y_{j,l})|\geq\epsilon/4\}|<M_j'.$$
($M_j'$ is the integer from Lemma 7.3.)
By eliminating a finite number of $l \in L_2$ we may assume that
$|x_{i,h}^*\otimes
y_{j,l}^*(T^{-1}P(I-P_{N_2})Tx_{i,h}\otimes
y_{j,l})|<\epsilon/4$ for all $l\in L_2\cup \{l_1\}.$ Let $h_2$ be any
element of $H_2'$.
Notice that for any $l\in L_2$, we have
$$|x_{i,h}^*\otimes
y_{j,l}^*(T^{-1}P(P_{N_3}-P_{N_2'})Tx_{i,h}\otimes
y_{j,l})|\geq\epsilon-\epsilon/4=3\epsilon/4$$ and thus we again have a block as
in Lemma 7.2.
 
 Otherwise choose $h_2$ and and an infinite set $L_2\subset L_1$
such that for all $l\in L_2$,$|x_{i,h_2}^*\otimes
y_{j,l}^*(T^{-1}P(I-P_{N_2})Tx_{i,h_2}\otimes
y_{j,l})|<\epsilon$. Let $N_3=N_2.$

The remainder of the proof proceeds by alternately choosing elements $h_k$
and $l_k$ as above. Notice that the number of integers $k$ for which
we find the set $H_k'$ of
cardinality $M_i$, $L_k'$ of cardinality $M_j$, respectively, is limited by
Lemma 7.2. Therefore, after some finite number of steps, we have that
there is an $N_{k_0}$, and infinite sets $L_{k_0}$ and $H_{k_0}$ such that 
for any $N\geq N_{k_0}$ and infinite subsets $L$ of $L_{k_0}$ and  $H$ of
$H_{k_0}$,
$$|\{l\in L_{k_0}:|x_{i,h}^*\otimes       
y_{j,l}^*(T^{-1}P(I-P_{N})Tx_{i,h}\otimes
y_{j,l})|\geq\epsilon\text{ for all }h\in H\}|<M_j$$
 and 
$$|\{h\in
H_{k_0}:|x_{i,h}^*\otimes       
y_{j,l}^*(T^{-1}P(I-P_{N})Tx_{i,h}\otimes
y_{j,l})|\geq\epsilon\text{ for all }l\in L\}|<M_i.$$
The sets $H=\{h_k:k\geq k_0\}$ and $L=\{l_k:k\geq k_0\}$ and $N_0=N_{k_0}$
satisfy the conclusion of the lemma.
\qed\enddemo.

The previous lemma gave us an estimate for a row in each factor. Now we
move to a rich set in one factor.

\proclaim{Lemma 7.6}
Let $(x_{i,h})$ and
$(y_{j,l})$ be equivalent to the standard basis of $X_{p,2,(w_{n,k})}$.
Suppose that $T$, $P_N$, and $P$ are as in Lemma 7.1,
$0<\epsilon <1$,  and $j$ is such that
$w_j=w_{j,l}<(1-\epsilon)/(3\|S_D\|\|T\|)$.
Let $(M_i)$ be the sequence of integers given by Lemma 7.2 for $\rho
=3\epsilon/4, \delta=1/2,$ and $\rho_1=\epsilon/4.$
Then there exists an integer
$N_0$, an infinite subset $L$ of $\N$, and a rich subset $K$ of $\NxN,$
such that for all $N\geq N_0$, $j\in \N$, there is a
set $L'\subset L,$ $|L'|<M_j$ and for all $l\notin L'$,
$\{h:(i,h)\in K,|x_{i,h}^*\otimes y_{j,l}^*(T^{-1}P(I-P_N)Tx_{i,h}\otimes
y_{j,l})|< \epsilon\}$ is infinite. Moreover, for all $l\in L,(i,h)\in K$,
$|x_{i,h}^*\otimes y_{j,l}^*(T^{-1}P(I-P_{N_0})Tx_{i,h}\otimes
y_{j,l})|<\epsilon.$
\endproclaim

\demo{Proof} By discarding the first few rows of
$\NxN$, i.e., $\{1,2,..,k\}\times \N$ for some $k$, and renumbering we may
assume that $w_i=w_{i,h}< (1-\epsilon)/(3\|T^{-1}\|\|T\|)$
for all $i,h\in \N.$ We will inductively construct $L$ and $K$
by using Lemma 7.5 and Lemma 7.2.

First by Lemma 7.5 there exist an integer $N_1$, 
 an infinite subset $L_1$ of $\N$ and an
infinite subset $H_1$ of $\{1\}\times \N$ such that for all $N\geq N_1$,
$$\{(h,l):|x_{1,h}^*\otimes y_{j,l}^*(T^{-1}P(I-P_N)Tx_{i,h}\otimes
y_{j,l})|<\epsilon\}$$
 does not contain a rectangle $A\times B$ with $A$
infinite and cardinality of $B$ greater than or equal to $M_j$ or  with
$B$
infinite and cardinality of $A$ greater than or equal to $M_1$. Also
$|x_{1,h}^*\otimes y_{j,l}^*(T^{-1}P(I-P_{N_1})Tx_{i,h}\otimes
y_{j,l})|<\epsilon$ for all $l\in L_1$ and $h\in H_1.$

Consider
$$\multline
L^N=\{l\in L_1:\exists 
 K\text{ rich such that }\\
|x_{i,h}^*\otimes
y_{j,l}^*(T^{-1}P(I-P_{N})Tx_{i,h}\otimes
y_{j,l})|<\epsilon \qquad\forall (i,h)\in K\}.
\endmultline$$
If for some $N_1'\geq N_1,$
$|L_1\setminus L^{N_1'}|\geq M_j$ then let $L_1'\subset
L_1\setminus L_1^{N_1'}$ with cardinality $M_j$. If we enumerate the elements of
$L_1'$ as $(l_i)_{i=1}^{M_j}$ we can produce a rich set which is bad for
all $l_i\in L_1'$ inductively as follows.

Because $l_1\notin L^{N_1'}$, there exists a rich set $K_1$
such that $$|x_{i,h}^*\otimes
y_{j,l_1}^*(T^{-1}P(I-P_{N_1'})Tx_{i,h}\otimes
y_{j,l_1})|<\epsilon \qquad\forall (i,h)\in K_1.$$
Because $K_1$ is a rich set and $l_2\notin L^{N_1'}$,
there is a rich subset $K_2$
of $K_1$ such that 
$$|x_{i,h}^*\otimes
y_{j,l_2}^*(T^{-1}P(I-P_{N_1'})Tx_{i,h}\otimes
y_{j,l_2})|<\epsilon \qquad\forall (i,h)\in K_2.$$
Similarly there exists a rich subset $K_3$ of $K_2$ such that 
$$|x_{i,h}^*\otimes
y_{j,l_3}^*(T^{-1}P(I-P_{N_1'})Tx_{i,h}\otimes
y_{j,l_3})|<\epsilon \qquad\forall (i,h)\in K_3.$$
Continuing in this way we find a decreasing sequence of rich sets
$(K_i)_{i=1}^{M_j}$ such that 
$$|x_{i,h}^*\otimes
y_{j,l_k}^*(T^{-1}P(I-P_{N_1'})Tx_{i,h}\otimes
y_{j,l_k})|<\epsilon \qquad\forall (i,h)\in K_s,s\leq k.$$
Let $K_1'=K_{M_j}.$ By at most $M_j$ applications of Lemma 7.3 there exists
an integer $N_1''>N_1'$ such that  for all $N\geq N_1''$, $l\in L_1'$, and $i\in
K^1,$
$$|\{h:|x_{i,h}^*\otimes
y_{j,l}^*(T^{-1}P(I-P_{N})Tx_{i,h}\otimes
y_{j,l})|\geq \epsilon/4\}|<M_i',$$
where $M_i'$ is given by Lemma 7.3.
Let $i_2$ be the smallest index in $K_1'{}^1.$ By Lemma 7.5 there is an
integer $N_2>N_1''$, an 
infinite subset $L_2$ of $L_1\setminus L_1'$ and an infinite subset $H_2$
of $\{h:(i_2,h)\in  K_1'\}$,
such that  for any $N\geq N_2,$
$\{(h,l):|x_{i,h}^*\otimes y_{j,l}^*(T^{-1}P(I-P_N)Tx_{i,h}\otimes
y_{j,l})|<\epsilon\}$ does not contain a rectangle $A\times B$ with $A$
infinite and cardinality of $B$ greater than or equal to $M_j$ or  with
$B$
infinite and cardinality of $A$ greater than or equal to $M_{i_2}$. Also
$|x_{i,h}^*\otimes y_{j,l}^*(T^{-1}P(I-P_{N_2})Tx_{i,h}\otimes
y_{j,l})|<\epsilon$ for all $l\in L_2$ and $h\in H_2.$
By discarding at most a finite number of elements from each row of $K_1'$
we may assume that for all $(i,h)\in K_1'$ and $l\in L_1'$ and for all
$(1,h)$ with $h\in H_1$ and $l\in L_2\cup L_1',$
$$|x_{i,h}^*\otimes
y_{j,l}^*(T^{-1}P(I-P_{N})Tx_{i,h}\otimes
y_{j,l})|< \epsilon.$$

Observe that for any $(i,h)\in K_1$ we have that for each $l\in L_1',$
$$|x_{i,h}^*\otimes y_{j,l}^*(T^{-1}P(P_{N_2}-P_{N_1})Tx_{i,h}\otimes
y_{j,l})|>\epsilon -\epsilon/4=3\epsilon/4.$$
Thus we have a block  with respect to $l$ as in
Lemma 7.2.

If we cannot find the set $L_1'$ of cardinality $M_j$
as above, choose $l_1\in L_1$ and let
$K_1$ be a rich subset of $L^{N_1}$ such that
$|x_{i,h}^*\otimes y_{j,l_1}^*(T^{-1}P(I-P_{N_1})Tx_{i,h}\otimes
y_{j,l_1})|<\epsilon$ for all $(i,h)\in K_1$. Note that we may assume that
$K_1$ is maximal and thus that $(1,h)\in K_1$ for all $h\in H_1.$
Let $L_1'=\{l_1\}.$

For each $i\in K_1^1$, $N\geq N_1$,
 and infinite set $L\subset L_2$
consider the set 
$$H(i,N,L)=\{h:|x_{i,h}^*\otimes
y_{j,l}^*(T^{-1}P(I-P_{N})Tx_{i,h}\otimes
y_{j,l})|\geq \epsilon\qquad\forall l\in L\}.$$
 If the cardinality of
$H(i,N,L)$ is at least $M_i$ for some $i=i_2,$ $L=L_1'$
 and $N=N_1'\geq N_1,$ let $H_2'$ be a 
subset of $H(i_2,N_1',L_1')$ with cardinality $M_i.$ By Lemma 7.5 there is an
integer $N_2>N_1'$, an
infinite subset $L_2$ of $L_1'$ and an infinite subset $H_2$
of $\{h:(i_2,h)\in  K_1'\}$,
such that  for any $N\geq N_2,$
$$\{(h,l):|x_{i,h}^*\otimes y_{j,l}^*(T^{-1}P(I-P_N)Tx_{i,h}\otimes
y_{j,l})|<\epsilon/2\}$$
 does not contain a rectangle $A\times B$ with $A$
infinite and cardinality of $B$ greater than or equal to $M_j$ or  with
$B$
infinite and cardinality of $A$ greater than or equal to $M_{i_2}$. Also
$|x_{i,h}^*\otimes y_{j,l}^*(T^{-1}P(I-P_{N_2})Tx_{i,h}\otimes
y_{j,l})|<\epsilon/2$ for all $l\in L_2$ and $h\in H_2.$ By passing to
infinite subsets of $H_1$ and $L_2$ we may assume that
$|x_{i,h}^*\otimes y_{j,l}^*(T^{-1}P(I-P_{N_2})Tx_{i,h}\otimes
y_{j,l})|<\epsilon$ for all $l\in L_2\cup \{l_1\}$ and $h\in H_1.$

Observe that for any $h\in H_2'$ we have that for each $l\in L_2,$
$$|x_{i_2,h}^*\otimes y_{j,l}^*(T^{-1}P(P_{N_2}-P_{N_1})Tx_{i_2,h}\otimes
y_{j,l})|>\epsilon -\epsilon/4=3\epsilon/4.$$
Thus we have a block with respect to $h$ as in
Lemma 7.2.

If neither $L_1'$ of cardinality $M_j$ as in the first case
nor $H_2'$ of cardinality $M_i$ as
in the second case can  be found, then choose an infinite subset
$L_2$ of $L_1$ and an infinite subset $H_2$ of $\{h:(i,h)\in K_1\}$ such
that $|x_{i,h}^*\otimes y_{j,l}^*(T^{-1}P(I-P_{N_1})Tx_{i,h}\otimes
y_{j,l})|<\epsilon$ for all $h\in H_2$ and $l\in L_2.$ Let $N_2=N_1.$
let $l_1$ be any element of $L_2$ and let $K_2$ be a  maximal rich subset
of $\NxN$ such that
$$|x_{i,h}^*\otimes y_{j,l}^*(T^{-1}P(I-P_{N_2})Tx_{i,h}\otimes
y_{j,l})|<\epsilon$$
 for all $(i,h)\in K_2.$

This completes the first full step of our construction. Notice that in each
case we have produced an integer $N_2$,
infinite sets $L_2,H_1, H_2$ of $\N$, a set $L_1'$
containing at least one
element and a rich set $K_2$ such that 
$$|x_{i,h}^*\otimes
y_{j,l}^*(T^{-1}P(I-P_{N_2})Tx_{i,h}\otimes
y_{j,l})|<\epsilon$$
 for all $l\in L_2\cup \{l_1\}$ and $(i,h)\in \{1\}\times
H_1\cup {i_2}\times H_2$ and 
for all $(i,h)\in K_2$ and $l\in L_1'$,
$$|x_{i,h}^*\otimes
y_{j,l}^*(T^{-1}P(I-P_{N_2})Tx_{i,h}\otimes
y_{j,l})|<\epsilon.$$

We will present one more step of the induction.

Consider $$L^N=\{l\in L_2:\exists
 K\text{ rich }\backepsilon |x_{i,h}^*\otimes
y_{j,l}^*(T^{-1}P(I-P_{N})Tx_{i,h}\otimes
y_{j,l})|<\epsilon\quad \forall (i,h)\in K\}.$$
If for some $N_2'\geq N_2,$
$|L_2\setminus L^{N_2'}|\geq M_j$ then let $L_2'\subset
L_2^{N_2'}\setminus L$ with cardinality $M_j$. If we enumerate the elements of
$L_2'$ as $(l_i)_{i=M_j+1}^{2M_j}$ we can produce a rich set $K_2'$
 which is bad for
all $l_i, M_j<i\leq 2M_j$ inductively as above.
 By at most $M_j$ applications of Lemma 7.3 there exists
an integer $N_2''>N_2'$ such that  for all $N\geq N_2''$, $l\in L_2'$, and
$i\in (K_2')^1,$
$$|\{h:|x_{i,h}^*\otimes
y_{j,l}^*(T^{-1}P(I-P_{N})Tx_{i,h}\otimes
y_{j,l})|\geq \epsilon/4\}|<M_i'.$$
Let $i_3$ be the smallest index in $(K_2')^1.$ By Lemma 7.5 there is an
integer $N_3>N_2''$, an
infinite subset $L_3$ of $L_2\setminus L_2'$ and an infinite subset $H_3$
of $\{h:(i_3,h)\in  K_2'\}$,
such that  for any $N\geq N_3,$
$\{(h,l):|x_{i,h}^*\otimes y_{j,l}^*(T^{-1}P(I-P_N)Tx_{i,h}\otimes
y_{j,l})|<\epsilon\}$ does not contain a rectangle $A\times B$ with $A$
infinite and cardinality of $B$ greater than or equal to $M_j$ or  with
$B$
infinite and cardinality of $A$ greater than or equal to $M_{i_3}$. Also
$|x_{i,h}^*\otimes y_{j,l}^*(T^{-1}P(I-P_{N_3})Tx_{i,h}\otimes
y_{j,l})|<\epsilon$ for all $l\in L_3$ and $h\in H_3.$
By discarding at most a finite number of elements from each row of $K_2'$
we may assume that for all $(i,h)\in K_2'$ and $l\in L_2'$ and for all
$(i,h)$ with $h\in H_i$ for $i=1,2$ and $l\in L_3\cup L_2'\cup L_1',$
$$|x_{i,h}^*\otimes
y_{j,l}^*(T^{-1}P(I-P_{N})Tx_{i,h}\otimes
y_{j,l})|< \epsilon.$$

Observe that for any $(i,h)\in K_2$ we have that for each $l\in L_2',$
$$|x_{i,h}^*\otimes y_{j,l}^*(T^{-1}P(P_{N_2}-P_{N_1})Tx_{i,h}\otimes
y_{j,l})|>\epsilon -\epsilon/4=3\epsilon/4.$$
Thus we have a new block  with respect to $l$ as in
Lemma 7.2. (If $N_2=N_1$ this case cannot occur.)

If we cannot find the set $L_2'$ of cardinality $M_j$
as above, choose $l_2\in L_2$ and let
$K_2$ be a rich subset of $L^{N_2}$ such that
$|x_{i,h}^*\otimes y_{j,l_2}^*(T^{-1}P(I-P_{N_1})Tx_{i,h}\otimes
y_{j,l_2})|<\epsilon$ for all $(i,h)\in K_2$. Note that we may assume that
$K_2$ is maximal and thus that $(i_s,h)\in K_2$ for all $h\in H_s, s=1,2.$
Let $L_2'=\{l_2\}.$

For each $i\in (K_2)^1$, $N\geq N_2$,
 and infinite set $L\subset L_2$
consider the set
$$H(i,N,L)=\{h:|x_{i,h}^*\otimes
y_{j,l}^*(T^{-1}P(I-P_{N})Tx_{i,h}\otimes
y_{j,l})|\geq \epsilon\qquad\forall l\in L\}.$$
 If the cardinality of
$H(i,N,L)$ is at least $M_i$ for some $i=i_3,$ $L=L_2'$
 and $N=N_2'\geq N_2,$ let $H_3'$ be a
subset of $H(i_3,N_2',L_2')$ with cardinality $M_i.$ By Lemma 7.5 there is
an
integer $N_3>N_2'$, an
infinite subset $L_3$ of $L_2'$ and an infinite subset $H_3$
of $\{h:(i_3,h)\in  K_2'\}$,
such that  for any $N\geq N_2,$
$$\{(h,l):|x_{i,h}^*\otimes y_{j,l}^*(T^{-1}P(I-P_N)Tx_{i,h}\otimes
y_{j,l})|<\epsilon/4\}$$
 does not contain a rectangle $A\times B$ with $A$
infinite and cardinality of $B$ greater than or equal to $M_j$ or  with
$B$
infinite and cardinality of $A$ greater than or equal to $M_{i_3}$. Also
$|x_{i,h}^*\otimes y_{j,l}^*(T^{-1}P(I-P_{N_3})Tx_{i,h}\otimes
y_{j,l})|<\epsilon/4$ for all $l\in L_3$ and $h\in H_3.$ By passing to
infinite subsets of $H_3$ and $L_3$ we may assume that
$|x_{i,h}^*\otimes y_{j,l}^*(T^{-1}P(I-P_{N_3})Tx_{i,h}\otimes
y_{j,l})|<\epsilon$ for all $l\in L_3\cup L_1'\cup L_2'$ and $h\in H_1\cup
H_2.$

Observe that for any $h\in H_3'$ we have that for each $l\in L_3,$
$$|x_{i_3,h}^*\otimes y_{j,l}^*(T^{-1}P(P_{N_2}-P_{N_1})Tx_{i_3,h}\otimes
y_{j,l})|>\epsilon -\epsilon/4=3\epsilon/4.$$
Thus we have a new block with respect to $h$ as in
Lemma 7.2.

If neither $L_2'$ of cardinality $M_j$ as in the first case
nor $H_3'$ of cardinality $M_i$ as
in the second case can  be found, then choose an infinite subset
$L_3$ of $L_2$ and an infinite subset $H_3$ of $\{h:(i,h)\in K_2\}$ such
that $|x_{i,h}^*\otimes y_{j,l}^*(T^{-1}P(I-P_{N_2})Tx_{i,h}\otimes
y_{j,l})|<\epsilon$ for all $h\in H_3$ and $l\in L_3.$ Let $N_3=N_2.$
let $l_2$ be any element of $L_3$ and let $K_3$ be a  maximal rich subset
of $\NxN$ such that
$$|x_{i,h}^*\otimes y_{j,l}^*(T^{-1}P(I-P_{N_3})Tx_{i,h}\otimes
y_{j,l})|<\epsilon$$
 for all $(i,h)\in K_3.$

We have now completed a second full step of the induction. Continuing in
this way we produce an increasing sequence of integers $(N_k)$, $(i_k)$ and
sequences of sets of integers $(L_k')$, $(H_k)$, $(L_k)$, and rich sets
$(K_k)$. The sequence $(N_k)$ must be eventually constant since each
increase in $N_k$ is produced when a new set $L_k'$ of cardinality $M_j$ is
found or a new set $H_k'$ of cardinality $M_{i_k}$ is found. If
$N_{k_1}<\dots <N_{k_s}$ and for each $r, 1\leq r\leq s,$ we have disjoint
sets $L_{k_r}'$, then any $(i,h)\in K_{s+1}$ will give us $s$ blocks as in
Lemma 7.2. But Lemma 7.2 gives a bound on the number of such blocks and thus
this cannot be the cause of the increase in $N_k$. Similarly, if we have
disjoint sets $H_{k_r}'$ of cardinality $M_{i_r}$ then choosing any $l\in
L_{s+1}$ will give us $s$ blocks as in Lemma 7.2. Therefore there is an
integer $k_0$ such that $N_k=N_{k_0}$ for all $k\geq k_0.$ Consequently, we
can let $N_0=N_{k_0}$, $L=\{l_r:r\geq k_0\}$ and $K=\cup_{r\geq k_0}
 \{i_r\}\times H_r.$
\qed\enddemo

\newpage

\head $X_p\otimes X_p$ is not in the scale $R_p^\alpha,\alpha<\omega_1$
\endhead

In this section we prove our main result and answer a question posed in
\cite{BRS}. Before beginning we need a few combinatorial lemmas.

\proclaim{Lemma 8.1} If $G\subset \NxN$ and there exists an integers
$M$ such that if $A\subset \N$ with cardinality $M$ and $B\subset \N$ which
is infinite, then $A\times B \cap G\neq \emptyset$ and $B\times A \cap
G\neq \emptyset$.
Then there exist infinite subsets of $\N$, $H,L$ such that $H\times
L\subset G.$
\endproclaim
\demo{Proof} Observe that the hypothesis implies the following.

Given $K\subset 
\N$ with cardinality greater than $M-1$ and any infinite subset $J$ of $\N$
then there is an element $n$ of $N$ and an infinite set $J'\subset J.$ such
that $(n,j)\in G$ for all $j\in J'.$

Indeed, enumerate the elements of $N$
as $n_1, n_2, \dots n_{M}.$  By hypothesis $\cap \{j\in J:(n_r,j)\notin
G\}$ is finite and
$$\multline
\cap \{j\in J:(n_r,j)\notin G\}
=J\setminus(\{j\in
J:(n_1,j)\in G\}\cup \{j\in J:(n_1,j)\notin G,\\
 (n_2,j)\in G\}\cup \dots\cup
\{j\in J:(n_r,j)\notin G\text{ for $r=1,2,\dots,M-1$, } (n_M,j)\in G\}.
\endmultline$$
thus one of the sets 
$$\multline
\{j\in J:(n_1,j)\notin G, (n_2,j)\in G\}\cup \dots
\{j\in J:(n_r,j)\notin G \\
\text{ for $r=1,2,\dots,s$, }(n_{s+1},j)\in G\}
\endmultline$$
is infinite.

To construct the sets $H,L$ we alternately select infinite sets $H_k$ and
$L_k$ and elements $h_k$ and $l_k$ such that
\roster
\item $(H_k)$ and $(L_k)$ are decreasing,
\item $(h_k,l)\in G$ for all $l\in L_{k+1}$,
\item $(h,l_k)\in G$ for all $h\in H_{k}$,
\item $h_{k}\in H_k$
\item $l_{k}\in L_k$
\endroster

Begin with $H_0=\N$ and $L_1=\N,$
and use the principle to find $l_1\in \{1,2,\dots M\}$
and an infinite set $H_1$ such that $(h,l_1)\in G$ for all $h\in H_1.$ Next
let $H_1'$ be a subset of $H_1$ with cardinality $M$ and let
$L_0=\N\setminus \{l_1\}.$ By the principle there exists an element $h_1$
of $H_1'$ and an infinite subset $L_2$ of $L_1$ such that $(h_1,l)\in G$
for all $l\in L_2.$ Next let $L_2'$ be a subset of $L_2\setminus \{l_1\}$
with cardinality $M$. By the principle there is an infinite
subset $H_2$ of $H_1\setminus \{h_1\}$ and $l_2\in L_2'$ such that
$(h,l_2)\in G$ for all $h\in H_2.$ Clearly this procedure will produce the
required sequences.
\qed\enddemo

\proclaim{Lemma 8.2} Suppose that $G \subset \N\times \NxN$ and $N\in \N$
then one of the following holds
\roster
\item There exists a rich subset $K\subset \NxN$ such
that for each $(j,l)\in K$, there are at least $N$ elements $n \in \N$ such
that $(n,j,l)\notin G.$
\item There exist $M\subset \N$, $M$ infinite, and a rich subset $K\subset
\NxN$ such
that for each $n\in M$ and $(j,l)\in K$, $(n,j,l)\in G.$
\endroster
\endproclaim

\demo{Proof} We begin by trying to directly satisfy the first alternative.
For each $j\in \N$ let $$L_j=\{l:\text{ there exists }F\subset \N, 
|F|\geq N \text{ such that } (n,j,l)\notin G \text{ for all }n\in F\}.$$
Let $J=\{j: |L_j|=\infty\}.$ If $J$ is infinite, then we define
$K=\{(j,l):j\in J, l\in L_j\}.$ $K$ is clearly rich and satisfies the first
criterion.

If $J$ is finite, discard all $(n,j,l)$ such that $j\in J$. Because
$\N\setminus L_j$ is infinite for all $j \notin J$, we can without loss of
generality assume that $J=\emptyset$ and $L_j=\emptyset$ for all $j.$
We need to remove another set of elements before we try to satisfy the
second criterion. 

Suppose that there exists some $n_1\in \N$ such that $J_1=\{j:
|\{l:(n_1,j,l)\notin
G\}|=\infty\}$ is infinite. Then we replace $\N\times \NxN$ by $\N\setminus
\{n_1\}\times \{(j,l):j\in J_1, (n_1,j,l)\notin G\}.$ If there exists $n_2
\in \N\setminus
\{n_1\}$ such that $J_2=\{j:j\in J_1,
|\{l:(n_i,j,l)\notin
G, i=1,2\}|=\infty\}$ we remove $n_2$ and consider only $j\in J_2$ and $l$
such that $(n_i,j,l)\notin
G, i=1,2.$ Because we have eliminated the pairs of coordinates which appear
$N$ or more times, this process must stop in at most $N-1$ steps. Thus
there exist some $k<N$, $N'=\N\setminus \{n_1,\dots,n_k\}$, an infinite
set $J'=J_k\subset \N$, and for each $j\in J'$, $L_j\subset \N$ infinite
such that
for each $l\in L_j$, $(n,j,l)\notin G$ for at most $N-1$ elements $n\in N'$
and for each $n\in N'$ there are only finitely many $j\in J'$ such that
$\{l: (n,j,l)\notin G\}$ is infinite.

We
now want to construct the sets $M$ and $K$ by a procedure like that used to
count the rationals. Choose $m_1 \in N'.$ Let $J'_1=\{j:
|\{l:(m_1,j,l)\notin G\}|<\infty\}$. Then by assumption $J'_1$ contains all
but finitely many elements of $J'.$ Let $L_{j,1}=\{l:(m_1,j,l)\in G\}$ for
each $j\in J'_1.$ Choose $j_1\in J'_1$ and then we can
find an infinite subset $L_1$ of
$L_{j_1,1}$ and a co-finite subset $N'_1$ of $N'$ such that for each $n\in
N'_1$ $(n,j_1,l)\notin G$ for at most finitely many $l\in L_{j_1,1}$ as follows.
If $L_1=L_{j_1,1}$ and $N_1'=N'$ work, we are done.
If not there is some $k_1\in N'$
such that $K_1=\{l:(k_1,j_1,l)\notin G\}$ is infinite. If $N'_1=N'\setminus
\{k_1\}$ and $L_1=K_1$ work we are done. If not there exists, $k_2\in
N_1'$, $k_2\neq k_1$, such that $K_2=\{l:(k_i,j_1,l)\notin G,i=1,2\}$ is
infinite. As in the earlier argument this can continue only at most $N-1$
times. Now choose $l_{j_1,1}\in L_1$. Let $N''_1=\{n:n\in N'_1,
(n,j_1,l_1)\in G\}.$ Again we have lost at most finitely many elements of
$N'_1.$

Next choose $m_2\in N_1''$, let $J'_2=\{j:j\in J'_1,
|\{l:(m_2,j,l)\notin G\}|<\infty\}$. Then by assumption $J'_2$ contains all
but finitely many elements of $J'_1.$
For $j\neq j_1$ let $L_{j,2}=\{l:(m_2,j,l)\in G\}$ for
each $j\in J'_2$  and $L_{j_1,2}=\{l:l\in L_1,(m_2,j,l)\in G\}$.
Note that $l_{j_1,1}\in L_{j_1,2}.$ Next choose $j_2\in J'_2, j_2\neq j_1.$
As above we can find a co-finite subset $N_2$ of $N''_1$ and find an
infinite subset $L_2$ of
$L_{j_2,2}$ such that for each $n\in
N'_2$, $(n,j_2,l)\notin G$ for at most finitely many $l\in L_{j_2,1}$. Now
choose $l_{j_1,2}\in L_{j_1,2}$, $l_{j_1,2}\neq l_{j_1,1}$
and $l_{j_2,1}\in L_2.$ Let
$N''_2=N'_2\setminus\{n:(n,j_1,l_{j_1,2})\notin G \text{ or
}(n,j_2,l_{j_2,1})\notin G\}.$ This removes at most a finite number of
elements from $N_2'$ but not $m_1$ or $m_2.$

The remainder of the argument consists of inductively 
choosing as we have done above new $m_i$'s, new
$j_i$'s, and corresponding $l_{j_i,k}$, so that in the end
$M=\{m_1,m_2,\dots\}$ and $K=\{(j_i,l_{j_i,k}):i,k\in \N\}$ satisfy the
second criterion. We leave the details to the reader.
\qed\enddemo

We are finally ready to prove our main result.

\proclaim{Theorem 8.3} Suppose that there is a constant $D$ such that for
all $n$ $Y_n$ is
a subspace of $L_p$ with a $D$-unconditional basis .
If $T$ is an isomorphism from $X_p\otimes X_p$ into
$(\sum Y_i)_{p,2,(w_n)}$ and
$P$ is a projection onto the range of $T$, then there exists an integer $N$
and a subspace $Z$ of $X_p \otimes X_p$, isomorphic to $X_p \otimes X_p$
such $P_N T|_Z$ is an isomorphism and $P_N T(Z)$ is complemented.
\endproclaim

\demo{Proof} We will use the standard basis of $X_{p,2,(w_{n,k})}$ where as
usual $w_{n,k}=w_n$ and $(w_n)$ decreases to 0. Let
$(x_{i,h})$ and $(y_{j,l})$ be two copies of the basis.
By Proposition 6.3 it is sufficient to find $N\in \N,\epsilon
>0,$ and rich subsets of
$\NxN$, $M,K$, such that
$$|x_{i,h}^*\otimes y_{j,l}^*(T^{-1}PP_NTx_{i,h}\otimes
y_{j,l})|>\epsilon$$
for all $(i,h)\in M,(j,l)\in K.$ In this proof we will take $\epsilon
=1/4.$ We proceed by induction as in the
proof of Proposition 6.3. Lemmas 7.2 and 7.4 will be used to show that a
uniform $N$ works. We will use the two player game approach with two
interwoven games as we did in the proof of Proposition 6.3. The proof is
essentially combinatorial so that topological condition in the definition
of the games will be 
irrelevant. (We could take $X_0$ to be the weak closure of the $X_p$ basis,
$(x_{m,n})$ or $(y_{j,l})$, and use constant functions,
but we will not define the functions at all.)
Lemmas 7.2 and 7.4 require that we use only small $w_k$
thus we immediately discard all indices $(i,h), (j,l)$ for which
$w_i$ or $w_j$ is larger than $(32
\|T\|\|T^{-1}\|\|P\|)^{-1}.$ This does not affect the isomorphic type of
the span of the remaining basis vectors. Therefore we will just assume that
the index sets are again $\NxN.$ Let $(M_m)$ be the sequence
of integers determined by Lemma 7.2 for $\rho=1/2,\rho_1=1/4$ and
$\delta=1/2.$

To reduce the notation a little define
$$f(N,i,h,j,l)=|x_{i,h}^*\otimes y_{j,l}^*(T^{-1}P(I-P_N)
Tx_{i,h}\otimes y_{j,l})|$$
and $$g(N,i,h,j,l)=|x_{i,h}^*\otimes y_{j,l}^*(T^{-1}PP_N
Tx_{i,h}\otimes y_{j,l})|$$
for all $i,h,j,l\in \N.$

Let $s_1$ be chosen according the game argument as in the
proof of Proposition 5.4 and let $S_{X,1}'$ be the set in $\Cal S_X$.
By Lemma 7.4 there is an integer $N_1$ and a
set
$S_{Y} \in \Cal S_Y$, such that
$$f(N_1,x_{\phi(s)},y_{\phi(t)})
<1/4$$
for all $t\in S_{Y}$ and
$s=s_1.$

Let
$$\multline
G=\{(h,j,l):\exists t\in S_{Y},s\in S_{X,1}'\backepsilon ,\phi(s)=(\phi(s_1)^1,h),\\
\phi(t)=(j,l),g(N_1,\phi(s),\phi(t))
\geq 1/4\}.
\endmultline$$
We apply Lemma 8.2 to $G$ with $N=M_{\phi(s_1)^1}$ and
with $\N\times S_{Y}$ in place of $\N\times\NxN$.

If alternative (2) occurs, let $H_1$ be the infinite subset and $S_{Y}'$ be
$\phi^{-1}(K)$, where $K$ is the rich subset. (We assume that $K$ is
maximal for $H_1$ and then that $H_1$  is maximal for this maximal $K$.) Let
$S_{X,1}''=S_{X,1}'\setminus \{(\phi(s_1)^1,h):h\notin H_1\}.$
Let $N_2=N_1.$ (This is a notational convenience.) Observe that we may
assume that $\phi(s_1)^2\in H_1$.

If alternative (1) occurs and not (2), let $K$ be the rich subset of
$\phi(S_{Y})$.
Then by Lemma 7.6 there is an integer $N_2$, an
infinite subset $H_1$ of $\{h:\exists s\in S_{X,1}' \backepsilon
\phi(s)=(\phi(s_1)^1,h)\}$ and a rich set
$K'\subset K$,
such that for any $N\geq N_2$, for all $j$, all  but at most $M_{\phi(s_1)}$
indices $h\in H_1$,
$f(N,\phi(s_1)^1,h,j,l)
<1/4$ for infinitely many $l$ with $(j,l)\in K'$. Moreover,
$f(N_2,\phi(s_1)^1,h,j,l)
<1/4$ for all $i\in H_1$, $(j,l)\in K'.$
Let $S_{Y}'=\phi^{-1}(K')$ and $S_{X,1}''=S_{X,1}'\setminus
\{(\phi(s_1)^1,h):h\notin H_1\}.$

Observe that in both cases we have that for the integer $N_2$,
$g(N_2,x_{\phi(s)},y_{\phi(t)})\geq 1/4$ for all $t\in S_{Y}',$
and all $s\in S_{X,1}''$ such that $\phi(s)^1=\phi(s_1)^1.$
This means that we can allow the second player in $X$-game to choose any
new element $s$ with $\phi(s)^1=\phi(s_1)^1$ from $S_{X,1}''$.
Also note that if alternative (2) failed, then for any $t\in S_Y'$, there
exist a set of $M_{\phi(s_1)}$ indices $H_1'\subset  \N$ such that
$$f(N_2,\phi(s_1)^1,h,\phi(t))<1/4$$ 
and 
$$g(N_1,\phi(s_1)^1,h,\phi(t))<1/4,$$
for all $h\in H_1'.$ Therefore
$$|x_{\phi(s_1)^1,h}^*\otimes y_{\phi(t)}^*(T^{-1}P(P_{N_2}-P_{N_1})
Tx_{\phi(s_1)^1,h}\otimes y_{\phi(t)})|>1/2,$$ for all $h\in H_1'$.
Thus we have one block for the $X$-game as in Lemma 7.2.

Now that we have started the construction, we can make further steps a
little more regular.

Let
$$T_0=\{t:t\in S_{Y}'\text{ and }
\phi(\{s:s\in S_{X,1}'',g(N_2,\phi(s),\phi(t))\geq 1/4\})
\text{ contains a rich set}\}.$$
If $\phi(T_0)$ contains a rich set, let $K$ be a maximal one and
let $S_{Y,1}=\phi^{-1}(K).$ If not, then $S_{Y}'\setminus T_0$ contains a
rich set, $K$. Let $S_{Y,1}=\phi^{-1}(K).$

Player 2 in the $Y$-game chooses $t_1\in S_{Y,1}$ and $S_{Y,1}'\in \Cal S_Y$
such that $S_{Y,1}'\subset S_{Y,1}.$

Let
$$\multline
G=\{(l,i,h):\exists t\in S_{Y,1}',s\in S_{X,1}'' ,\phi(t)=(\phi(t_1)^1,l),\\
\phi(s)=(i,h),g(N_2,\phi(s),\phi(t))
\geq 1/4\}.
\endmultline$$
We apply Lemma 8.2 to $G$ with $N=M_{\phi(t_1)^1}$ and
with $\N\times S_{X,1}''$ in place of $\N\times\NxN$.

If alternative (2) occurs, let $L_1$ be the infinite subset and 
$$S_{X,1}''' =\{s\in S_{X,1}'':\phi(s)^1=\phi(s_1)^1\}\cup\phi^{-1}(K),$$
where $K$ is the (maximal) rich subset. (We know that if $s\in S_{X,1}''$ and
$\phi(s)^1=\phi(s_1)^1$, then $g(N_2,\phi(s),\phi(t))
\geq 1/4$ and thus $(l,\phi(s))\in G.$)
Let
$S_{Y,1}''=S_{Y,1}'\setminus \{(\phi(t_1)^1,l):l\notin L_1\}.$
Let $N_3=N_2.$ (This is again a notational convenience.) Observe that we may
assume that $\phi(t_1)^2\in L_1$ by making $L_1$ maximal.

If alternative (1) occurs and not (2), let $K$ be the rich subset of
$\phi(S_{X,1}'')$.
Then by Lemma 7.6 there is an integer $N_3>N_2$, an
infinite subset $L_1$ of $\N$ and a rich set                  
$K'\subset K$,
such that for any $N\geq N_3$, for all $i\in \phi(S_{X,1}'')^1$,
all  but at most $M_{\phi(t_1)}$
indices $l\in L_1$,
$f(N,i,h,\phi(t_1)^1,l)
<1/4$ for infinitely many $h$ with $(i,h)\in K'$. Moreover,
$f(N_3,i,h,\phi(t_1)^1,l)
<1/4$ for all $l\in L_1$, $(i,h)\in K'.$
Let 
$$S_{X,1}'''=
\phi^{-1}(K')$$
and 
$$S_{Y,1}''=S_{Y,1}'\setminus
\{(\phi(t_1)^1,l):l\notin L_1\}.$$

Observe that in both cases we have that for the integer $N_3$,
$g(N_3,\phi(s),\phi(t))\geq 1/4$ for all $s\in S_{X,1}''',$
and all $t\in S_{Y,1}''$ such that $\phi(t)^1=\phi(t_1)^1.$
This means that we can allow the second player in $Y$-game to choose any
new element $t$ with $\phi(t)^1=\phi(t_1)^1$ from $S_{Y,1}''$.
Also note that if alternative (2) failed, then for any $s\in S_{X,1}'''$, there
exist a set of $M_{\phi(t_1)}$ indices $L_1'\subset  \N$ such that
$$f(N_3,\phi(s),\phi(t_1)^1,l)<1/4$$
and
$$g(N_2,\phi(s),\phi(t_1)^1,l)<1/4,$$
for all $l\in L_1'.$ Therefore
$$|x_{\phi(s)}^*\otimes y_{\phi(t_1)^1,l}^*(T^{-1}P(P_{N_3}-P_{N_2})
Tx_{\phi(s)}\otimes y_{\phi(t_1)^1,l})|>1/2,$$ for all $l\in L_1'$.
Thus we have one block for the $Y$-game as in Lemma 7.2.

The idea of the proof is to continue to play the two games as above. Each
time we are forced to take alternative (1) in the application of Lemma 8.2,
we produce a new block as in Lemma 7.2. However Lemma 7.2 tells us that this
cannot go on happening. Thus eventually only alternative (2) occurs in each
game and we are free to construct the required basic sequences.
In order to play the games according to the rules we will fatten the sets
$S_{X,r},S_{Y,r}$, so that they will contain the previously chosen
elements. However in the end we will discard the elements chosen before the
integer $N_k$ has been fixed.

We will do a few more steps of the induction in order to cover a few cases
which have yet to arise.

Let
$$\multline
S_0=\{s:s\in S_{X,1}'''\text{ and }\\
\phi(\{t:t\in S_{Y,1}'',f(N_3,\phi(s),\phi(t))<3/4\})
\text{ contains a rich set}\}.
\endmultline$$
If $\phi(S_0)$ contains a rich set, let $K$ be a maximal one and 
let $S_{X,2}=\{s_1\}\cup\{s:s\in S_{X,1}'',\phi(s)^1=\phi(s_1)^1\}\cup
\phi^{-1}(K).$ If not, then $S_{X,1}'''\setminus S_0$ contains a
rich set, $K$. Let $S_{X,2}=\{s_1\}\cup\{s:s\in
S_{X,1}'',\phi(s)^1=\phi(s_1)^1\}\cup \phi^{-1}(K).$

Player 2 in the $X$-game chooses $s_2\in S_{X,2}$ and $S_{X,2}'\in \Cal
S_X$
such that $S_{X,2}'\subset S_{X,2}.$

If it happens that
$\phi(s_1)^1=\phi(s_2)^1$ and $N_3=N_2$, there is nothing to do. Thus we
simply let $S_{Y,2}'''=S_{Y,1}''.$ Otherwise,
let
$$\multline
G=\{(h,j,l):\exists t\in S_{Y,1}'',s\in S_{X,2}' ,\phi(s)=(\phi(s_2)^1,h),\\
\phi(t)=(j,l),f(N_3,\phi(s),\phi(t))
<3/4\}.
\endmultline$$
We apply Lemma 8.2 to $G$ with $N=M_{\phi(s_2)^1}$ and
with $(\phi(S_{X,2}')\cap(\{\phi(s_2)^1\}\times\N))\times S_{Y,1}''$
in place of $\N\times\NxN$.

If alternative (2) occurs, let $H_2$ be the infinite subset and $S_{Y,1}'''$ be
$\phi^{-1}(K)$, where $K$ is the maximal rich subset. Let
$S_{X,2}''=S_{X,2}'\setminus \{(\phi(s_2)^1,h):h\notin H_2\}.$
Let $N_4=N_3.$ Observe that we may
assume that $\phi(s_2)^2\in H_2$.

If alternative (1) occurs and not (2), let $K$ be the rich subset of
$\phi(S_{Y,1}'')$.
Then by Lemma 7.6 there is an integer $N_4>N_3$, an
infinite subset $H_2$ of $\N$ and a rich set
$K'\subset K$,
such that for any $N\geq N_4$, for all $j$, all  but at most
$M_{\phi(s_2)}$
indices $h\in H_1$,
$f(N,\phi(s_2)^1,h,j,l)
<1/4$ for infinitely many $l$ with $(j,l)\in K'$. Moreover,
$f(N_4,\phi(s_2)^1,h,j,l)
<1/4$ for all $h\in H_2$, $(j,l)\in K'.$
Let $S_{Y,1}'''=\phi^{-1}(K')$ and $S_{X,2}''=S_{X,2}'\setminus
\{(\phi(s_2)^1,h):h\notin H_2\}.$

Observe that in all three cases we have that for the integer $N_4$,
$$g(N_4,x_{\phi(s)},y_{\phi(t)})\geq 1/4$$
 for all $t\in S_{Y,1}''',$
and all $s\in S_{X,2}''$ such that $\phi(s)^1=\phi(s_2)^1.$
This means that we can allow the second player in $X$-game to choose any
new element $s$ with $\phi(s)^1=\phi(s_r)^1$, $r=2$ (or $r=1,2$ if
$N_4=N_3$)
from $S_{X,2}''$.
Also note that if alternative (2) failed, then for any $t\in S_Y'$, there
exist a set of $M_{\phi(s_2)}$ indices $H_2'\subset  \N$ such that
$$f(N_4,\phi(s_2)^1,h,\phi(t))<1/4$$
and
$$g(N_3,\phi(s_2)^1,h,\phi(t))<1/4,$$
for all $h\in H_2'.$ Therefore
$$|x_{\phi(s_2)^1,h}^*\otimes y_{\phi(t)}^*(T^{-1}P(P_{N_4}-P_{N_3})
Tx_{\phi(s_2)^1,h}\otimes y_{\phi(t)})|>1/2,$$ for all $h\in H_2'$.
Thus for any $t\in S_{Y,1}'''\subset S_Y'$,
we have one or two blocks (depending on whether alternative (2) has
failed two times in $X$-game) for the $X$-game as in Lemma 7.2.

Let
$$\multline
T_0=\{t:t\in S_{Y,1}'''\text{ and }\\
\phi(\{s:s\in S_{X,2}'',g(N_4,\phi(s),\phi(t))\geq 1/4\})
\text{ contains a rich set}\}.
\endmultline$$
If $\phi(T_0)$ contains a rich set, let $K$ be a maximal rich subse and
let $S_{Y,2}=\phi^{-1}(K)\cup\{t_1\}\cup\{t\in
S_{Y,1}''':\phi(t)^1=\phi(t_1)^1\}$.
If not, then $S_{Y,1}'''\setminus T_0$ contains a
rich set, $K$. Let $S_{Y,2}=\phi^{-1}(K)\cup\{t_1\}\cup\{t\in
S_{Y,1}''':\phi(t)^1=\phi(t_1)^1\}$.

Player 2 in the $Y$-game chooses $t_2\in S_{Y,2}$ and $S_{Y,2}'\in \Cal
S_Y$
such that $S_{Y,2}'\subset S_{Y,2}.$

If $\phi(s_2)^1=\phi(s_1)^1$ and $N_3=N_4$, then let $S_{X,2}'''=S_{X,2}''.$
otherwise let
$$\multline
G=\{(i,h,l):\exists t\in S_{Y,2}',s\in S_{X,2}''
,\phi(t)=(\phi(t_2)^1,l),\\
\phi(s)=(i,h),g(N_4,\phi(s),\phi(t))
\geq 1/4\}.
\endmultline$$
We apply Lemma 8.2 to $G$ with $N=M_{\phi(t_2)^1}$ and
with $(\phi(S_{Y,2}')\cap(\{\phi(t_2)^1\}\times\N))
\times S_{X,2}''$ in place of $\N\times\NxN$.

If alternative (2) occurs, let $L_2$ be the infinite subset and $S_{X,2}'''$
be
$\phi^{-1}(K)$, where $K$ is the maximal rich subset. Let
$S_{Y,2}''=S_{Y,2}'\setminus \{(\phi(t_2)^1,l):l\notin L_2\}.$
Let $N_5=N_4.$ (This is again a notational convenience.) Observe that we
may
assume that $\phi(t_2)^2\in L_2$ by making $L_2$ maximal.

If alternative (1) occurs and not (2), let $K$ be the rich subset of
$\phi(S_{X,2}'')$.
Then by Lemma 7.6 there is an integer $N_5$, an
infinite subset $L_2$ of $\N$ and a rich set
$K'\subset K$,
such that for any $N\geq N_5$, for all $i\in \phi(S_{X,2}'')^1$,
all  but at most $M_{\phi(t_2)}$
indices $l\in L_2$,
$f(N,i,h,\phi(t_2)^1,l)
<1/4$ for infinitely many $h$ with $(i,h)\in K'$. Moreover,
$f(N_5,i,h,\phi(t_2)^1,l)
<1/4$ for all $l\in L_2$, $(i,h)\in K'.$
Let $S_{X,2}'''=\phi^{-1}(K')$ and $S_{Y,2}''=S_{Y,2}'\setminus
\{(\phi(t_2)^1,l):l\notin L_2\}.$

Observe that in all three cases we have that for the integer $N_5$,
$g(N_5,\phi(s),\phi(t))\geq 1/4$ for all $s\in S_{X,2}''',$
and all $t\in S_{Y,2}''$ such that $\phi(t)^1=\phi(t_2)^1.$
This means that we can allow the second player in $Y$-game to choose any
new element $t$ with $\phi(t)^1=\phi(t_2)^1$ from $S_{Y,2}''$.
Also note that if alternative (2) failed, then for any $s\in S_{X,2}'''$,
there
exist a set of $M_{\phi(t_2)}$ indices $L_2'\subset  \N$ such that
$$f(N_5,\phi(s),\phi(t_2)^1,l)<1/4$$
and
$$g(N_4,\phi(s),\phi(t_2)^1,l)<1/4,$$
for all $l\in L_1'.$ Therefore
$$|x_{\phi(s)}^*\otimes y_{\phi(t_2)^1,l}^*(T^{-1}P(P_{N_5}-P_{N_4})
Tx_{\phi(s)}\otimes y_{\phi(t_2)^1,l})|>1/2,$$ for all $l\in L_2'$.
Thus for any $s\in S_{X,2}'''\subset S_{X,1}'''$,
we have one or two blocks (depending on whether alternative (2) has
failed two times in $Y$-game) for the $Y$-game as in Lemma 7.2.

We have now established the pattern of the induction. As we have noted
earlier alternative (2) of Lemma 8.2 can not fail infinitely many times,
thus for some $k_0$, $N_k=N_{k_0}$ for all $k\geq k_0.$ Also
$f(N_{k_0},\phi(s_k),\phi(t_{k'}))<3/4$ for all $k,k'\geq k_0.$
Therefore $g(N_{k_0},\phi(s_k),\phi(t_{k'}))>1-3/4=1/4$ for all $k,k'\geq
k_0.$ Because the games must yield a sequence in $\Cal S_X$, $\Cal S_Y$,
respectively, and removing finitely many rows from such a set is still
such a set, $(x_{\phi(s_k)})_{k\in S}$ and $(y_{\phi(t_k)})_{k
\in T}$, where $S=\{s_r:\phi(s_r)^1\neq \phi(s_k)^1 \qquad\forall k<k_0\}$ and
$T=\{t_r:\phi(t_r)^1\neq \phi(t_k)^1 \qquad\forall k<k_0\}$,
are equivalent to bases of $X_p$ and Proposition 6.3 concludes the proof.
\qed\enddemo

\proclaim{Corollary 8.4} For all $\alpha < \omega_1,$
$X_p\otimes X_p$ is not isomorphic to a
complemented subspace of $R^\alpha_p.$
\endproclaim

\demo{Proof} We know that  for each $\alpha <\omega_1$,
$R^\alpha_p$ is isomorphic to  an $(p,2)$ sum of spaces $R_p^{\alpha_n}$
where for each $n$, $R_p^{\alpha_n}$ is not isomorphic to $R^\alpha_p$.
If $X_p\otimes X_p$ were isomorphic to a complemented subspace of some
$R_p^\gamma$, $\gamma <\omega_1,$ let $\alpha$ be the  smallest such
ordinal. Since $R^\alpha_p$ is isomorphic to $(\sum R_p^{\alpha_n})_{p,2},$
Theorem 8.3 implies that $X_p\otimes X_p$ is isomorphic to a complemented
subspace of
$(\sum_{n=1}^N R_p^{\alpha_n})_{p,2},$ for some $N\in \N.$ However this
space is isomorphic to $R_p^{\alpha'}$ for some $\alpha'$ such that
$\alpha'+\omega\leq \alpha.$ This contradicts the choice of $\alpha.$
\qed\enddemo

\newpage

\head 9. Final remarks and open problems \endhead

In this paper we have answered some of the questions posed in \cite{BRS}, but
there are many more questions that are raised by this work. 

\item{1.} In Section 1 we note that the best projection onto
$\otimes^k X_p$ has
norm which grows with $k$. Thus it is natural to ask:
 Is $\otimes^k X_p$ isomorphic to a well complemented subspace of
$L_p$? More precisely, is there a constant $C$ and subspaces $Y_k, k\in \N$
of $L_p$ such that $\otimes^k X_p$ is isomorphic to $Y_k$ and there is a
projection $P_k$ of $L_p$ onto $Y_k$ with $\|P_k\|\leq C$?

\item{2.} In \cite{S} Schechtman uses spaces of the form 
$$(\sum (\sum \dots(\sum
\ell_{r_1})_{r_2}\dots)_{r_{n-1}})_{r_n}$$
 in order to distinguish the
spaces $\otimes^k X_q$, $k\in \N.$ If $2\geq r_1>r_2>\dots>r_1>q$, is
$(\sum (\sum \dots(\sum \ell_{r_1})_{r_2}\dots)_{r_{n-1}})_{r_n}$
isomorphic to a subspace of $\Rq{\omega n}$? For $n=1$ this is a result of
Rosenthal \cite{R2} and it is not hard to see
that for $n=2$ and $r_1=2$ or $r_2=q$
 that it is
true. However we do not know whether there is any $\alpha<\omega_1,$ such
that for $2>r_1>r_2>q$, $(\sum
\ell_{r_1})_{r_2}$ is isomorphic to a subspace of $\Rq{\alpha}$. Notice that
were it the case that there is no such $\alpha$, then Corollary 8.4 would
follow.

\item{3.} The proof of Theorem 8.3 that we have presented uses the
unconditional basis assumption very sparingly. We had hoped to eliminate it
altogether. Is the assumption that each $Y_n$ have a $D$-unconditional
basis necessary? Can the proof of Theorem 8.3 be simplified substantially
if we make more use of the assumption that the spaces $Y_n$ have
unconditional bases?

\item{4.} In Section 5. we introduce a general framework for gliding hump
type arguments, but we do not carry the work very far. What are good
classes $\Cal S$ for the natural bases of 
spaces such as $\otimes^k X_p$, $(\sum (\sum
\dots(\sum\ell_{r_1})_{r_2}\dots)_{r_{n-1}})_{r_n}$, spaces modeled on
trees, etc.?

\item{5.} This paper shows that at least for certain questions the spaces
$R_p^\alpha$ are similar enough to $X_p$ that the techniques originated by
Rosenthal can be adapted for use with these spaces. J. T. Woo \cite{Wo1}, 
\cite{Wo2},
 showed that
$X_p$ is just one of a collection of modular sequence spaces with similar
properties. Many of the arguments in this paper are really about multiple
norm spaces. Thus it is likely that much of it would generalize to a class
spaces where $p$ and $2$ are replaced by $p$ and $r$ or perhaps by spaces
which are defined by families  of indices.

\item{6.}
Suppose that $P$ is a projection on $X_p\otimes
X_p$, is there a complemented subspace $Z$ of $X_p\otimes X_p$ which
is isomorphic to $X_p\otimes X_p$ and is contained in the the range of $P$
or the range of $I-P.$ Because of Propostions 6.3 and 6.8,
we think that it is very likely that this is true. The main difficulty
remaining seems to be a combinatorial problem. Suppose that $G \subset
\N^4$ and $\phi$ is a bijection from $\N$ onto $\NxN$ as in Section 5. Are
there infinite subsets $K,L$ of $\N$, such that
$$\{(\phi(k),\phi(l)):o(k,K)>o(l,L)\}$$ or
$$\{(\phi(k),\phi(l)):o(k,K)<o(l,L)\}$$ is contained in $G$ or
$\N^4\setminus G$ and $\phi(K)$ and $\phi(L)$ are rich? If these questions
have affirmative answers then it may be possible to show that
$X_p\otimes X_p$ is primary.

\enddocument

\bye